\documentclass{article}
\usepackage{graphicx} 
\usepackage{amssymb}
\usepackage{amsfonts}
\usepackage{url}
\usepackage{amsmath,latexsym}
\usepackage{bm}

\usepackage[linecolor=red!40!white, backgroundcolor=blue!05!white,bordercolor=blue!40!white,textsize=small]{todonotes}

\tikzset{notestyleraw/.append style={align=justify}}

\newtheorem{theorem}{Theorem}

\newtheorem{corollary}{Corollary}

\newtheorem{definition}{Definition}

\newtheorem{lemma}{Lemma}

\newtheorem{proposition}{Proposition}
\newtheorem{remark}{Remark}

\newenvironment{proof}[1][Proof]{\textbf{#1.} }{\ \rule{0.5em}{0.5em}}

\title{Two-temperature fluid models for a polyatomic gas based on kinetic theory for
nearly resonant collisions}
\author{Kazuo Aoki \thanks{Department of Mathematics,
National Cheng Kung University, Tainan 70101, Taiwan}
\;\,and\, Niclas Bernhoff \thanks{Department of Mathematics
and Computer Science, Karlstad University, 65188 Karlstad, Sweden}}
\date{\today}

\begin{document}

\maketitle

\abstract
A polyatomic ideal gas with weak interaction between the
translational and internal modes is considered. For the purpose of describing the
behavior of such a gas, a Boltzmann equation is proposed in the form that
the collision integral is a linear combination of inelastic and elastic (or resonant)
collisions, and its basic properties are discussed. Then, in the case where the
elastic collisions are dominant, fluid dynamic equations of Euler and Navier--Stokes
type including two temperatures, i.e., translational and internal temperatures,
as well as relaxation terms are systematically obtained by means of the
Chapman--Enskog expansion. The obtained equations are different depending
on the degree of weakness of the interaction between the translational and
internal modes.

\section{Introduction}\label{sec:intro}

Multi-temperature fluid models have been widely used for high-speed and
high-temperature flows of polyatomic gases \cite{P-90,NK-09}. Because these flows are generally in highly nonequilibrium, fluid models must be based on kinetic theory.
However, it is not an easy task to derive multi-temperature fluid models
systematically from kinetic theory. This is mainly because 
the Boltzmann equation for polyatomic gases is very complex due to
energy exchange between translational and internal modes and
between different internal modes during molecular collisions \cite{MBKK-90,CC-91,G-98,C-00,NK-09}.
The Boltzmann equation based on state-to-state models
\cite{NK-09,AK-13,KKKN-21} can accurately describe all these exchanges
using different collision integrals and can, in principle, provide multi-temperature
fluid models. However, these models are effective for specific gases for which
the data for the transition probability are available and thus
lack generality with respect to gas species. Therefore, various simplified kinetic
models have been proposed so far. Throughout this paper, we do not consider
gas mixtures and restrict ourselves to
a single polyatomic (including diatomic) ideal gas.

One of such simplified models is to use model kinetic equations of relaxation type, such as
the Bhatnagar--Gross--Krook (BGK) model and the Ellipsoidal Statistical (ES) model,
instead of the Boltzmann equation. In fact, various model equations of this type
have been proposed
(e.g., \cite{M-64, H-66,GS-99,S-99,ALPP-00,BrS-09,BS-17,KKA-19,MM-20,DMM-21,AMR-24,BrP-25})
and successfully used in many applications. It should
also be mentioned that rigorous mathematical studies of these models have
also been conducted (e.g., \cite{Y-19,PY-19}). 

Another approach is to keep the Boltzmann equation as is but to use simplified
models for transition probabilities in the collision integrals. The simplest model
introduces an additional variable, which is either discrete or continuous,
representing the total energy of the internal modes. The approach using the
continuous variable was first introduced for the purpose of numerically simulating
collision processes \cite{BL-75}, but later the corresponding collision integral
was constructed explicitly \cite{BDLP-94}. This motivated recent mathematical studies
of the Boltzmann equation with a single additional variable, and some important
results have been obtained (e.g.,
\cite{DPS-21,BBG-22,GP-23,Be-23a,Be-23b,DOPT-23,BBS-23,BRS-24,BST-24,BBPG-24,BBMS-25}).

In this paper, attention is focused on the formal, but systematic, derivation of two-temperature
fluid models from kinetic theory. In spite of the importance of the topic, the number
of published papers on it has been limited \cite{W-CU-51,BG-11,ABGK-20,BG-22},
mainly due to the complexity of kinetic models
for polyatomic gases mentioned above. In \cite{ABGK-20},  a two-temperature fluid model of Navier--Stokes
type was systematically derived from the kinetic ES model \cite{ALPP-00}. Furthermore,
the boundary conditions for the two-temperature Navier--Stokes model have
been established using Knudsen-layer analysis \cite{KABGM-21}. The advantage of
starting from the ES model is that the resulting two-temperature Navier--Stokes
model has explicit parameter dependence, so that its application to practical
problems is easy. On the other hand, the derivation of multi-temperature fluid
models from Boltzmann-type models, rather than relaxation-type kinetic models
such as BGK or ES models, is a fundamental problem that has only been partially
resolved (e.g., \cite{W-CU-51,BG-11,BG-22}). 

The aim of the present paper is to establish two-temperature fluid models on
the basis of a Boltzmann-type kinetic equation, rather than the models of relaxation type.
We employ the Boltzmann
equation with an additional continuous variable corresponding to the total
energy of the internal modes. We first propose a model of the collision kernel
that is a linear combination of a standard (or inelastic) collision kernel with
coefficient $\theta$ and a resonant (or elastic) collision kernel with coefficient
$(1-\theta)$. Here, resonant collisions are collisions in which there is no energy
exchange between the translational and internal modes \cite{BBS-23,BRS-24,BBMS-25},
and $\theta$
is a parameter indicating perfectly resonant collisions when $\theta=0$. Then,
assuming that the Knudsen number Kn (the ratio of the mean free path of
the gas molecules to the characteristic length) is small, we consider the case
when the resonant collisions are dominant, that is, when $\theta$ is small.
This corresponds to a polyatomic gas in which the interaction between the
translational and internal modes
is weak; in other words, the relaxation of internal modes is slow.
We derive fluid equations of Euler and Navier--Stokes types, which include
translational and internal temperatures as well as relaxation terms,
by means of the Chapman--Enskog expansion \cite{CC-91,S-07}
for two different cases: (i) $\theta$ is of the order of $\text{Kn}^2$; and (ii) $\theta$
is of the order of $\text{Kn}$.

It should be remarked here that various higher-order macroscopic equations with
two (or multi) temperatures, different from Euler or Navier--Stokes type models,
have been constructed (e.g., \cite{ATRS-12,PRS-13,RS14,TARS-16,ARS-18,RS-21}).
Some of them are based on extended or irreversible
thermodynamics, where information from kinetic theory is partially taken
into account, and others are based on moment equations
derived directly from the Boltzmann equation. In any approach, one needs
appropriate closure assumptions, which characterize the resulting macroscopic equations. 

The paper is organized as follows. In Sec.~\ref{sec:kinetic-model}, the Boltzmann model
used here is presented and its basic properties are summarized. For example,
the equilibrium solution (Sec.~\ref{sec:col-inv}), the corresponding linearized collision
operator (Sec.~\ref{sec:linear-col-op}), and its Fredholm properties (Sec.~\ref{sec:Fredholm})
are discussed. In particular, a specific collision kernel, which is the basis of the
subsequent analysis, is introduced in Sec.~\ref{sec:Fredholm}. Section \ref{sec:nearly-resonant}
is devoted to the derivation of two-temperature fluid models.
In Sec.~\ref{sec:prelim}, necessary preliminaries are given.
Then, the case of $\theta = O(\text{Kn}^2)$ and that of
$\theta=O(\text{Kn})$ are studied in detail in Secs.~\ref{sec:CE-1} and
\ref{sec:CE-2}, respectively. Finally, concluding remarks are given in Sec.~\ref{sec:concluding}.

\section{Kinetic model}\label{sec:kinetic-model}

In this section, the kinetic model that will be considered in this paper
is explained.

\subsection{Velocity-energy distribution function and
macroscopic quantities}

Let us consider an ideal polyatomic (or diatomic) rarefied gas.
Let $t \in \mathbb{R}_+$ be the time variable, $\bm{x}$ (or $x_i$) $\in \mathbb{R}^3$
the position vector in the physical space, $\bm{\xi}$ (or $\xi_i$) $\in \mathbb{R}^3$
the molecular velocity, and $I \in \mathbb{R}_+$ the total energy associated
with the internal modes per molecule. We denote by
\begin{align*}
f(t,\, \bm{x},\, \bm{\xi},\, I)\, \text{d}\bm{x}\, \text{d}\bm{\xi}\, \text{d}I,
\end{align*}
the total number of gas molecules contained in an infinitesimal volume
$\text{d}\bm{x}\, \text{d}\bm{\xi}\, \text{d}I$
around a point $(\bm{x},\, \bm{\xi},\, I)$ in the seven-dimensional
space, consisting of $\bm{x}$, $\bm{\xi}$, and $I$, at time $t$.
We may call $f(t,\, \bm{x},\, \bm{\xi},\, I)$, which is the number density
in the seven-dimensional space, the velocity-energy
distribution function of the gas molecules. 

Let $\delta\, ( \ge 2)$ be the number of internal degrees of freedom, which
is constant but not necessarily an integer. Under the assumption that
the equipartition law holds, the ratio of the specific heats $\gamma$ is
expressed as
\begin{align*}
\gamma = \frac{c_\text{p}}{c_\text{v}} = \frac{\delta + 5}{\delta + 3},
\end{align*}
where $c_\text{p}$ is the specific heat at constant pressure and
$c_\text{v}$ is that at constant volume.

To define macroscopic quantities of the gas, we introduce the real Hilbert space
$L^{2} \left( \text{d}\bm{\xi}\, \text{d} I \right)$, with inner product
\begin{equation*}
( f, g ) =\int_{\mathbb{R}^{3}\times \mathbb{R}_{+}} fg\, 
\text{d}\bm{\xi} \text{d}I 
\quad \text{ for }\;\; f, g \in L^{2} ( \text{d}\bm{\xi}
\text{d}I ).
\end{equation*}
Let us denote by $m$ the mass of a molecule and by $k_\text{B}$ the
Boltzmann constant.
Let $n$ be the molecular number density, $\rho$ the mass density,
$\bm{u}$ (or $u_i$) the flow velocity, $e$ the internal energy per molecule,
$T$ the temperature, $p$ the pressure, $p_{ij}$ the stress tensor, and
$\bm{q}$ (or $q_i$) the heat-flow vector. Then, they are defined by
\begin{align}\label{macro}
\begin{aligned}
& n = (1,\, f), \qquad \rho = mn = (m,\, f), \qquad u_i = \frac{1}{n} (\xi_i,\, f),
\\
& e = e_\text{tr} + e_\text{int}, \qquad
e_\text{tr} = \frac{1}{n} \left( \frac{m}{2}|\bm{\xi}-\bm{u}|^2,\, f \right), \qquad
e_\text{int} = \frac{1}{n} (I,\, f),
\\
& T = \frac{3T_\text{tr} + \delta T_\text{int}}{3 + \delta}, \qquad
T_\text{tr} = \frac{2}{3 k_\text{B}} e_\text{tr}, \qquad
T_\text{int} = \frac{2}{\delta k_\text{B}} e_\text{int},
\\
& p_{ij} = \big( m(\xi_i - u_i) (\xi_j - u_j),\, f \big),
\\
&q_i = q_{(\text{tr}) i} + q_{(\text{int}) i},
\qquad q_{(\text{tr}) i}  = \left( (\xi_i - u_i) \frac{m}{2} |\bm{\xi} - \bm{u}|^2,\, f \right),
\\
& q_{(\text{int}) i} = \big( (\xi_i - u_i) I,\, f \big),
\end{aligned}
\end{align}
where $e_\text{tr}$ and $e_\text{int}$ are, respectively, the contribution
of the translational motion and that of the internal modes to
the internal energy $e$ per molecule, and $T_\text{tr}$ and $T_\text{int}$
are, respectively, the temperature associated with the translational
motion and that associated with the internal modes. We will call $T_\text{tr}$
the translational temperature and $T_\text{int}$ the internal
temperature.

\subsection{Boltzmann equation and collision operator}\label{sec:B-EQ}

\subsubsection{Boltzmann equation}

The evolution of the velocity-energy distribution function is,
in the absence of external forces, described by the Boltzmann equation of the form
\begin{equation} \label{BE-1}
\frac{\partial f}{\partial t} + \bm{\xi } \bm{\cdot} \frac{\partial f}{\partial \bm{x}}
=Q_{\theta }\left( f,\, f \right),
\end{equation}
where the collision operator $Q_{\theta }=Q_{\theta }\left( f,\, f \right) $ is
a quadratic bilinear operator that accounts for the change of velocities and of energy of the internal modes of the molecules due to binary collisions (assuming that the
gas is rarefied, so that other collisions are negligible). The collision
operator $Q_\theta$ will be detailed in the following.

\subsubsection{Binary collisions}

A collision can be represented by two pre-collisional pairs, each pair
consisting of a molecular velocity and an energy of the internal modes,
$\left( \bm{\xi},\, I \right) $ and $\left( \bm{\xi}_{\ast },\, I_{\ast} \right)$,
and two corresponding post-collisional pairs,
$\left( \bm{\xi }^{\prime },\, I^{\prime } \right)$ and
$\left( \bm{\xi}_{\ast }^{\prime },\, I_{\ast }^{\prime } \right) $.
The notation for pre- and
post-collisional pairs may, of course, be interchanged as well. Due to
momentum and total energy conservation, the following relations have to be
satisfied by the pairs:
\begin{eqnarray}\label{conservation}
\begin{aligned}
\bm{\xi }+\bm{\xi }_{\ast }&  = \bm{\xi }^{\prime }+%
\bm{\xi }_{\ast }^{\prime }
\\
\frac{m}{2} \left\vert \bm{\xi } \right\vert ^{2}+\frac{m}{2}%
\left\vert \bm{\xi }_{\ast } \right\vert ^{2}+I+I_{\ast }
& = \frac{m}{2} \left\vert \bm{\xi }^{\prime } \right\vert ^{2}+\frac{m}{2}%
\left\vert \bm{\xi }_{\ast }^{\prime } \right\vert ^{2}+I^{\prime
}+I_{\ast }^{\prime }.
\end{aligned}
\end{eqnarray}
The momentum conservation can be expressed as the conservation of
the velocity of the center of mass, i.e.,
\begin{align*}
\bm{G} = \bm{G}', \qquad \bm{G} := \frac{\bm{\xi} + \bm{\xi}_*}{2}, \qquad
\bm{G}': = \frac{\bm{\xi}' + \bm{\xi}_*'}{2},
\end{align*}
and the energy conservation can also be expressed through the conservation of the
total energy in the center of mass frame, i.e.,
\begin{equation}\label{EC-2}
E = E', \quad
E:=\frac{m}{4}\left\vert \bm{g} \right\vert ^{2}+I+I_{\ast }, \quad
E' :=\frac{m}{4}\left\vert \bm{g}^\prime \right\vert ^{2}+I^{\prime
}+I_{\ast }^{\prime },
\end{equation}
where the relative velocities before and after the collision
are introduced: 
\begin{equation*}
\bm{g}:=\bm{\xi }-\bm{\xi }_{\ast } \quad \text{ and } \quad \bm{g}%
^{\prime }:=\bm{\xi }^{\prime }-\bm{\xi }_{\ast }^{\prime }.
\end{equation*}
Incidentally, the gap of the energy of the internal modes for the collision is denoted by
\begin{equation*}
\Delta I:=I^{\prime }+I_{\ast }^{\prime }-I-I_{\ast }.
\end{equation*}

The collision during which there is no energy exchange between the
translational mode and the internal modes is called a resonant
(or elastic) collision.
Therefore, the energy conservation holds for the translational and
internal modes separately. That is, the following relations hold: 
\begin{equation}\label{conservation-theta=0}
\frac{m}{2} \left\vert \bm{\xi } \right\vert ^{2}+\frac{m}{2}%
\left\vert \bm{\xi }_{\ast } \right\vert ^{2}=\frac{m}{2}\left\vert 
\bm{\xi }^{\prime } \right\vert ^{2}+\frac{m}{2}\left\vert 
\bm{\xi }_{\ast }^{\prime } \right\vert ^{2} \quad
\text{ and } \quad
I+I_{\ast}=I^{\prime }+I_{\ast }^{\prime },
\end{equation}
or
\begin{equation*}
\left\vert \bm{g} \right\vert = \left\vert \bm{\xi }-\bm{%
\xi }_{\ast } \right\vert = \left\vert \bm{\xi }^{\prime }-\bm{%
\xi }_{\ast }^{\prime } \right\vert =\left\vert \bm{g}^{\prime
}\right\vert \quad \text{ and } \quad \Delta I=0. 
\end{equation*}
Resonant collisions will play an important role in the present paper.

\subsubsection{Collision operator}

Let $F$ be a function of $t$, $\bm{x}$, $\bm{\xi}$, and $I$. The model of the
collision operator $Q_\theta (f, f)$ in the Boltzmann equation \eqref{BE-1}
that is adopted here is defined via the following
bilinear operator (cf.~\cite{Be-23b}):
\begin{align}\label{col-op}
Q_{\theta }(f, F) &= \frac{1}{2} \int_{\left( \mathbb{R}^{3}\times \mathbb{R}_{+}\right)^{3}}
\left( \frac{ f' F_*'+f_* 'F' }{\left( I^{\prime}I_{\ast }^{\prime }\right)^{\delta /2-1}}
-\frac{ f F_* + f_* F }{\left( II_{\ast } \right) ^{\delta /2-1}} \right) \,
\notag \\
& \qquad \times 
W_{\theta }(\bm{\xi },\bm{\xi }_{\ast },I,I_{\ast} \left\vert
\bm{\xi }^{\prime }, \bm{\xi }_{\ast }^{\prime},I^{\prime },I_{\ast }^{\prime }\right. )\,
\text{d}\bm{\xi }_{\ast } \text{d}\bm{\xi }^{\prime }
\text{d}\bm{\xi }_{\ast }^{\prime } \text{d}I_{\ast}\text{d}I^{\prime }\text{d}I_{\ast }^{\prime },
\end{align}
where $W_\theta$ is the transition probability for the collision
$\{(\bm{\xi}, I),\,(\bm{\xi}_*, I_*)\}$ $\to$ $\{(\bm{\xi}', I'),\, (\bm{\xi}_*', I_*')\}$.
Here and below, the following conventional abbreviations
are used:
\begin{equation*}
h = h (\bm{\xi}, I), \;\;
h_{\ast }=h \left( \bm{\xi }_{\ast },I_{\ast }\right), \;\;
h^{\prime }=h \left( \bm{\xi }^{\prime},I^{\prime }\right), \;\; 
h_{\ast }^{\prime }=h\left( \bm{\xi }_{\ast }^{\prime },I_{\ast }^{\prime }\right),
\end{equation*}
for an arbitrary function $h$ of $\bm{\xi}$ and $I$, which
may depend on $t$ and $\bm{x}$;
and it should be recalled that $\delta $ ($\delta \geq 2$) denotes
the number of internal degrees of freedom. 

The transition probability $W_\theta$ in the operator~\eqref{col-op}, which depends
on the parameter $\theta$ specified later, is assumed
to be of the form
\begin{align}\label{Wtheta}
& W_{\theta } (\bm{\xi }, \bm{\xi }_{\ast }, I ,I_{\ast} \left\vert
\bm{\xi }^{\prime }, \bm{\xi }_{\ast }^{\prime}, I^{\prime }, I_{\ast }^{\prime } \right. )
\nonumber \\
& \quad = 4m\left( II_{\ast }\right) ^{\delta /2-1}
\frac{\left\vert \bm{g} \right\vert }
{\left\vert \bm{g}^{\prime } \right\vert }
\bm{\delta}_{3} \left( \bm{\xi }+\bm{\xi }_{\ast }-\bm{\xi }^{\prime }-\bm{\xi }_{\ast }^{\prime }\right) 
\notag \\
& \qquad
\times \bm{\delta}_{1}\left( \frac{m}{2}\left( \left\vert \bm{\xi } \right\vert ^{2}
+\left\vert \bm{\xi }_{\ast}\right\vert ^{2}
- \left\vert \bm{\xi }^{\prime } \right\vert^{2}-\left\vert \bm{\xi }_{\ast }^{\prime }
\right\vert ^{2} \right) -\Delta I \right) \sigma_\theta,
\end{align}
where $\bm{\delta}_3$ and $\bm{\delta}_1$ are the Dirac delta
function in $\mathbb{R}^3$ and $\mathbb{R}$, respectively, and
$\sigma_\theta$ is the scattering cross
section depending on $\theta$ and is expressed as
\begin{align}\label{sigmatheta}
\sigma_\theta =\sigma_\theta \left( \left\vert \bm{g} \right\vert,
\left\vert \cos \phi \right\vert, I, I_{\ast }, I^{\prime }, I_{\ast }^{\prime } \right)
>0 \;\; \text{a.e.},
\end{align}
with
\begin{align*}
\cos \phi = \bm{g}\cdot \bm{g}^{\prime }/
\big( \left\vert \bm{g}\right\vert \left\vert \bm{g}^{\prime }\right\vert \big).
\end{align*}
The parameter $\theta$ ($0 \le \theta \le 1$) is such that the
probability $W_\theta$ reduces to that for standard inelastic
collisions when $\theta=1$ and to
resonant collisions, in which $\Delta I = 0$ holds, when $\theta=0$.

The form of the collision operator $Q_\theta (f, f)$ [cf.~\eqref{col-op}],
proposed in \cite{Be-23b}, is inspired by the probabilistic formulation
for a monatomic gas \cite{K-69,BPS-90}. Furthermore, the form of the
transition probability \eqref{Wtheta},
also proposed in \cite{Be-23b}, is designed in consistency with the
conventional Borgnakke-Larsen representation for
standard collisions (see Sec.~\ref{sec:BL-model}).

It is assumed that the scattering cross section $\sigma_\theta$ satisfies
the microreversibility condition
\begin{align}\label{mr}
&\left( I I_{\ast }\right) ^{\delta /2-1} \left\vert \bm{g} \right\vert^{2}
\sigma_\theta \left( \left\vert \bm{g} \right\vert, \left\vert \cos \phi \right\vert,
I, I_{\ast }, I^{\prime }, I_{\ast }^{\prime } \right)
\notag \\
& \quad = \left( I^{\prime } I_{\ast }^{\prime } \right) ^{\delta /2-1}
\left\vert \bm{g}^{\prime } \right\vert^{2} \sigma_\theta
\left( \left\vert \bm{g}^{\prime } \right\vert, \left\vert \cos \phi \right\vert,
I^{\prime }, I_{\ast}^{\prime }, I, I_{\ast } \right),
\end{align}
and the symmetry relations
\begin{align}\label{sr}
\sigma_\theta \left( \left\vert \bm{g} \right\vert ,\left\vert \cos \phi \right\vert,
I, I_{\ast }, I^{\prime }, I_{\ast }^{\prime } \right)
& = \sigma_\theta \left( \left\vert \bm{g} \right\vert, \left\vert \cos \phi \right\vert,
I, I_{\ast }, I_{\ast }^{\prime }, I^{\prime } \right)
\notag \\
& = \sigma_\theta \left( \left\vert \bm{g} \right\vert, \left\vert \cos \phi \right\vert,
I_{\ast }, I, I_{\ast }^{\prime }, I^{\prime } \right).
\end{align}
The latter is to fulfill the invariance under interchange of molecules in a collision.

The form of the transition probability \eqref{Wtheta} and the properties
\eqref{mr} and \eqref{sr} for the scattering cross section lead to
the following relations:
\begin{align}\label{rel1}
W_{\theta }(\bm{\xi }, \bm{\xi }_{\ast }, I, I_{\ast} \left\vert
\bm{\xi }^{\prime }, \bm{\xi }_{\ast }^{\prime}, I^{\prime }, I_{\ast }^{\prime } \right. )
& = W_{\theta }(\bm{\xi }^{\prime }, \bm{\xi }_{\ast }^{\prime }, I^{\prime }, I_{\ast }^{\prime}
\left\vert \bm{\xi }, \bm{\xi }_{\ast }, I, I_{\ast } \right. )
\notag \\
& = W_{\theta }(\bm{\xi }, \bm{\xi }_{\ast }, I, I_{\ast } \left\vert
\bm{\xi }_{\ast}^{\prime }, \boldsymbol{\xi }^{\prime }, I_{\ast }^{\prime }, I^{\prime} \right. )
\notag \\
& = W_{\theta }(\bm{\xi }_{\ast }, \bm{\xi }, I_{\ast }, I \left\vert
\bm{\xi }_{\ast}^{\prime }, \bm{\xi }^{\prime }, I_{\ast }^{\prime }, I^{\prime} \right. ).
\end{align}

Now, we assume that $\sigma_\theta$ has the following form:
\begin{align}\label{sigmatheta-2}
\sigma_\theta = \theta \sigma_\text{s} + (1-\theta) \sigma_\text{r}\,
\bm{\delta}_1 (\Delta I),
\end{align}
where $\sigma_\text{s}$ and $\sigma_\text{r}$ are independent of $\theta$ and are
assumed to have the form \eqref{sigmatheta} and satisfy the relations \eqref{mr}
and \eqref{sr}. Obviously, $\sigma_\text{s}$ and $\sigma_\text{r}$ are,
respectively, the collision cross section for standard collisions and that for
resonant collisions. Correspondingly, $W_\theta$ is written as
\begin{align}\label{Wtheta-2}
W_\theta = \theta W_\text{s} + (1-\theta) W_\text{r},
\end{align}
where $W_\text{s}$ and $W_\text{r}$ are independent of $\theta$
and are, respectively, the transition probability for standard collisions and
that for resonant collisions.
Then, applying known properties of the Dirac delta function, 
$W_\text{s}$ and $W_\text{r}$ may be transformed into the following form:
\begin{subequations}\label{Ws-Wr-2}
\begin{align}
& W_\text{s}(\bm{\xi }, \bm{\xi }_{\ast }, I, I_{\ast} \left\vert
\bm{\xi }^{\prime }, \bm{\xi }_{\ast }^{\prime}, I^{\prime }, I_{\ast }^{\prime } \right. )
\nonumber \\
& \; = \frac{m}{2}\left( I I_{\ast }\right)^{\delta /2-1} \sigma_\text{s}
\frac{\left\vert \bm{g} \right\vert }{\left\vert \bm{g}^{\prime } \right\vert }
\bm{\delta }_{3} \left( \bm{G}-\bm{G}^{\prime } \right) 
\bm{\delta }_{1} \left( \frac{m}{4}\left( \left\vert \bm{g}\right\vert^{2}
-\left\vert \bm{g}^{\prime } \right\vert^{2} \right)
-\Delta I \right)
\nonumber \\
& \; = \frac{m}{2} \left( I I_{\ast } \right) ^{\delta /2-1} \sigma_\text{s}
\frac{\left\vert \bm{g} \right\vert }{\left\vert \bm{g}^{\prime} \right\vert }
\bm{\delta }_{3} \left( \bm{G}-\bm{G}^{\prime} \right) \bm{\delta }_{1}
\left( E-E^{\prime } \right)
\nonumber \\
& \;= \left( I I_{\ast } \right) ^{\delta /2-1} \sigma_\text{s}
\frac{\left\vert \bm{g} \right\vert }{\left\vert \bm{g}^{\prime } \right\vert^{2}}
\bm{\delta }_{3} \left( \bm{G}-\bm{G}^{\prime } \right) \bm{\delta}_{1}
\left( \sqrt{\left\vert \bm{g} \right\vert^{2}-\frac{4}{m} \Delta I}
-\left\vert \bm{g}^{\prime } \right\vert \right),
\label{Ws-Wr-2-a} \\
& W_\text{r}(\bm{\xi }, \bm{\xi }_{\ast }, I, I_{\ast} \left\vert
\bm{\xi }^{\prime }, \bm{\xi }_{\ast }^{\prime}, I^{\prime }, I_{\ast }^{\prime } \right. )
\notag \\
& \; =
\left( I I_{\ast } \right) ^{\delta /2-1} \sigma_\text{r}
\frac{\left\vert \bm{g} \right\vert }{\left\vert \bm{g}^{\prime } \right\vert^{2}}
\bm{\delta }_{3} \left( \bm{G}-\bm{G}^{\prime } \right) \bm{\delta}_{1}
\left( \sqrt{\left\vert \bm{g} \right\vert^{2}-\frac{4}{m} \Delta I}
-\left\vert \bm{g}^{\prime } \right\vert \right)
\bm{\delta}_1 (\Delta I)
\nonumber \\
& \; =
\left( I I_{\ast } \right) ^{\delta /2-1} \sigma_\text{r}
\left\vert \bm{g} \right\vert^{-1}
\bm{\delta }_{3} \left( \bm{G}-\bm{G}^{\prime } \right) \bm{\delta}_{1}
\left( \left\vert \bm{g} \right\vert
-\left\vert \bm{g}^{\prime } \right\vert \right)
\bm{\delta}_1 (\Delta I). 
\label{Ws-Wr-2-b}
\end{align}
\end{subequations}

For later convenience, we introduce the bilinear operators $Q_\text{s}$ and
$Q_\text{r}$ based on $W_\text{s}$ and $W_\text{r}$, respectively, that is,
\begin{align}\label{Qs-Qr}
Q_\text{s} (f, F) := Q_\theta (f, F)\;\, \text{with}\;\, \theta =1, \quad
Q_\text{r} (f, F) := Q_\theta (f, F)\;\, \text{with}\;\, \theta =0,
\end{align}
and write
\begin{align*}
Q_\theta  (f, F) = \theta Q_\text{s} (f, F) 
+ (1 - \theta) Q_\text{r} (f, F).
\end{align*}
%

\subsection{Borgnakke-Larsen-type model}\label{sec:BL-model}

Borgnakke--Larsen \cite{BL-75} proposed a phenomenological procedure
to simulate the collision process of polyatomic gas molecules by Monte-Carlo
methods. This approach has been widely used in practical computations
using the direct simulation Monte Carlo (DSMC) method \cite{B-76,B-94}.
The Boltzmann collision operator along the lines of the Borgnakke--Larsen
procedure has also been established \cite{BDLP-94} and has been
a target of mathematical study (e.g.,
\cite{DPS-21,BBG-22,GP-23,Be-23a,Be-23b,DOPT-23,BBS-23,BRS-24,BST-24,BBPG-24,BBMS-25}).

In this procedure, it is assumed that, after a collision, the total energy $E$ in the center of mass frame
[see \eqref{EC-2}] is transmitted to the
kinetic energy $(m/4)|\bm{\xi}' - \bm{\xi}_*'|^2$ with the rate $R$ ($\in [0,\, 1]$)
and to the energy of the internal modes $I' + I_*'$ with the rate $1 - R$, that is,
\begin{align}\label{BL-R}
\frac{m}{4} |\bm{\xi}' - \bm{\xi}_*'|^2 = \frac{m}{4} | \bm{g}' |^2 = R E, \qquad
I' + I_*' = (1 - R) E.
\end{align}
The first equation can be written as
$\bm{\xi}' - \bm{\xi}_*'=\bm{g}'=2 \sqrt{RE/m}\, \bm{\sigma}$
with a unit vector $\bm{\sigma}$ ($\in \mathbb{S}^2$). Thus, the
post collisional velocities $\bm{\xi}'$ and $\bm{\xi}_*'$ are expressed as
\begin{align*}
\bm{\xi}' = \bm{G} + \sqrt{\frac{RE}{m}}\, \bm{\sigma}, \qquad
\bm{\xi}_*' = \bm{G} - \sqrt{\frac{RE}{m}}\, \bm{\sigma}, \qquad
\bm{\sigma} = \frac{\bm{g}'}{|\bm{g}'|}.
\end{align*}
In addition, it is assumed that the energy $(1 - R)E$ is divided between
$I'$ and $I_*'$ with the rates $r$ ($\in [0,\, 1]$) and $1-r$, respectively, i.e.,
\begin{align}\label{BL-r}
I' = r (1 - R) E, \qquad I_*' = (1 - r)(1 - R)E.
\end{align}

For resonant collisions, the following relations hold:
\begin{align*}
\bm{\xi}' = \bm{G} + |\bm{g}| \bm{\sigma}/2, \qquad
\bm{\xi}_*' = \bm{G} - |\bm{g}| \bm{\sigma}/2.
\end{align*}
We then assume that the total energy $I + I_*$ of the
internal modes, which is conserved in the collision, is divided between
$I'$ and $I_*'$ with the rates $r$ ($\in [0,\, 1]$) and $1-r$, respectively, after the
collision. That is,
\begin{align*}
I' = r (I + I_*), \qquad I_*' = (1 - r) (I + I_*).
\end{align*}

The numbers $R$ and $r$ thus introduced play the roles of variables
in the Borgnakke-Larsen representation.
With the help of these new variables, the collision operator $Q_\text{s} (f, f)$
can be transformed into the conventional form and $Q_\text{r} (f, f)$
into the corresponding form. 

For this transformation, a series of changes of integration variables is performed.
To be more specific,

\begin{itemize}

\item
$(\bm{\xi}_*,\, \bm{\xi}',\, \bm{\xi}_*',\, I_*,\, I',\, I_*') \to  
(\bm{\xi}_*,\, \bm{g}',\, \bm{G}',\, I_*,\, I',\, I_*')$ with the help of
$\bm{g}' = \bm{\xi}' - \bm{\xi}_*'$ and $\bm{G}' =(\bm{\xi}' + \bm{\xi}_*')/2$;

\item
$(\bm{\xi}_*,\, \bm{g}',\, \bm{G}',\, I_*,\, I',\, I_*') \to  
(\bm{\xi}_*,\, |\bm{g}'|,\, \bm{\sigma},\, \bm{G}',\, I_*,\, I',\, I_*')$
with the help of $\bm{\sigma} = \bm{g}'/|\bm{g}'|$;

\item
$(\bm{\xi}_*,\, |\bm{g}'|,\, \bm{\sigma},\, \bm{G}',\, I_*,\, I',\, I_*') \to  
(\bm{\xi}_*,\, \bm{\sigma},\, \bm{G}',\, I_*,\, R,\, r,\, E')$.
Since the delta function $\bm{\delta}_1( E-E')$ in $W_\text{s}$
indicates $E'=E$, one can write $I' = r(1 - R)E'$, $I_*' = (1-r)(1-R)E'$, and
$|\bm{g}'|^2 = (4/m)RE'$ instead of relations~\eqref{BL-R} and \eqref{BL-r}.
These relations should be used for the above change of variables,
in which $E'$ appears as a new variable.

\end{itemize}

\noindent
By calculating the Jacobian at each step, we obtain
\begin{align*}
\text{d}\bm{\xi}_* \text{d}\bm{\xi}' \text{d}\bm{\xi}_*' \text{d}I_* \text{d}I' \text{d}I_*'
= \frac{2}{m} (1 - R) E'^2 |\bm{g}'|\, \text{d}\bm{\xi}_* \text{d} \bm{\sigma}
\text{d}\bm{G}' \text{d} I_* \text{d}R \text{d}r \text{d}E', 
\end{align*}
where $|\bm{g}'| = \sqrt{(4/m) RE'}$. This relation and the second of the equalities~\eqref{Ws-Wr-2-a} lead to the following expression of $Q_\text{s}$:
\begin{align}\label{Qs-BL}
Q_\text{s} (f,F) &= \frac{1}{2} \int_{\left( \mathbb{R}^{3}\times \mathbb{R}_{+}\right)^{3}}
\left( \frac{ f' F_*' + f_*' F' }{\left( I^{\prime}I_{\ast }^{\prime }\right)
^{\delta /2-1}}-\frac{ f F_* + f_* F }{\left( II_{\ast } \right) ^{\delta /2-1}} \right) \,
\nonumber \\
& \qquad \qquad \qquad \qquad \times 
W_\text{s}\,
\text{d}\bm{\xi }_{\ast } \text{d}\bm{\xi }^{\prime }
\text{d}\bm{\xi }_{\ast }^{\prime } \text{d}I_{\ast} \text{d} I^{\prime }\text{d}I_{\ast }^{\prime }
\nonumber \\
& = \frac{1}{2} \int_{ [0, 1]^2 \times \mathbb{S}^2 \times \mathbb{R}^{3} \times \mathbb{R}_{+} }
\left( \frac{ f' F_*' + f_*' F' }{\left( I^{\prime}I_{\ast }^{\prime }\right)
^{\delta /2-1}}-\frac{ f F_* + f_* F }{\left( II_{\ast } \right) ^{\delta /2-1}} \right) \,
\nonumber \\
& \qquad \qquad \;\; \times 
\left( I I_{\ast } \right)^{\delta /2-1}
|\bm{g}| \sigma_\text{s} (1-R) E^2
\text{d}R\, \text{d}r\, \text{d}\bm{\sigma}\, \text{d}\bm{\xi}_* \, \text{d}I_*.
\end{align}
In the last representation, the fixed variables are $(\bm{\xi}, I)$, and the integration variables
are $(R, r, \bm{\sigma}, I_*, \bm{\xi}_*)$. Noting that $\sigma_\text{s}$ is a function
of $|\bm{g}|$, $|\cos \phi| = |\bm{g}\bm{\cdot}\bm{\sigma}|/|\bm{g}|$, $I$, $I_*$, $I'$, and $I_*'$
and that $f' = f(\bm{\xi}', I')$, $f_*' = f_*' (\bm{\xi}_*', I_*')$, etc.~($t$ and $\bm{x}$ are omitted),
we notice that the integrand contains the variables $\bm{g}$, $\bm{\xi}'$, $\bm{\xi}_*'$,
$I'$, $I_*'$, and $E$ in addition to the fixed and integration variables. Therefore,
$\bm{g}$, $\bm{\xi}'$, $\bm{\xi}_*'$, $I'$, $I_*'$, and $E$ have to be expressed
in terms of the fixed and integration variables, that is,
\begin{align}\label{variables}
\begin{aligned}
& \bm{g} = \bm{\xi} - \bm{\xi}_*,
\\
& \bm{\xi}' = \frac{\bm{\xi} + \bm{\xi}_*}{2} + \sqrt{\frac{RE}{m}}\, \bm{\sigma}, \qquad
\bm{\xi}_*' = \frac{\bm{\xi} + \bm{\xi}_*}{2} - \sqrt{\frac{RE}{m}}\, \bm{\sigma},
\\
& I' = r (1 - R) E, \qquad I_*' = (1 - r)(1 - R)E,
\\
& E = \frac{m}{4} |\bm{\xi} - \bm{\xi}_*|^2 + I + I_*.
\end{aligned}
\end{align}
It should be noted that the operator \eqref{Qs-BL} can be transformed into the
conventional form \cite{GP-23,BST-24}
\begin{align*}
& Q_\text{s} (f, F) = \frac{1}{2}
\int_{ [0, 1]^2 \times \mathbb{S}^2 \times \mathbb{R}^{3} \times \mathbb{R}_{+} }
\left( \frac{ f' F_*' + f_*' F' }{\left( I^{\prime}I_{\ast }^{\prime }\right)
^{\delta /2-1}}-\frac{ f F_* + f_* F }{\left( II_{\ast } \right) ^{\delta /2-1}} \right) \,B_\text{s}
\left( I I_{\ast } \right)^{\delta /2-1}
\nonumber \\
& \qquad \qquad \qquad \times
[r (1 - r)]^{\delta/2-1} (1 - R)^{\delta-1} R^{1/2}\,
\text{d}R\, \text{d}r\, \text{d}\bm{\sigma}\, \text{d}\bm{\xi}_* \, \text{d}I_*,
\end{align*} 
by letting
\begin{align}\label{Bs-BL}
B_\text{s} = \frac{ \sigma_\text{s} |\bm{g}| E^2}
{ [r (1-r)]^{\delta/2-1} (1-R)^{\delta-2} R^{1/2} }.
\end{align}

Similarly, the operator $Q_\text{r}$ for resonant collisions can be
transformed in the following way:
\begin{align}\label{Qr-BL}
Q_\text{r} (f, F) &= \frac{1}{2}
\int_{\left( \mathbb{R}^{3}\times \mathbb{R}_{+}\right)^{3}}
\left( \frac{ f' F_*' + f_*' F' }{\left( I^{\prime}I_{\ast }^{\prime }\right)
^{\delta /2-1}}-\frac{ f F_* + f_* F }{\left( II_{\ast } \right) ^{\delta /2-1}} \right) \,
\nonumber \\
& \qquad \qquad \qquad \qquad \times 
W_\text{r}\,
\text{d}\bm{\xi }_{\ast } \text{d}\bm{\xi }^{\prime }
\text{d}\bm{\xi }_{\ast }^{\prime } \text{d}I_{\ast} \text{d} I^{\prime }\text{d}I_{\ast }^{\prime }
\nonumber \\
& = \frac{1}{2}
\int_{ [0, 1] \times \mathbb{S}^2 \times \mathbb{R}^{3} \times \mathbb{R}_{+} }
\left(\frac{ f' F_*' + f_*' F' }{\left( I^{\prime}I_{\ast }^{\prime }\right)
^{\delta /2-1}}-\frac{ f F_* + f_* F }{\left( II_{\ast } \right) ^{\delta /2-1}} \right) \,
\nonumber \\
& \qquad \qquad \qquad \quad \times 
\left( I I_{\ast } \right)^{\delta /2-1}
|\bm{g}| \sigma_\text{r} (I + I_*)
\text{d}r\, \text{d}\bm{\sigma}\, \text{d}\bm{\xi}_* \, \text{d}I_*,
\end{align}
where  $\sigma_\text{r}$ is a function
of $|\bm{g}|$, $|\cos \phi| = |\bm{g}\bm{\cdot}\bm{\sigma}|/|\bm{g}|$, $I$, $I_*$, $I'$, and $I_*'$, and
\begin{align*}
& \bm{g} = \bm{\xi} - \bm{\xi}_*,
\\
& \bm{\xi}' = \frac{\bm{\xi} + \bm{\xi}_*}{2} + \frac{|\bm{\xi} - \bm{\xi}_*|}{2} \bm{\sigma}, \qquad
\bm{\xi}_*' = \frac{\bm{\xi} + \bm{\xi}_*}{2} - \frac{|\bm{\xi} - \bm{\xi}_*|}{2} \bm{\sigma},
\\
& I' = r (I + I_*), \qquad I_*' = (1 - r)(I + I_*).
\end{align*}
Note that the operator~\eqref{Qr-BL} is recast in the following form:
\begin{align*}
Q_\text{r} (f, F) & = \frac{1}{2}
\int_{ [0, 1] \times \mathbb{S}^2 \times \mathbb{R}^{3} \times \mathbb{R}_{+} }
\left( \frac{ f' F_*' + f_*' F' }{\left( I^{\prime}I_{\ast }^{\prime }\right)
^{\delta /2-1}}-\frac{ f F_* + f_* F }{\left( II_{\ast } \right) ^{\delta /2-1}} \right) \,
\nonumber \\
& \qquad \qquad \qquad \quad \times 
B_\text{r} \left( I I_{\ast } \right)^{\delta /2-1}
[r(1-r)]^{\delta/2-1}
\text{d}r\, \text{d}\bm{\sigma}\, \text{d}\bm{\xi}_* \, \text{d}I_*,
\end{align*}
where
\begin{align}\label{Br-BL}
B_\text{r} = \frac{ \sigma_\text{r} |\bm{g}| (I + I_*)}{[r(1-r)]^{\delta/2-1}}.
\end{align}
%

\subsection{Collision invariants and equilibrium distributions}\label{sec:col-inv}

In this section, the properties of the collision operator
$Q_\theta (f, f)$ are discussed. For non-resonant collisions ($\theta \ne 0$),
they are basically the same as those for the standard collision
operator $Q_\text{s} (f, f)$ [i.e., $Q_\theta (f, f)$ with $\theta=1$] discussed
in \cite{Be-23b}. Although the case of resonant collisions ($\theta=0$)
has to be treated separately in some cases, the treatment is straightforward.
Therefore, we mainly summarize the results without proof.

Let us define the measure $\text{d}A_\theta$ by 
\begin{equation}\label{dAtheta}
\text{d}A_{\theta }=W_{\theta }(\bm{\xi }, \bm{\xi }_{\ast},
I, I_{\ast } \left\vert \bm{\xi }^{\prime }, \bm{\xi }_{\ast}^{\prime },
I^{\prime }, I_{\ast }^{\prime } \right. )\,
\text{d}\bm{\xi }\,\text{d}\bm{\xi }_{\ast } \text{d}\bm{\xi }^{\prime } \text{d}\bm{\xi }_{\ast }^{\prime }
\text{d}I \text{d}I_{\ast } \text{d}I^{\prime } \text{d}I_{\ast }^{\prime }.
\end{equation}
Then, the weak form $\left( Q_\theta (f, F), g\right)$ of the bilinear
operator $Q_\theta (f, F)$ is expressed as
\begin{equation*}
\left( Q_{\theta }(f, F),g \right) = \frac{1}{2}
\int_{\left( \mathbb{R}^{3}\times \mathbb{R%
}_{+}\right) ^{4}}\left( \frac{ f' F_*' +f_*' F' }{\left(
I^{\prime }I_{\ast }^{\prime }\right) ^{\delta /2-1}}-\frac{ f F_* + f_* F }{%
\left( II_{\ast }\right) ^{\delta /2-1}}\right) g\,\text{d}A_{\theta },
\end{equation*}
where $g=g\left( \bm{\xi }, I \right) $ is any function
such that the integral is defined.

The following lemma follows directly from the relations \eqref{rel1}:
\begin{lemma}
\label{L0}
The measure $\mathrm{d}A_{\theta }$ is invariant under the interchanges
of variables
\begin{align}\label{tr1}
\begin{aligned}
(i)& \qquad\left( \bm{\xi }, \bm{\xi }_{\ast }, I, I_{\ast }\right)
\leftrightarrow \left( \bm{\xi }^{\prime }, \bm{\xi }_{\ast}^{\prime },
I^{\prime }, I_{\ast }^{\prime } \right),
\\
(ii)& \qquad\left( \bm{\xi }, I \right) \leftrightarrow \left( \bm{\xi }_{\ast },
I_{\ast }\right), 
\\
(iii)& \qquad\left( \bm{\xi }^{\prime }, I^{\prime } \right) \leftrightarrow
\left( \bm{\xi }_{\ast }^{\prime }, I_{\ast }^{\prime } \right),
\end{aligned}
\end{align}
respectively.
\end{lemma}
This leads to the following proposition:
\begin{proposition}
\label{P1}
Let $g=g\left( \bm{\xi }, I \right)$ be such that
\begin{equation*}
\int_{\left( \mathbb{R}^{3}\times \mathbb{R}_{+}\right)^{4}}
\left( \frac{ f' F_*' + f_*' F' }{\left( I^{\prime }I_{\ast }^{\prime }\right)^{\delta /2-1}}
-\frac{ f F_* + f_* F }{\left( II_{\ast }\right)^{\delta /2-1}} \right) g\,\mathrm{d}A_{\theta }
\end{equation*}
is defined. Then, it holds that
\begin{align}\label{sym-relation}
\left( Q_{\theta }(f, F), g \right)
& = \frac{1}{8} \int_{\left( \mathbb{R}^{3}\times \mathbb{R}_{+}\right)^{4}}
\left( \frac{ f' F_*' + f_*' F' }{\left( I^{\prime }I_{\ast }^{\prime }\right)^{\delta /2-1}}
- \frac{ f F_* + f_* F }{\left( II_{\ast }\right)^{\delta /2-1}} \right) \nonumber
\\
& \qquad \qquad \qquad \qquad \quad \times
\left( g + g_{\ast } - g^{\prime } - g_{\ast }^{\prime }\right) \, \mathrm{d}A_{\theta }.
\end{align}
\end{proposition}

In accordance with Proposition \ref{P1}, we introduce the concept of a collision invariant for the collision operator $Q_\theta (f, f)$ as
\begin{definition}
A function $g=g\left( \bm{\xi }, I \right) \ $is a collision invariant if
\begin{equation}\label{CI}
\left( g+g_{\ast }-g^{\prime }-g_{\ast }^{\prime }\right) W_{\theta }(%
\boldsymbol{\xi },\boldsymbol{\xi }_{\ast },I,I_{\ast }\left\vert 
\boldsymbol{\xi }^{\prime },\boldsymbol{\xi }_{\ast }^{\prime },I^{\prime
},I_{\ast }^{\prime }\right. )=0, \;\text{ a.e.}
\end{equation}
holds.
\end{definition}
When $\theta \ne 0$, it is obvious that $1,$ $\xi _i,$ ($i=1,\,2,\,3$), and
$m \left\vert \bm{\xi } \right\vert^{2}+2I$ are collision invariants
due to the conservation of mass, momentum, and total energy [cf.~\eqref{conservation}].
However, it should be noted that, when $\theta =0$, each of
$\left\vert \bm{\xi } \right\vert^{2}$ and $I$ is a collision invariant,
since not only the conservation of the total energy but also the
separate conservation of the kinetic energy and the energy of the
internal modes holds [cf.~\eqref{conservation-theta=0}].
In fact, we have the following proposition
corresponding to Proposition 2 in \cite{Be-23b}~(see \cite{BDLP-94,BRS-24}; cf.~\cite{AC-90}).

\begin{proposition}
\label{P2}
The vector space of collision invariants is generated by 
\begin{equation*}
\left\{ 1, \xi _1, \xi _2, \xi _3, m \left\vert \bm{\xi } \right\vert^{2} + 2 I \right\},
\end{equation*}
in the non-resonant case $\theta \ne 0$ and%
\begin{equation*}
\left\{ 1 ,\xi _1, \xi _2, \xi _3, \left\vert \boldsymbol{\xi }\right\vert^{2}, I \right\},
\end{equation*}
in the resonant case $\theta =0$.
\end{proposition}

In addition, following the line of Sec.~2.2 in \cite{Be-23b}, we have the
following properties related to $Q_\theta (f, f)$.
\begin{proposition}
Let $\mathcal{W}_\theta [f]$ be the functional defined by
\begin{align*}
\mathcal{W}_\theta [f] = \left( Q_{\theta }(f, f), \log \left(
I^{1-\delta /2} f \right) \right).
\end{align*}
Then, it follows that
\begin{align*}
\mathcal{W}_\theta [f] \le 0.
\end{align*}
\end{proposition}
\begin{proposition}\label{prop-4}
The following (i), (ii), and (iii) are equivalent.

\begin{enumerate}

\item[(i)] $\mathcal{W}_\theta [f] = 0.$

\item[(ii)] $Q_\theta (f, f) = 0.$

\item[(iii)] $f$ is the equilibrium distribution
$M_\text{s}$ (for $\theta \ne 0$) or $M_\text{r}$ (for $\theta = 0$) given as follows:
\begin{align}\label{Ms}
M_\mathrm{s} = \frac{n I^{\delta/2-1}}{(2\pi k_\mathrm{B} T/m)^{3/2}
(k_\mathrm{B} T)^{\delta/2}
\Gamma (\delta/2)}
\exp \left( - \frac{m |\bm{\xi} - \bm{u}|^2 + 2I}{2k_\mathrm{B} T}
\right), \quad (\theta \ne 0),
\end{align}
where [cf.~relations~\eqref{macro}]
\begin{align*}
& n = (1, M_\mathrm{s}), \qquad \bm{u} = \frac{1}{n} (\bm{\xi}, M_\mathrm{s}), \qquad
T = T_\mathrm{tr} = T_\mathrm{int},
\\
& T_\mathrm{tr} = \frac{2}{3 k_\mathrm{B}} n \left( \frac{m}{2} |\bm{\xi} - \bm{u}|^2, M_\mathrm{s} \right), \qquad
T_\mathrm{int} = \frac{2}{\delta k_\mathrm{B} n} \left( I, M_\mathrm{s} \right),
\end{align*}
and $\Gamma(s) = \int_0^\infty x^{s-1} e^{-x}$ is the gamma function; and
\begin{align}\label{Mr}
M_\mathrm{r} & = \frac{n I^{\delta/2-1}}{(2\pi k_\mathrm{B} T_\mathrm{tr}/m)^{3/2}
(k_\mathrm{B} T_\mathrm{int})^{\delta/2}
\Gamma (\delta/2)}
\exp \left( - \frac{m |\bm{\xi} - \bm{u}|^2}{2k_\mathrm{B} T_\mathrm{tr}}
- \frac{I}{k_\mathrm{B} T_\mathrm{int}}
\right), \quad (\theta=0),
\end{align}
where [cf.~relations~\eqref{macro}]
\begin{align*}
n = (1, M_\mathrm{r}), \quad \bm{u} = \frac{1}{n} (\bm{\xi}, M_\mathrm{r}), \\
T_\mathrm{tr} = \frac{2}{3 k_\mathrm{B} n} \left( \frac{m}{2} |\bm{\xi} - \bm{u}|^2, M_\mathrm{r} \right),\quad
T_\mathrm{int} = \frac{2}{\delta k_\mathrm{B} n} \left( I, M_\mathrm{r} \right).
\end{align*}

\end{enumerate}

\end{proposition}

\noindent
The distribution $M_\text{s}$
indicates the local equilibrium state in the non-resonant case ($\theta \ne 0$) with the molecular
number density $n$, the flow velocity $\bm{u}$, and the single temperature $T$; and
$M_\text{r}$ indicates the local equilibrium state in the resonant case
($\theta=0$) with the molecular number density $n$,
the flow velocity $\bm{u}$, and two distinct temperatures, i.e., the translational temperature
$T_\text{tr}$ and the internal temperature $T_\text{int}$. 
Note that $M_\text{r}$ reduces to $M_\text{s}$ when $T_\text{tr} = T_\text{int} = T$.

\begin{remark}
Introducing the $\mathcal{H}$-functional
\begin{equation*}
\mathcal{H}\left[ f\right] =\left( f,\log \left( I^{1-\delta /2}f\right)
\right) \text{,}
\end{equation*}%
an $\mathcal{H}$-theorem can be obtained (cf.~\cite{BDLP-94, BBBD-18, GP-23}).
\end{remark}

\subsection{Linearized collision operator}\label{sec:linear-col-op}

Recall that there are two local equilibrium distributions: $M_\text{s}$ for $\theta \ne 0$ and $M_\text{r}$ for $\theta = 0$ [see~\eqref{Ms}~and~\eqref{Mr}]. Let $M$
stand for $M_\text{s}$ when $\theta \ne 0$ and $M_\text{r}$ when $\theta=0$, i.e.,
\begin{align}\label{M}
M = \left\{
\begin{array}{l}
M_\text{s} \;\;\;\; (\theta \ne 0), \\
M_\text{r} \;\;\;\; (\theta = 0).
\end{array}
\right.
\end{align}
We consider deviations from $M$ as
\begin{equation*} 
f=M (1 + h),
\end{equation*}
and define the linearized collision operator $\mathcal{L}_{\theta }$ by
\begin{equation}\label{dec2}
\mathcal{L}_{\theta }h=-2M^{-1}
Q_{\theta }(M, M h)
= \nu _{\theta } h - K_{\theta} \left( h\right),
\end{equation}
where
\begin{subequations}\label{dec1}
\begin{align}
\nu _{\theta } &=\int_{\left( \mathbb{R}^{3}\times \mathbb{R}_{+}\right)^{3}}\!
\frac{M_\ast }{\left( II_{\ast }\right)^{\delta /2-1}}W_{\theta }
\text{d}\bm{\xi }_{\ast } \text{d}\bm{\xi }^{\prime } \text{d}\bm{\xi }_{\ast }^{\prime }
\text{d}I_{\ast } \text{d}I^{\prime } \text{d}I_{\ast }^{\prime },
\label{dec1-a} \\
K_{\theta }\left( h\right) &=\int_{\left( \mathbb{R}^{3}\times \mathbb{R}_{+}\right)^{3}}
\frac{\left( M M_\ast M^{\prime } M_\ast^\prime \right)^{1/2}}
{M \left( II_{\ast }I^{\prime }I_{\ast }^{\prime }\right)
^{\delta /4-1/2}}W_{\theta }
\nonumber \\
& \qquad \qquad \qquad \times
\left( h^{\prime }+h_{\ast }^{\prime }-h_{\ast}\right) \,
\text{d}\bm{\xi }_{\ast } \text{d}\bm{\xi }^{\prime } \text{d}\bm{\xi }_{\ast }^{\prime }
\text{d}I_{\ast } \text{d}I^{\prime } \text{d}I_{\ast }^{\prime }.
\label{dec1-b}
\end{align}
\end{subequations}
The following lemma follows immediately by Lemma \ref{L0}.
\begin{lemma}
\label{L1}
The measure
\begin{equation*}
\mathrm{d} \widetilde{A}_{\theta }
= \frac{\left( M M_\ast M^\prime M_\ast^{\prime} \right)^{1/2}}
{\left( II_{\ast }I^{\prime }I_{\ast }^{\prime }\right)^{\delta /4-1/2}}
\mathrm{d}A_{\theta }
\end{equation*}
is invariant under the interchanges of variables~\eqref{tr1},
respectively.
\end{lemma}

The weak form of the linearized collision operator $\mathcal{L}_{\theta }$
reads
\begin{equation*}
\left( \mathcal{L}_{\theta }h, M g \right)
=\int_{\left( \mathbb{R}^{3}\times \mathbb{R}_{+}\right)^{4}}
\left( h+h_{\ast }-h^{\prime }-h_{\ast }^{\prime}\right) g\,\text{d}\widetilde{A}_{\theta },
\end{equation*}
for $g=g\left( \boldsymbol{\xi },I\right) $ such that the integral is
defined. Applying Lemma \ref{L1}, we obtain the following lemma.

\begin{lemma}
\label{L2}
Let $g=g\left( \boldsymbol{\xi },I\right) $ be such that
\begin{equation*}
\int_{\left( \mathbb{R}^{3}\times \mathbb{R}_{+}\right) ^{4}}
\left( h+h_{\ast }-h^{\prime }-h_{\ast }^{\prime }\right) g\,
\mathrm{d}\widetilde{A}_{\theta}
\end{equation*}
is defined. Then
\begin{equation*}
\left( \mathcal{L}_{\theta }h, M g \right)
=\frac{1}{4}\int_{\left( \mathbb{R}^{3}\times \mathbb{R}_{+}\right)^{4}}
\left( h+h_{\ast }-h^{\prime }-h_{\ast}^{\prime }\right)
\left( g+g_{\ast }-g^{\prime }-g_{\ast }^{\prime }\right)
\mathrm{d}\widetilde{A}_{\theta }.
\end{equation*}
\end{lemma}
Therefore, it follows the following proposition:
\begin{proposition}\label{prop:positivity-kernel}
The linearized collision operator is symmetric and nonnegative,
with respect to the weighted inner product
$\left(\, \mathcal{\cdot }\,,\, M \cdot \,\right)$, i.e.,
\begin{equation*}
\left( \mathcal{L}_{\theta }h, M g \right)
=\left( M h, \mathcal{L}_{\theta}g\right)
\; \text{ and }\;
\left( \mathcal{L}_{\theta }h, M h \right) \geq 0,
\end{equation*}
and the kernel of $\mathcal{L}_{\theta }$, $\ker \mathcal{L}_{\theta }$, is
generated by%
\begin{equation*}
\left\{ 1,\xi _{x},\xi _{y},\xi _{z},m\left\vert \boldsymbol{\xi }%
\right\vert ^{2}+2I\right\},
\end{equation*}%
in the non-resonant case $0<\theta \leq 1$, where $M=M_\mathrm{s}$, and%
\begin{equation*}
\left\{ 1,\xi _{x},\xi _{y},\xi _{z},\left\vert \boldsymbol{\xi }\right\vert
^{2},I\right\},
\end{equation*}%
in the resonant case $\theta =0$, where $M=M_\mathrm{r}$.
\end{proposition}
\smallskip
\begin{proof}
By Lemma \ref{L2}, it is immediate that 
$\left( \mathcal{L}_{\theta}h, M g\right)
=\left( M h,\mathcal{L}_{\theta }g\right) $ and
$\left( \mathcal{L}_{\theta }h, M h\right) \geq 0.$
Furthermore, $h\in \ker \mathcal{L}_{\theta}$ indicates
$\left( \mathcal{L}_{\theta }h, M h\right) =0$, which
means that $h$ satisfies relation~\eqref{CI},
i.e., $h$ is a collision invariant. {Conversely,
if $h$ is a collision invariant, then $h\in \ker \mathcal{L}_{\theta}$
due to equalities~\eqref{dec2} and \eqref{dec1}. Thus, the} last
part of the lemma follows by Proposition \ref{P2}.
\end{proof}

\noindent
Here, we introduce the following notation, which will be used later:
\begin{align}\label{perp}
\begin{aligned}
& M \ker \mathcal{L}_\theta := 
\left\{
\begin{array}{l}
\rm{Span}\{M_\text{s}, M_\text{s}\xi_x, M_\text{s}\xi_y, M_\text{s}\xi_z, M_\text{s}(m|\bm{\xi}|^2 + 2I) \}
\qquad (\theta \ne 0),
\\[1mm]
\rm{Span}\{M_\text{r}, M_\text{r}\xi_x, M_\text{r}\xi_y, M_\text{r}\xi_z, M_\text{r}|\bm{\xi}|^2, M_\text{r}I) \}
\qquad (\theta = 0),
\end{array}
\right.
\\[2mm]
& (M \ker \mathcal{L}_\theta)^\perp :=
\left(
\begin{array}{l}
\text{orthogonal complement of}\; M \ker \mathcal{L}_\theta\; \text{in}\;
L^{2}\left( M \mathrm{d}\bm{\xi} \mathrm{d}I \right)\;
\\
\text{with respect to the inner product}\; (\;\cdot\;, \;\cdot\,)
\end{array}
\right).
\end{aligned}
\end{align}
%

\subsection{Fredholmness of the linearized collision operator}\label{sec:Fredholm}

The discussion so far has been based on a general form of
the bilinear operator~\eqref{col-op} with~\eqref{Wtheta}.
To proceed further, we need to specify models
for $\sigma_\text{s}$ and $\sigma_\text{r}$ in~\eqref{sigmatheta-2}.
Hereafter, the following $\sigma_\text{s}$ and $\sigma_\text{r}$ are assumed:
\begin{subequations}\label{sigma-model}
\begin{align}
& \sigma_\text{s}=C_\text{s}\left( \frac{m}{4}\right) ^{\left( \beta +1\right) /2}
\left( I+I_{\ast }\right) ^{\alpha } \left( I^{\prime }+I_{\ast }^{\prime}\right) ^{\alpha }
\dfrac{\left( I^{\prime }I_{\ast }^{\prime }\right)^{\delta /2-1}}
{E^{\delta +\alpha +\left( \beta +1\right) /2}}\left\vert \bm{g}\right\vert^{\beta -1}
\left\vert \bm{g}^{\prime }\right\vert^{\beta +1},
\label{sigma-model-a} \\
& \sigma_\text{r} = C_\text{r} \frac{\left( I' I_*' \right)^{\delta/2-1}}{(I + I_*)^{\delta -1-\alpha}}
|\bm{g}|^{\beta - 1},
\label{sigma-model-b}
\end{align}
\end{subequations}
where
\begin{align}\label{Cs-Cr}
C_\text{s}=\frac{\Gamma \big( \delta +\alpha +( \beta +3 ) /2\big) }
{\Gamma \big( ( \beta +3 ) /2 \big) \Gamma \left( \delta +\alpha \right) }C_\text{r},
\end{align}
$C_\text{s}$ and $C_\text{r}$ are positive constants,
and $\alpha$ and $\beta$ are real numbers such that
$\alpha \in [0, \delta /2 )$ and $\beta \in [0, 1]$.
If $\sigma_\text{s}$ and $\sigma_\text{r}$ given by~\eqref{sigma-model}
are used in~\eqref{Bs-BL} and \eqref{Br-BL}, then $B_\text{s}$
and $B_\text{r}$ are obtained, respectively, in the following form:
\begin{subequations}\label{Bs-Br-model}
\begin{align}
& B_\text{s}=C_\text{s}\left( \frac{m}{4}\right)^{\beta /2}\left( I+I_{\ast }\right)^{\alpha }
\left\vert \bm{g}\right\vert ^{\beta }\dfrac{\left( I^{\prime}+I_{\ast }^{\prime }\right) ^{\alpha }
\left\vert \bm{g}^{\prime}\right\vert ^{\beta }}
{E^{\alpha +\beta /2}}=C_\text{s}\left( I+I_{\ast }\right)^{\alpha }
\left\vert \bm{g}\right\vert ^{\beta }R^{\beta/2 }\left(1-R\right) ^{\alpha },
\\
& B_\text{r} = C_\text{r} (I + I_*)^\alpha |\bm{g}|^\beta.
\end{align}
\end{subequations}
The forms of $\sigma_\text{s}$ and $\sigma_\text{r}$ (thus those of
$B_\text{s}$ and $B_\text{r}$) are chosen for convenience of
later mathematical analysis rather than for physical reasons.
One might say that the kernels $B_\text{s}$ and $B_\text{r}$ given by~\eqref{Bs-Br-model} is a generalization of the variable hard-sphere molecules for a monatomic gas
(the case of $\beta=1$ corresponds to
a generalization of the hard-sphere molecules). 
Then, we have the results summarized in the following
[recall that $M$ is defined by~\eqref{M}].
  
Combining the results in \cite{Be-23b, Be-24a} with the
compactness results in the resonant case (cf.~\cite{B-24, BBS-23}),
we obtain the following result:

\begin{theorem}
\label{Thm1}
The operator $K_{\theta }$ [see~\eqref{dec1-b}] is a self-adjoint compact operator on $L^{2}\left( M \mathrm{d}\bm{\xi} \mathrm{d}I \right)$.
\end{theorem}

\noindent
Noting that the sum of two self-adjoint operators, at least one of
which is bounded, is self-adjoint itself, one arrives at the following conclusion:

\begin{corollary}
\label{Cor1}
The linearized collision operator $\mathcal{L}_{\theta }$ is a closed,
densely defined, and self-adjoint operator on
$L^{2}\left( M \mathrm{d}\bm{\xi} \mathrm{d}I \right)$.
\end{corollary}

\noindent
Then, the following decomposition of the linearized collision operator is obtained.

\begin{theorem}
\label{Thm2} The linearized collision operator $\mathcal{L}_{\theta }$
can be expressed in the form
\begin{align*}
\mathcal{L}_{\theta }=\Lambda _{\theta }-K_{\theta },
\end{align*}
where $\Lambda_{\theta }$ is the positive multiplication operator
defined by $\Lambda _{\theta }f=\nu _{\theta }f$ with
$\nu _{\theta}=\nu _{\theta }(\left\vert \bm{\xi}\right\vert, I)$
defined by \eqref{dec1-a}, and $K_{\theta }$ is the compact operator
on $L^{2} \left( M \mathrm{d}\bm{\xi} \mathrm{d}I \right)$ defined by
 \eqref{dec1-b}.
Moreover, there exist positive numbers $\nu _{\theta }^{-}$ and $\nu_{\theta }^{+}$,
$0<\nu _{\theta }^{-}<\nu _{\theta }^{+}$, such that for all
$\boldsymbol{\xi }\in \mathbb{R}^{3}$ and for all $\theta \in [0, 1]$,
\begin{equation}
\nu _{\theta}^{-}\left( 1+\left\vert \boldsymbol{\xi }\right\vert \right) ^{\beta
}\left( 1+I\right) ^{\alpha }\leq \nu _{\theta}(\left\vert \boldsymbol{\xi }%
\right\vert ,I)\leq \nu _{\theta}^{+}\left( 1+\left\vert \boldsymbol{\xi }%
\right\vert \right) ^{\beta }\left( 1+I\right) ^{\alpha }.
\label{cocfa}
\end{equation}
\end{theorem}

\smallskip
\noindent
The bounds \eqref{cocfa} are obtained by standard arguments
(see Appendix \ref{sec:nu-bounds}).

The multiplication operator $\Lambda _{\theta }$ is a Fredholm operator if
and only if it is coercive. Since the set of Fredholm operators is closed
under the addition of compact operators, we obtain the following result.

\begin{corollary}
\label{Cor2}
The linearized collision operator $\mathcal{L}_{\theta }$ with parameters
$\left( \alpha ,\beta \right) \in \left[ 0,\delta/2\right) \times \left[ 0,1\right]$
is a Fredholm operator on $L^{2}\left( M \mathrm{d}\bm{\xi}\mathrm{d}I \right)$
with domain 
\begin{equation*}
\mathcal{D}\left( \mathcal{L}_{\theta}\right) =L^{2}\left( \left( 1+\left\vert \bm{\xi }\right\vert \right)^{\beta }
\left( 1+I \right)^{\alpha } M \mathrm{d}\bm{\xi} \mathrm{d}I \right),
\end{equation*}
for all $\theta \in [0, 1]$.
\end{corollary}
\begin{remark}\label{solvability}
Consider the integral equation $\mathcal{L}_\theta h = g$, where $h(\bm{\xi}, I)$
is an unknown function and $g(\bm{\xi}, I)$ a given function. According to Corollary \ref{Cor2},
the integral equation has a unique solution
$h(\bm{\xi}, I) \in L^{2}\left( M \mathrm{d}\bm{\xi}\mathrm{d}I \right) \cap (M \ker \mathcal{L}_\theta)^\perp$
if and only if $g (\bm{\xi}, I) \in \mathcal{D} \left( \mathcal{L}_\theta \right) \cap (M \ker \mathcal{L}_\theta)^\perp$.
\end{remark}
%
\section{Nearly resonant collisions and two-temperature fluid models}\label{sec:nearly-resonant}

In this section, we consider the case where resonant collisions are dominant, that is, the interaction between the translational and internal modes are weak, and derive fluid-dynamic equations with two temperatures by appropriate parameter settings.

\subsection{Preliminaries}\label{sec:prelim}

\subsubsection{Collision frequency and mean free path}

As a preparation, we first define the collision frequency and the mean free path of the gas molecules. If the gain and loss terms in the collision operator
$Q_\theta (f, f)$ [cf.~\eqref{col-op}] are assumed to be separable, the collision frequency $\nu(\bm{\xi}, I)$ is given by the coefficient of $-f$ in the loss term, i.e.,
\begin{align}\label{nu}
\nu (\bm{\xi}, I) & = 
\int_{\left( \mathbb{R}^{3}\times \mathbb{R}_{+}\right)^{3}}
\frac{f_{\ast }}{\left( II_{\ast } \right) ^{\delta /2-1}} 
W_{\theta }(\bm{\xi },\bm{\xi }_{\ast },I,I_{\ast} \left\vert
\bm{\xi }^{\prime }, \bm{\xi }_{\ast }^{\prime},I^{\prime },I_{\ast }^{\prime }\right. )\,
\text{d}\bm{\xi }_{\ast } \text{d}\bm{\xi }^{\prime }
\text{d}\bm{\xi }_{\ast }^{\prime } \text{d}I_{\ast}\text{d}I^{\prime }\text{d}I_{\ast }^{\prime }.
\end{align}
Let us denote by $n_0$ and $T_0$ the reference number density and
temperature, respectively, and by $M_0 (|\bm{\xi}|, I)$ the equilibrium distribution
$M_\text{s}$ at number density $n_0$, temperature $T_0$, and flow
velocity $0$, that is,
\begin{align*}
M_0 (|\bm{\xi}|, I) =  \frac{n_0 I^{\delta/2-1}}{(2\pi k_\mathrm{B} T_0/m)^{3/2}
(k_\mathrm{B} T_0)^{\delta/2}
\Gamma (\delta/2)}
\exp \left( - \frac{m |\bm{\xi} |^2 + 2I}{2k_\mathrm{B} T_0}
\right),
\end{align*}
The reference collision frequency $\nu_0 (\bm{\xi}, I)$ is defined
by~\eqref{nu} with $f_* = M_{0*} = M_0 (|\bm{\xi}_*|, I_*)$, i.e.,
\begin{align*}
\nu_0 (\bm{\xi}, I)  = 
\int_{\left( \mathbb{R}^{3}\times \mathbb{R}_{+}\right)^{3}}
\left( II_{\ast } \right) ^{1-\delta /2} M_{0 *} W_{\theta }\,
\text{d}\bm{\xi }_{\ast } \text{d}\bm{\xi }^{\prime }
\text{d}\bm{\xi }_{\ast }^{\prime } \text{d}I_{\ast}\text{d}I^{\prime }\text{d}I_{\ast }^{\prime }.
\end{align*}
If the average of $\nu_0 (\bm{\xi}, I)$ with respect to the equilibrium
distribution $M_0 (|\bm{\xi}|, I)$ is denoted by $\overline{\nu}_0$, it is
written as
\begin{align*}
\overline{\nu}_0 = \frac{1}{n_0} \int_{\mathbb{R}^3 \times \mathbb{R}_+}
\nu_0 (\bm{\xi}, I) M_0 (|\bm{\xi}|, I) \text{d}\bm{\xi } \text{d}I
= n_0 W_{\theta 0},
\end{align*}
where
\begin{align*}
W_{\theta 0} = \frac{1}{n_0^2} \int_{\left( \mathbb{R}^{3}\times \mathbb{R}_{+}\right)^{3}}
\left( II_{\ast } \right) ^{1-\delta /2} M_0 M_{0 *} \, \text{d} A_\theta,
\end{align*}
with $\text{d} A_\theta$ defined by~\eqref{dAtheta}. Then, we define
the reference mean free time $\tau_0$ and the reference mean free path $l_0$ by
\begin{align*}
\tau_0 = \frac{1}{\overline{\nu}_0} = \frac{1}{n_0 W_{\theta 0}}, \qquad
l_0 = \xi_0 \tau_0 = \frac{\xi_0}{n_0 W_{\theta 0}},
\end{align*}
where $\xi_0 = \sqrt{k_\text{B} T_0/m}$, which is of the order of the average
thermal speed of the gas molecules at temperature $T_0$, is the reference
speed.

\subsubsection{Nondimensionalization}\label{sec:dimless}

In addition to the reference number density $n_0$, reference temperature $T_0$,
and reference speed $\xi_0$ already appeared, we introduce the reference pressure
$p_0=k_\text{B} n_0 T_0$, reference time $t_0$, and reference length $L_0$.
Then, the dimensionless quantities $\, (\widehat{t},$ $\widehat{\bm{x}},$ $\widehat{\bm{\xi}},$
$\widehat{I},$
$\widehat{f},$ $\widehat{n},$ $\widehat{\rho},$ $\widehat{\bm{u}},$ $\widehat{e},$ $\widehat{e}_\text{tr},$
$\widehat{e}_\text{int},$ $\widehat{T},$ $\widehat{T}_\text{tr},$ $\widehat{T}_\text{int},$
$\widehat{p}_{ij},$ $\widehat{q}_i,$ $\widehat{q}_{\text{(tr)}i},$ $\widehat{q}_{\text{(int)}i},$
$\widehat{W}_\theta)$ 
corresponding to  $(t,$ $\bm{x},$ $\bm{\xi},$ $I,$ $f,$ $n,$ $\rho,$ $\bm{u},$
$e,$ $e_\text{tr},$ $e_\text{int},$ $T,$ $T_\text{tr},$ $T_\text{int},$ $p_{ij},$
$q_i,$ $q_{\text{(tr)}i},$ $q_{\text{(int)}i},$ $W_\theta)$ are
introduced by the following relations:
\begin{align}\label{dimless-variables}
\begin{aligned}
& \widehat{t} = \frac{t}{t_0}, \qquad \widehat{\bm{x}} = \frac{\bm{x}}{L_0}, \qquad
\widehat{\bm{\xi}} = \frac{\bm{\xi}}{\xi_0}, \qquad
\widehat{I} = \frac{I}{k_\text{B} T_0} = \frac{I}{m\xi_0^2},
\\
& \widehat{f} = \frac{m \xi_0^5}{n_0} f, \qquad \widehat{n} = \frac{n}{n_0}, \qquad
\widehat{\rho} = \frac{\rho}{mn_0}, \qquad \widehat{\bm{u}} = \frac{\bm{u}}{\xi_0},
\\
& \left(\widehat{e},\, \widehat{e}_\text{tr},\, \widehat{e}_\text{int} \right)
= \frac{1}{m \xi_0^2} \left( e,\, e_\text{tr},\, e_\text{int} \right), \qquad
\left( \widehat{T},\, \widehat{T}_\text{tr},\, \widehat{T}_\text{int} \right)
= \frac{1}{T_0} \left( T,\, T_\text{tr},\, T_\text{int} \right),
\\
& \widehat{p}_{ij} = \frac{p_{ij}}{p_0} , \qquad
\left( \widehat{q}_i,\, \widehat{q}_{\text{(tr)}i},\, \widehat{q}_{\text{(int)}i} \right)
= \frac{1}{p_0 \xi_0} \left( q_i,\, q_{\text{(tr)}i},\, q_{\text{(int)}i} \right),
\\
& \widehat{W}_\theta = \frac{\xi_0^6}{(m\xi_0^2)^{\delta-4} W_{\theta 0}} W_\theta.
\end{aligned}
\end{align}
The variables ($\bm{\xi}_*$, $\bm{\xi}'$, $\bm{\xi}_*'$) and ($I_*$, $I'$, $I_*'$)
involved in binary collisions are nondimensionalized in the same way
as $\bm{\xi}$ and $I$, and the resulting dimensionless variables
are denoted by ($\widehat{\bm{\xi}}_*$, $\widehat{\bm{\xi}}'$, $\widehat{\bm{\xi}}_*'$)
and ($\widehat{I}_*$, $\widehat{I}'$, $\widehat{I}_*'$), respectively.

By the use of relations~\eqref{dimless-variables}, the dimensionless version
of relations~\eqref{macro} is obtained as follows:
\begin{align}\label{macro-dimless}
\begin{aligned}
& \widehat{n} = \widehat{\rho} = (1,\, \widehat{f}), \qquad \widehat{u}_i = \frac{1}{\widehat{n}} (\xi_i,\, \widehat{f}),
\\
& \widehat{e} = \widehat{e}_\text{tr} + \widehat{e}_\text{int}, \qquad
\widehat{e}_\text{tr} = \frac{1}{\widehat{n}} \left( \frac{1}{2}|\widehat{\bm{\xi}}-\widehat{\bm{u}}|^2,\, \widehat{f} \right), \qquad
\widehat{e}_\text{int} = \frac{1}{\widehat{n}} (\widehat{I},\, \widehat{f}),
\\
& \widehat{T} = \frac{3\widehat{T}_\text{tr} + \delta \widehat{T}_\text{int}}{3 + \delta}, \qquad
\widehat{T}_\text{tr} = \frac{2}{3} \widehat{e}_\text{tr}, \qquad
\widehat{T}_\text{int} = \frac{2}{\delta} \widehat{e}_\text{int},
\\
& \widehat{p}_{ij} = \big( (\widehat{\xi}_i - \widehat{u}_i) (\widehat{\xi}_j - \widehat{u}_j),\, \widehat{f} \big),
\\
&\widehat{q}_i = \widehat{q}_{(\text{tr}) i} + \widehat{q}_{(\text{int}) i},
\qquad \widehat{q}_{(\text{tr}) i}  = \left( \frac{1}{2}(\widehat{\xi}_i - \widehat{u}_i)  |\widehat{\bm{\xi}} - \widehat{\bm{u}}|^2,\, \widehat{f} \right),
\\
&
\widehat{q}_{(\text{int}) i} = \big( (\widehat{\xi}_i - \widehat{u}_i) \widehat{I},\, \widehat{f} \big),
\end{aligned}
\end{align}
where $(\widehat{f}, \widehat{g})$ indicates the inner product of dimensionless
functions $\widehat{f}$ and $\widehat{g}$ of $\widehat{\bm{\xi}}$ and $\widehat{I}$ in
the dimensionless Hilbert space $L^2(\text{d}\widehat{\bm{\xi}}\, \text{d}\widehat{I})$, i.e.,
\begin{align*}
\big( \widehat{f}, \widehat{g} \big) =\int_{\mathbb{R}^{3}\times \mathbb{R}_{+}}
\widehat{f} \widehat{g}\, 
\text{d}\widehat{\bm{\xi}}\, \text{d} \widehat{I} 
\quad \text{ for }\;\; \widehat{f}, \widehat{g} \in L^{2} \big( \text{d}\widehat{\bm{\xi}}\,
\text{d}\widehat{I} \big).
\end{align*}
Note that the same symbol $(\cdot\, ,\, \cdot)$ is used for the inner product
in $L^2(\text{d}\bm{\xi}\, \text{d}I)$ and that in $L^2(\text{d}\widehat{\bm{\xi}}\, \text{d}\widehat{I})$.

Similarly, the dimensionless version of equation~\eqref{BE-1} with~\eqref{col-op} is derived as
\begin{equation} \label{BE-dimless}
\text{Sh} \frac{\partial \widehat{f}}{\partial \widehat{t}} + \widehat{\bm{\xi }} \bm{\cdot}
\frac{\partial \widehat{f}}{\partial \widehat{\bm{x}}}
= \frac{1}{\epsilon}\, \widehat{Q}_{\theta } ( \widehat{f}, \widehat{f} ),
\end{equation}
where
\begin{align}\label{col-op-dimless}
\widehat{Q}_{\theta }(\widehat{f}, \widehat{f}) &=\int_{\left( \mathbb{R}^{3}\times \mathbb{R}_{+}\right)^{3}}
\left( \frac{\widehat{f}^{\prime } \widehat{f}_{\ast }^{\prime }}{( \widehat{I}^{\prime} \widehat{I}_{\ast }^{\prime })
^{\delta /2-1}}-\frac{\widehat{f} \widehat{f}_{\ast }}{( \widehat{I} \widehat{I}_{\ast } ) ^{\delta /2-1}} \right) \,
\notag \\
& \qquad \times 
\widehat{W}_{\theta }(\widehat{\bm{\xi }}, \widehat{\bm{\xi }}_{\ast }, \widehat{I}, \widehat{I}_{\ast} \big|
\widehat{\bm{\xi }}^{\prime }, \widehat{\bm{\xi }}_{\ast }^{\prime}, \widehat{I}^{\prime }, \widehat{I}_{\ast }^{\prime } )\,
\text{d}\widehat{\bm{\xi }}_{\ast } \text{d}\widehat{\bm{\xi }}^{\prime }
\text{d}\widehat{\bm{\xi }}_{\ast }^{\prime } \text{d}\widehat{I}_{\ast}\text{d}\widehat{I}^{\prime }\text{d}\widehat{I}_{\ast }^{\prime },
\end{align}
and
\begin{align*}
\text{Sh} = \frac{L_0}{t_0 \xi_0}, \qquad
\epsilon = \frac{\xi_0}{L_0 n_0 W_{\theta 0}} = \frac{l_0}{L_0}. 
\end{align*}
Here, $\text{Sh}$ is the Strouhal number and $\epsilon$ is the Knudsen number.
Furthermore, with the help of the properties of the Dirac delta function, it follows
from expression~\eqref{Wtheta} and the last of the relations~\eqref{dimless-variables}
that
\begin{align}\label{Wtheta-dimless}
& \widehat{W}_{\theta } (\widehat{\bm{\xi }}, \widehat{\bm{\xi }}_{\ast }, \widehat{I}, \widehat{I}_{\ast} \big|
\widehat{\bm{\xi }}^{\prime }, \widehat{\bm{\xi }}_{\ast }^{\prime}, \widehat{I}^{\prime }, \widehat{I}_{\ast }^{\prime } )
\nonumber \\
& \quad = 4 ( \widehat{I} \widehat{I}_{\ast } ) ^{\delta /2-1}
\frac{\left\vert \widehat{\bm{g}} \right\vert }
{\left\vert \widehat{\bm{g}}^{\prime } \right\vert }
\bm{\delta}_{3} \big( \widehat{\bm{\xi }}+\widehat{\bm{\xi }}_{\ast }-\widehat{\bm{\xi }}^{\prime }-\widehat{\bm{\xi }}_{\ast }^{\prime }\big) 
\notag \\
& \qquad
\times \bm{\delta}_{1}\left( \frac{1}{2}\big( \big| \widehat{\bm{\xi }} \big|^{2}
+\big| \widehat{\bm{\xi }}_{\ast} \big|^{2}
- \big| \widehat{\bm{\xi }}^{\prime } \big|^{2}-\big| \widehat{\bm{\xi }}_{\ast }^{\prime }
\big|^{2} \big)
-\Delta \widehat{I} \right) \widehat{\sigma}_\theta,
\end{align}
with
\begin{align*}
\widehat{\sigma}_\theta = \frac{m^2 \xi_0^5}{W_{\theta 0}} \sigma_\theta, \qquad
\Delta \widehat{I} = \widehat{I}' + \widehat{I}_*' - \widehat{I}' - \widehat{I}_*'.
\end{align*}

Other relations that appeared in Secs.~\ref{sec:B-EQ}--\ref{sec:Fredholm} are
also appropriately nondimensionalized. Here, we only show the results corresponding
to expressions~\eqref{sigmatheta-2} and \eqref{sigma-model}--\eqref{Bs-Br-model}.
The scattering cross section~\eqref{sigmatheta-2} is nondimensionalized as
\begin{align}\label{sigma-dimless}
\widehat{\sigma}_\theta = \theta \widehat{\sigma}_\text{s} + (1 - \theta) \widehat{\sigma}_\text{r}
\bm{\delta}_1 (\Delta \widehat{I}),
\end{align}
with
\begin{align*}
\widehat{\sigma}_\text{s} = \frac{m^2 \xi_0^5}{W_{\theta 0}} \sigma_\text{s}, \qquad
\widehat{\sigma}_\text{r} = \frac{m \xi_0^3}{W_{\theta 0}} \sigma_\text{r}.
\end{align*}
For the models of $\sigma_\text{s}$ and $\sigma_\text{r}$ introduced in~\eqref{sigma-model} and \eqref{Cs-Cr}, the corresponding
$\widehat{\sigma}_\text{s}$ and $\widehat{\sigma}_\text{r}$ become as follows:
\begin{subequations}\label{sigma-model-dimless}
\begin{align}
& \widehat{\sigma}_\text{s}=\widehat{C}_\text{s} \cdot 2^{-( \beta +1)}
\big( \widehat{I}+ \widehat{I}_{\ast }\big) ^{\alpha } \big( \widehat{I}^{\prime }
+ \widehat{I}_{\ast }^{\prime}\big) ^{\alpha }
\dfrac{\big( \widehat{I}^{\prime } \widehat{I}_{\ast }^{\prime }\big)^{\delta /2-1}}
{\widehat{E}^{\delta +\alpha +\left( \beta +1\right) /2}}\left\vert \widehat{\bm{g}}\right\vert^{\beta -1}
\left\vert \widehat{\bm{g}}^{\prime }\right\vert^{\beta +1},
\\
& \widehat{\sigma}_\text{r} = \widehat{C}_\text{r} \cdot \frac{\big( \widehat{I}' \widehat{I}_*' \big)^{\delta/2-1}}
{\big(\widehat{I} + \widehat{I}_*\big)^{\delta -1-\alpha}}
|\widehat{\bm{g}}|^{\beta - 1},
\end{align}
\end{subequations}
where
\begin{align}\label{Cs-Cr-dimless}
\widehat{C}_\text{s} = \frac{m^\alpha \xi_0^{2\alpha+\beta}}{W_{\theta 0}} C_\text{s}, \qquad
\widehat{C}_\text{r} = \frac{m^\alpha \xi_0^{2\alpha+\beta}}{W_{\theta 0}} C_\text{r}, \qquad
\widehat{C}_\text{s}=\frac{\Gamma \big( \delta +\alpha +( \beta +3 ) /2\big) }
{\Gamma \big( ( \beta +3 ) /2 \big) \Gamma \left( \delta +\alpha \right) }
\widehat{C}_\text{r}.
\end{align}
%

\subsubsection{Parameter setting and convention}\label{sec:parameter}

We have derived the dimensionless version of the Boltzmann equation~\eqref{BE-dimless} with a collision operator given by~\eqref{col-op-dimless}, \eqref{Wtheta-dimless},
\eqref{sigma-dimless}, \eqref{sigma-model-dimless}, and \eqref{Cs-Cr-dimless}.
In this paper, we assume that
\begin{align}\label{parameters}
\text{Sh} = 1, \qquad \epsilon \ll 1, \qquad \theta \ll 1.
\end{align}
Here,  $\text{Sh}=1$ corresponds to the so-called fluid time scaling and
$\epsilon \ll 1$ corresponds to the near fluid regime. The assumption
$\theta \ll 1$ indicates that the resonant collisions are dominant, that is,
the relaxation of the internal modes is slow.
In the following subsections, we consider the case
of $\theta \approx \epsilon^2$ and that of $\theta \approx \epsilon$ separately.

Now, let us compare the expression of the dimensional macroscopic
quantities \eqref{macro} and that of the dimensionless ones \eqref{macro-dimless}.
Then, we notice that the relations~\eqref{macro-dimless} is formally obtained from the relations~\eqref{macro} by letting $m=k_\text{B}=1$ (and putting a hat $\,\widehat{}\,$ on
each physical quantity). Similarly, the same operation formally transforms the
dimensional Boltzmann equation~\eqref{BE-1},
\eqref{col-op}, \eqref{Wtheta}, \eqref{sigmatheta-2}, \eqref{sigma-model},
and \eqref{Cs-Cr}, into its dimensionless version~\eqref{BE-dimless}, \eqref{col-op-dimless}, \eqref{Wtheta-dimless},
\eqref{sigma-dimless}, \eqref{sigma-model-dimless}, and \eqref{Cs-Cr-dimless},
if $1/\epsilon$ is put on the right-hand side.

Taking advantage of this fact, we carry out our analysis using the dimensional
Boltzmann equation with $1/\epsilon$, i.e.,
\begin{equation} \label{BE-2}
\frac{\partial f}{\partial t} + \bm{\xi } \bm{\cdot} \frac{\partial f}{\partial \bm{x}}
= \frac{1}{\epsilon}Q_{\theta }\left( f,\, f \right)
= \frac{1}{\epsilon} \left[ \, \theta Q_\text{s} (f, f) + (1 - \theta) Q_\text{r} (f, f)\, \right],
\end{equation}
and the equations and relations for the dimensional variables appeared in
Sec.~\ref{sec:kinetic-model}. However, in the following Secs.~\ref{sec:CE-1}
and \ref{sec:CE-2}, it should be interpreted that $m=k_\text{B}=1$
and all the variables are dimensionless, unless otherwise stated.
In this way, we can omit the cumbersome hats on the dimensionless
quantities and recover the dimensional formulas from the dimensionless
ones immediately by letting $\epsilon=1$.

\subsubsection{Transport equations}\label{sec:transport}

It is obvious from the relations~\eqref{conservation} and \eqref{conservation-theta=0},
and equality~\eqref{sym-relation} that the following relations hold:
\begin{align}\label{cons-Qr-Qs}
\begin{aligned}
& \left( 1, Q_\text{s}(f, f) \right) = \left( \bm{\xi}, Q_\text{s} (f, f) \right)
= \left( m|\bm{\xi}|^2+2I, Q_\text{s}(f, f) \right) = 0,
\\
& \left(1, Q_\text{r}(f, f) \right) = \left( \bm{\xi}, Q_\text{r} (f, f) \right)
= \left( |\bm{\xi}|^2, Q_\text{r}(f, f) \right) = \left( I, Q_\text{r} (f, f) \right) = 0.
\end{aligned}
\end{align}
Let us multiply equation~\eqref{BE-2} by $(m, m\bm{\xi}, (1/2)m |\bm{\xi}|^2, I)$
and integrate with respect to $\bm{\xi}$ and $I$ over $\mathbb{R}^3$
and $\mathbb{R}_+$, respectively. Then,
taking account of the properties \eqref{cons-Qr-Qs}, we obtain
the following transport equations:
\begin{subequations}\label{transport}
\begin{align}
& \frac{\partial \rho}{\partial t} + \frac{\partial}{\partial x_j} (\rho u_j) = 0,
\label{transport-a} \\
& \frac{\partial}{\partial t} (\rho u_i) + \frac{\partial}{\partial x_j} 
(\rho  u_i u_j + p_{ij}) = 0,
\label{transport-b} \\
& \frac{\partial}{\partial t} \left[ \rho \left( \frac{3}{2} \frac{k_\text{B}}{m} T_\text{tr} + \frac{1}{2} |\bm{u}|^2 \right) \right]
+ \frac{\partial}{\partial x_j}  \left[ \rho u_j \left( \frac{3}{2} \frac{k_\text{B}}{m} T_\text{tr} + \frac{1}{2} |\bm{u}|^2 \right)
+ p_{ij} u_i + q_{(\text{tr}) j} \right]
\nonumber \\
& \qquad = \frac{\theta}{\epsilon} \left( \frac{1}{2}m |\bm{\xi}|^2, Q_\text{s} (f, f) \right)  
= - \frac{\theta}{\epsilon} \left( I, Q_\text{s} (f, f) \right),
\label{transport-c} \\
& \frac{\partial}{\partial t} \left( \frac{\delta}{2} \frac{k_\text{B}}{m} \rho T_\text{int} \right)
+ \frac{\partial}{\partial x_j} \left( \frac{\delta}{2} \frac{k_\text{B}}{m} \rho u_j T_\text{int}
+ q_{(\text{int}) j} \right)
= \frac{\theta}{\epsilon} \left( I, Q_\text{s} (f, f) \right),
\label{transport-d}
\end{align}
\end{subequations}
where the macroscopic quantities $\rho$, $u_i$, $T_\text{tr}$, $T_\text{int}$, etc.
are defined by relations~\eqref{macro}. Here and in what follows, the summation convention
(the Einstein convention) is used.
Equations \eqref{transport-a} and \eqref{transport-b}
indicate the mass and momentum conservations, respectively, and equations~\eqref{transport-c}
and \eqref{transport-d} the transport of the translational energy and that of the energy of the internal modes, respectively.

\subsection{Case of $\theta = O(\epsilon^2)$}\label{sec:CE-1}

We first consider the case of $\theta = O(\epsilon^2)$ and let
\begin{align}\label{scale-1}
\theta = \kappa \epsilon^2,
\end{align}
where $\kappa$ is a positive constant \cite{W-CU-51}.
Then, equation~\eqref{BE-2} reads
\begin{align}\label{BE-3}
\frac{\partial f}{\partial t} + \xi_j \frac{\partial f}{\partial x_j}
= \frac{1-\kappa \epsilon^2}{\epsilon}\, Q_\text{r} (f, f)
+ \kappa \epsilon Q_\text{s}(f, f).
\end{align}
%

\subsubsection{Chapman--Enskog expansion and zeroth-order solution}

Let us consider the Chapman--Enskog expansion
\begin{align}\label{CE-0}
f = f^{(0)} + \epsilon f^{(1)} + \epsilon^2 f^{(2)} + \cdots,
\end{align}
and substitute it into equation~\eqref{BE-3}. Then, the $O(1/\epsilon)$ term
gives
\begin{align*}
Q_\text{r} (f^{(0)}, f^{(0)}) = 0,
\end{align*}
so that $f^{(0)}$ is the two-temperature equilibrium distribution $M_\text{r}$ [see~\eqref{Mr}], i.e.,
\begin{align}\label{f0}
f^{(0)} = M_\mathrm{r} & = \frac{n I^{\delta/2-1}}{(2\pi k_\mathrm{B} T_\mathrm{tr}/m)^{3/2}
(k_\mathrm{B} T_\mathrm{int})^{\delta/2}
\Gamma (\delta/2)}
\nonumber \\
& \qquad \qquad \qquad \qquad \times
\exp \left( - \frac{m |\bm{\xi} - \bm{u}|^2}{2k_\mathrm{B} T_\mathrm{tr}}
- \frac{I}{k_\mathrm{B} T_\mathrm{int}}
\right).
\end{align}
This suggests that $\rho$, $\bm{u}$, $T_\text{tr}$, and $T_\text{int}$ are unexpanded.
Therefore, the following conditions are imposed for the higher-order terms
$f^{(1)}$, $f^{(2)}$, $\dots$:
\begin{align}\label{CE-constraint}
\big( 1, f^{(m+1)} \big) = \big( \bm{\xi}, f^{(m+1)} \big) = \big( |\bm{\xi}|^2, f^{(m+1)} \big)
= \big( I, f^{(m+1)} \big) = 0,  \qquad (m=0,\, 1,\, 2,\, \dots).
\end{align}

Letting $f = M_\text{r} + O(\epsilon)$ in $p_{ij}$ and $q_i$ in~\eqref{macro},
we have
\begin{align}\label{pq-0th}
p_{ij} = \frac{k_\text{B}}{m} \rho T_\text{tr} \delta_{ij} + O(\epsilon), \qquad
q_i = 0 + O(\epsilon),
\end{align}
where $\delta_{ij}$ is the Kronecker delta.
Substituting identities~\eqref{pq-0th} into system~\eqref{transport} with scaling~\eqref{scale-1}
and neglecting the terms of $O(\epsilon)$ lead to 
\begin{subequations}\label{Euler-1}
\begin{align}
& \frac{\partial \rho}{\partial t} + \frac{\partial}{\partial x_j} (\rho u_j) = 0,
\label{Euler-1-a} \\
& \frac{\partial}{\partial t} (\rho u_i) + \frac{\partial}{\partial x_j} 
\left( \rho  u_i u_j + \frac{k_\text{B}}{m} \rho T_\text{tr} \delta_{ij} \right) = 0,
\label{Euler-1-b} \\
& \frac{\partial}{\partial t} \left[ \rho \left( \frac{3}{2} \frac{k_\text{B}}{m} T_\text{tr} + \frac{1}{2} |\bm{u}|^2 \right) \right]
+ \frac{\partial}{\partial x_j}  \left[ \rho u_j \left( \frac{5}{2} \frac{k_\text{B}}{m} T_\text{tr} + \frac{1}{2} |\bm{u}|^2 \right)
\right] = 0,
\label{Euler-1-c} \\
& \frac{\partial}{\partial t} ( \rho T_\text{int} )
+ \frac{\partial}{\partial x_j} ( \rho u_j T_\text{int} )
= 0. 
\label{Euler-1-d}
\end{align}
\end{subequations}
Equations~\eqref{Euler-1-a}--\eqref{Euler-1-c}
are the Euler equations for $\rho$, $\bm{u}$, and $T_\text{tr}$, and~\eqref{Euler-1-d}
determines $T_\text{int}$. Note that there is no direct interaction between $T_\text{tr}$
and $T_\text{int}$. Equations \eqref{Euler-1} correspond to the Euler equations
in the case of resonant collisions \cite{W-CU-51,BRS-24}.

\subsubsection{First-order solution}

Equation \eqref{BE-3} then gives the equation containing
the terms of $O(1)$ and higher.
Letting $f^{(1)} = M_\text{r} h$ and recalling~\eqref{Mr}, we can write the equation
in the following form:
\begin{align}\label{EQ-h}
\mathcal{L}_\text{r} h = - \frac{1}{M_\text{r}} \left( \frac{\partial M_\text{r}}{\partial t}
+ \xi_j \frac{\partial M_\text{r}}{\partial x_j} \right)
+ O(\epsilon),
\end{align}
where $\mathcal{L}_\text{r} h = \mathcal{L}_\theta h$ with $\theta=0$ [see~\eqref{dec2}], i.e.,
\begin{align}\label{Lr}
\mathcal{L}_\text{r} h = -2 M_\text{r}^{-1} Q_\text{r} (M_\text{r}, M_\text{r} h).
\end{align}
The derivative terms on the right-hand side of~\eqref{EQ-h} can be calculated
explicitly. Then, the time-derivative terms $\partial \rho/\partial t$,
$\partial \bm{u}/\partial t$, $\partial T_\text{tr}/\partial t$, and $\partial T_\text{int}/\partial t$,
arising from $\partial M_\text{r}/\partial t$, are replaced
by the space derivative terms of the macroscopic quantities with the help of equations~\eqref{Euler-1}.
Thus we obtain the following expression [note that system~\eqref{Euler-1} contains the
error of $O(\epsilon)$]:
\begin{align*}
\frac{1}{M_\text{r}} \left( \frac{\partial M_\text{r}}{\partial t} + \xi_j \frac{\partial M_\text{r}}{\partial x_j}
\right)
& = \frac{m}{2k_\text{B} T_\text{tr}} 
\left( \frac{\partial u_i}{\partial x_j} + \frac{\partial u_j}{\partial x_i} \right) A_{ij} (\bm{c})
+ \frac{1}{T_\text{tr}} \frac{\partial T_\text{tr}}{\partial x_j} B_j (\bm{c})
\nonumber \\
& \quad
+ \frac{1}{T_\text{int}} \frac{\partial T_\text{int}}{\partial x_j}
C_j (\bm{c}, I)
+ O(\epsilon),
\end{align*}
where
\begin{align}\label{def-ABC}
A_{ij} (\bm{c}) =  c_i c_j - \frac{1}{3} |\bm{c}|^2 \delta_{ij}, \quad
B_i (\bm{c}) = c_i \left( \frac{m |\bm{c}|^2}{2k_\text{B} T_\text{tr}} - \frac{5}{2} \right), \quad
C_i (\bm{c}, I) = c_i  \left( \frac{I}{k_\text{B} T_\text{int}} - \frac{\delta}{2} \right),
\end{align}
and $\bm{c}$ (or $c_i$) indicates the peculiar velocity, i.e.,
\begin{align*}
\bm{c} = \bm{\xi} - \bm{u}, \qquad (\text{or} \;\; c_i = \xi_i - u_i ).
\end{align*}
Using the above result in equation~\eqref{EQ-h} and neglecting the terms of $O(\epsilon)$, we obtain
the integral equation for $h$:
\begin{align}\label{EQ-h-2}
\mathcal{L}_\text{r} h
& = - \frac{m}{2k_\text{B} T_\text{tr}} 
\left( \frac{\partial u_i}{\partial x_j} + \frac{\partial u_j}{\partial x_i} \right) A_{ij} (\bm{c})
- \frac{1}{T_\text{tr}} \frac{\partial T_\text{tr}}{\partial x_j} B_j (\bm{c})
- \frac{1}{T_\text{int}} \frac{\partial T_\text{int}}{\partial x_j}
C_j (\bm{c}, I).
\end{align}
It should be noted that both sides of~\eqref{EQ-h-2} are functions of $\bm{\xi}$ and $I$, and $\bm{c}$ is
used just for brevity on the right-hand side.

Note that the following equalities hold:
\begin{align*}
& \int_{\mathbb{R}^3} (1,\, c_i,\, |\bm{c}|^2) \left( c_i c_j - \frac{1}{3} |\bm{c}|^2 \delta_{ij} \right)
e^{-m|\bm{c}|^2/(2k_\text{B} T_\text{tr})} \text{d}\bm{c} = 0,
\\ 
& \int_{\mathbb{R}^3}  (1,\, c_i,\, |\bm{c}|^2)\, c_j \left( \frac{m|\bm{c}|^2}{2k_\text{B} T_\text{tr}} - \frac{5}{2} \right)
e^{-m|\bm{c}|^2/(2k_\text{B} T_\text{tr})} \text{d}\bm{c} = 0,
\\
& \int_{\mathbb{R}^3}  (1,\, |\bm{c}|^2)\, c_j \, e^{-m|\bm{c}|^2/(2k_\text{B} T_\text{tr})} \text{d}\bm{c} = 0,
\\
& \int_{\mathbb{R}_+} I^{\delta/2 -1} \left( \frac{I}{k_\text{B} T_\text{int}} - \frac{\delta}{2} \right) e^{-I/(k_\text{B} T_\text{int})} \text{d}I = 0.
\end{align*}
Thus, because of $\text{d}\bm{\xi} = \text{d}\bm{c}$, it is
obvious that $A_{ij} (\bm{c})$, $B_i (\bm{c})$, and $C_i (\bm{c}, I)$ belong to
$(M_\text{r} \text{ker} \mathcal{L}_\text{r})^\perp$
[cf.~notation~\eqref{perp} with $\theta=0$], i.e.,
\begin{align*}
( \Psi,\, M_\text{r} ) = ( \Psi,\, \bm{\xi} M_\text{r} )
= ( \Psi,\, |\bm{\xi}|^2 M_\text{r} ) = ( \Psi,\, I M_\text{r} )
= 0, \qquad ( \Psi = A_{ij},\; B_i,\; \text{and}\,\; C_i ).
\end{align*}
Therefore, equation~\eqref{EQ-h-2} is solvable due to Corollary \ref{Cor2} or
Remark \ref{solvability} (for $\theta=0$).
If we let
\begin{align*}
h = - \frac{m}{2k_\text{B} T_\text{tr}} 
\left( \frac{\partial u_i}{\partial x_j} + \frac{\partial u_j}{\partial x_i} \right) \widetilde{A}_{ij} (\bm{c}, I)
- \frac{1}{T_\text{tr}} \frac{\partial T_\text{tr}}{\partial x_j} \widetilde{B}_j (\bm{c}, I)
- \frac{1}{T_\text{int}} \frac{\partial T_\text{int}}{\partial x_j} \widetilde{C}_j (\bm{c}, I),
\end{align*}
then we have the integral equations for $\widetilde{A}_{ij}$, $\widetilde{B}_i$, and $\widetilde{C}_i$, i.e.,
\begin{align*}
\mathcal{L}_\text{r} \widetilde{A}_{ij} = A_{ij}, \qquad 
\mathcal{L}_\text{r} \widetilde{B}_i = B_i, \qquad
\mathcal{L}_\text{r} \widetilde{C}_i = C_i.
\end{align*}

\noindent
Since the operator $\mathcal{L}_\text{r}$, in the $\bm{c}$ variable, is isotropic
in the sense of Sec.~A.2.6 in \cite{S-07}, the solutions 
$\widetilde{A}_{ij}$, $\widetilde{B}_i$, and $\widetilde{C}_i$ can be obtained
in the following form, in accordance with the form
of the inhomogeneous terms, as in the case of a monatomic gas
(cf.~Appendix A.2.9 in \cite{S-07} and \cite{DG-94}):
\begin{align*}
\widetilde{A}_{ij} (\bm{c}, I) = A_{ij} (\bm{c}) \mathcal{A} (|\bm{c}|, I), \qquad
\widetilde{B}_i = c_i \mathcal{B} (|\bm{c}|, I), \qquad
\widetilde{C}_i = c_i \mathcal{C} (|\bm{c}|, I),
\end{align*}
where $\mathcal{A} (|\bm{c}|, I)$, $\mathcal{B} (|\bm{c}|, I)$, and $\mathcal{C} (|\bm{c}|, I)$
are functions of $|\bm{c}|$ and $I$.

In summary, the solution $h$ is obtained in the following form:
\begin{align}\label{sol-h}
h & = - \frac{m}{2k_\text{B} T_\text{tr}} 
\left( \frac{\partial u_i}{\partial x_j} + \frac{\partial u_j}{\partial x_i} \right) A_{ij} (\bm{c}) \mathcal{A} (|\bm{c}|, I)
- \frac{1}{T_\text{tr}} \frac{\partial T_\text{tr}}{\partial x_j} c_j \mathcal{B} (|\bm{c}|, I)
\nonumber \\
& \;\;\;\;
- \frac{1}{T_\text{int}} \frac{\partial T_\text{int}}{\partial x_j} c_j \mathcal{C} (|\bm{c}|, I),
\end{align}
where $\mathcal{A} (|\bm{c}|, I)$, $\mathcal{B} (|\bm{c}|, I)$, and $\mathcal{C} (|\bm{c}|, I)$ are, respectively,
the solutions of the following equations:
\begin{align}\label{EQ-decomp}
\mathcal{L}_\text{r} \big( A_{ij} (\bm{c}) \mathcal{A} (|\bm{c}|, I) \big) = A_{ij} (\bm{c}),
\quad \,
\mathcal{L}_\text{r} \big( c_i \mathcal{B} (|\bm{c}|, I) \big) = B_i (\bm{c}),
\quad \,
\mathcal{L}_\text{r} \big( c_i \mathcal{C} (|\bm{c}|, I) \big) = C_i (\bm{c}, I).
\end{align}

\noindent
Here, it is recalled that $f^{(1)} = M_\text{r} h$ should satisfy the constraints
\eqref{CE-constraint}. It is obvious that the first term [the term
containing $\mathcal{A} (\bm{c}, I)$] on the right-hand
side of equation~\eqref{sol-h} satisfies~\eqref{CE-constraint}. In order for
the other terms to satisfy the constraints~\eqref{CE-constraint}, the following conditions
should be imposed on $\mathcal{B} (|\bm{c}|, I)$ and
$\mathcal{C} (|\bm{c}|, I)$:
\begin{align*}
\big( c_i,\, c_j \mathcal{B} (|\bm{c}|, I) M_\text{r} \big) = 0, \qquad
\big( c_i,\, c_j \mathcal{C} (|\bm{c}|, I) M_\text{r} \big) = 0, 
\end{align*}
or, with $c=|\bm{c}|$,
\begin{align}\label{sub-cond}
\int_{\mathbb{R}_+ \times \mathbb{R}_+} c^4 \left(
\begin{array}{l}
\mathcal{B} (c, I) \\ \mathcal{C} (c, I)
\end{array}
\right) I^{\delta/2-1}
\exp \left( -\frac{m c^2}{2k_\text{B} T_\text{tr}} - \frac{I}{k_\text{B} T_\text{int}} \right)
\text{d}c\, \text{d}I = 0.
\end{align}
%

\subsubsection{Constitutive laws at Navier--Stokes level}

Now we have the solution up to the first order in $\epsilon$, i.e.,
$f = M_\text{r} ( 1 + \epsilon h ) + O(\epsilon^2)$. The stress tensor $p_{ij}$
up to the corresponding order can be obtained by substituting this $f$ into
$p_{ij}$ in relations~\eqref{macro}. That is,
\begin{align}\label{pij}
p_{ij} & = m \big( c_i c_j, M_\text{r} (1 + \epsilon h) \big) + O(\epsilon^2)
\nonumber \\
& = \frac{k_\text{B}}{m} \rho T_\text{tr} \delta_{ij}
- \epsilon\, m \frac{m}{2k_\text{B} T_\text{tr}}
\left[ \int_{\mathbb{R}^3 \times \mathbb{R}_+}
c_i c_j \left( c_k c_l - \frac{1}{3} |\bm{c}|^2 \delta_{kl} \right) \mathcal{A} (|\bm{c}|, I) M_\text{r}\,
\text{d}\bm{c}\, \text{d}I \right]
\nonumber \\
& \hspace*{75mm}
\times
\left( \frac{\partial u_k}{\partial x_l} + \frac{\partial u_l}{\partial x_k} \right)
+ O(\epsilon^2)
\nonumber \\
& = \frac{k_\text{B}}{m} \rho T_\text{tr} \delta_{ij} - \epsilon \Lambda_\mu (\rho, T_\text{tr}, T_\text{int})
\left( \frac{\partial u_i}{\partial x_j} + \frac{\partial u_j}{\partial x_i}
- \frac{2}{3} \frac{\partial u_k}{\partial x_k} \delta_{ij} \right)
+ O(\epsilon^2),
\end{align}
where we have let
\begin{align}\label{Lambda-1}
\Lambda_\mu (\rho, T_\text{tr}, T_\text{int}) = & \frac{8}{15\sqrt{\pi}}
\left( \frac{m}{2k_\text{B} T_\text{tr}} \right)^{5/2} \frac{\rho}{(k_\text{B} T_\text{int})^{\delta/2} \Gamma (\delta/2)}
\nonumber \\
& \!\! \times
\int_0^\infty \left[ \int_0^\infty c^6 \mathcal{A}(c, I) \exp\left( -\frac{mc^2}{2k_\text{B} T_\text{tr}} \right) \text{d}c \right]
I^{\delta/2 - 1} \exp \left( - \frac{I}{k_\text{B} T_\text{int}} \right) \text{d}I.
\end{align}

Similarly, the heat-flow vector up to $O(\epsilon)$ can be obtained from $q_{i}$ in relations~\eqref{macro}, that is,
\begin{align*}
q_i & = q_{(\text{tr}) i} + q_{\text{(int}) i},
\end{align*}
and
\begin{subequations}\label{qi}
\begin{align}
q_{(\text{tr}) i} & = \frac{m}{2} \left( c_i |\bm{c}|^2, M_\text{r} (1 + \epsilon h) \right) + O(\epsilon^2) 
\nonumber \\
& = - \epsilon \frac{m}{2} \left[ \int_{\mathbb{R}^3 \times \mathbb{R}_+}
c_i c_j |\bm{c}|^2 \mathcal{B} (|\bm{c}|, I) M_\text{r} d\bm{c}dI \cdot \frac{1}{T_\text{tr}} \frac{\partial T_\text{tr}}{\partial x_j} \right.
\nonumber \\
& \qquad \qquad \qquad \left.
+ \int_{\mathbb{R}^3 \times \mathbb{R}_+}
c_i c_j |\bm{c}|^2 \mathcal{C} (|\bm{c}|, I) M_\text{r} d\bm{c}dI \cdot \frac{1}{T_\text{int}} \frac{\partial T_\text{int}}{\partial x_j} \right]
 + O(\epsilon^2)\nonumber 
\\
& = - \epsilon \Lambda_\text{tr}^\text{tr} (\rho, T_\text{tr}, T_\text{int}) \frac{1}{T_\text{tr}} \frac{\partial T_\text{tr}}{\partial x_i}
- \epsilon \Lambda_\text{int}^\text{tr} (\rho, T_\text{tr}, T_\text{int}) \frac{1}{T_\text{int}} \frac{\partial T_\text{int}}{\partial x_i}
 + O(\epsilon^2),
\label{qi-a} \\[2mm]
q_{(\text{int}) i} & = \left( c_i I, M_\text{r} (1 + \epsilon h) \right) + O(\epsilon^2) 
\nonumber \\
& = - \epsilon \left[ \int_{\mathbb{R}^3 \times \mathbb{R}_+}
c_i c_j I\, \mathcal{B} (|\bm{c}|, I) M_\text{r} d\bm{c}dI \cdot \frac{1}{T_\text{tr}} \frac{\partial T_\text{tr}}{\partial x_j} \right.
\nonumber \\
& \qquad \qquad \qquad \left.
+ \int_{\mathbb{R}^3 \times \mathbb{R}_+}
c_i c_j I\, \mathcal{C} (|\bm{c}|, I) M_\text{r} d\bm{c}dI \cdot \frac{1}{T_\text{int}} \frac{\partial T_\text{int}}{\partial x_j} \right]
 + O(\epsilon^2)\nonumber 
\\
& = - \epsilon \Lambda_\text{tr}^\text{int} (\rho, T_\text{tr}, T_\text{int}) \frac{1}{T_\text{tr}} \frac{\partial T_\text{tr}}{\partial x_i}
- \epsilon \Lambda_\text{int}^\text{int} (\rho, T_\text{tr}, T_\text{int}) \frac{1}{T_\text{int}} \frac{\partial T_\text{int}}{\partial x_i}
 + O(\epsilon^2),
 \label{qi-b}
\end{align}
\end{subequations}
where we have let
\begin{subequations}\label{Lambda-2}
\begin{align}
\left[
\begin{array}{l}
\Lambda_\text{tr}^\text{tr} (\rho, T_\text{tr}, T_\text{int}) \\[2mm]
\Lambda_\text{int}^\text{tr} (\rho, T_\text{tr}, T_\text{int})
\end{array} \right]
& = \frac{2}{3\sqrt{\pi}} \left( \frac{m}{2 k_B T_\text{tr}} \right)^{3/2}
\frac{\rho}{(k_B T_\text{int})^{\delta/2} \Gamma (\delta/2)}
\nonumber \\
& \quad  \times
\int_0^\infty \left\{ \int_0^\infty c^6
\left[ \begin{array}{l}
\mathcal{B} (c, I)\\[2mm]
\mathcal{C} (c, I)
\end{array} \right]
 \exp\left( -\frac{mc^2}{2k_B T_\text{tr}} \right) \text{d}c \right\}
 \nonumber \\
& \qquad  \qquad  \qquad \qquad \qquad  \qquad \qquad \times
I^{\delta/2 - 1} \exp \left( - \frac{I}{k_B T_\text{int}} \right) \text{d}I,
\label{Lambda-2-a} \\
\left[
\begin{array}{l}
\Lambda_\text{tr}^\text{int} (\rho, T_\text{tr}, T_\text{int}) \\[2mm]
\Lambda_\text{int}^\text{int} (\rho, T_\text{tr}, T_\text{int})
\end{array}
\right]
& = \frac{4}{3\sqrt{\pi}} \frac{1}{m} \left( \frac{m}{2 k_B T_\text{tr}} \right)^{3/2}
\frac{\rho}{(k_B T_\text{int})^{\delta/2} \Gamma (\delta/2)}
\nonumber \\
& \quad  \times
\int_0^\infty \left\{ \int_0^\infty c^4 
\left[
\begin{array}{l}
\mathcal{B} (c, I) \\[2mm]
\mathcal{C} (c, I)
\end{array}
\right]
\exp\left( -\frac{mc^2}{2k_B T_\text{tr}} \right) \text{d}c \right\}
 \nonumber \\
& \qquad  \qquad  \qquad \qquad \qquad  \qquad \qquad \times
I^{\delta/2} \exp \left( - \frac{I}{k_B T_\text{int}} \right) \text{d}I.
\label{Lambda-2-b}
\end{align}
\end{subequations}

It should be noted that the bulk viscosity does not occur in the stress
tensor $p_{ij}$ up to the order of $\epsilon$. In addition, both heat-flow
vectors $q_{\text{(tr)}i}$ and $q_{\text{(int)}i}$ contain terms proportional
to $-\partial T_\text{tr}/\partial x_i$ and $-\partial T_\text{int}/\partial x_i$.
Thus, they show the effect of cross diffusion.

\subsubsection{Source term and two-temperature Navier--Stokes equations}

Now, let us consider the source term (i.e., the right-hand side) of equation~\eqref{transport-d},
which is also the source term in equation~\eqref{transport-c}.
Recalling the scaling~\eqref{scale-1} and using the expansion \eqref{CE-0}, it
can be written as
\begin{align}\label{Qs-f-f}
\frac{\theta}{\epsilon} \big( I, Q_\text{s} (f, f) \big)
& = \epsilon \kappa \, \big( I, Q_\text{s} (f^{(0)}, f^{(0)}) \big) + O(\epsilon^2),
\nonumber \\
&  \;= \epsilon \kappa \, \big( I, Q_\text{s} (M_\text{r}, M_\text{r}) \big) + O(\epsilon^2).
\end{align}
For the collision kernel given by \eqref{sigma-model-a},
the term $\big( I, Q_\text{s} (M_\text{r}, M_\text{r}) \big)$
can be calculated explicitly, as shown in Appendix \ref{sec:source}, and is reduced to
the following form [see~\eqref{source-term-2}]:
\begin{align}\label{Qs-Mr-Mr}
\big( I, Q_\text{s} (M_\text{r}, M_\text{r}) \big) = \mathcal{F} (\rho, T_\text{tr}, T_\text{int})
(T_\text{tr} - T_\text{int}),
\end{align}
where $\mathcal{F} (\rho, T_\text{tr}, T_\text{int})$ is a function of $\rho$,
$T_\text{tr}$, and $T_\text{int}$ given by~\eqref{source-term-aux}, i.e.,
\begin{align*}
\mathcal{F} (\rho, T_\text{tr}, T_\text{int})
= C\, \frac{k_\text{B}^{\alpha + 1 +\beta/2}}{m^{2 + \beta/2}}\,
\rho^2 T_\text{tr}^{\beta/2} T_\text{int}^\alpha,
\end{align*}
with
\begin{align*}
C = 2^{\beta +2} \sqrt{\pi }\,
\frac{\Gamma \left( \delta +\alpha +1 \right)
\Gamma^2 \left( \delta /2 \right) \Gamma \big( \left( \beta+5\right) /2\big)}
{\left[ \delta +\alpha +\left( \beta +3 \right)/2 \right]
\Gamma^{2} \left( \delta \right) }\, C_\text{r}.
\end{align*}

Substituting expressions~\eqref{pij} with \eqref{Lambda-1}, \eqref{qi} with \eqref{Lambda-2},
and \eqref{Qs-f-f} with \eqref{Qs-Mr-Mr}, into the system~\eqref{transport}
and neglecting the terms of $O(\epsilon^2)$, we have the following
equations:
\begin{subequations}\label{NS}
\begin{align}
& \frac{\partial \rho}{\partial t} + \frac{\partial}{\partial x_j} (\rho u_j) = 0,
\label{NS-a} \\
& \frac{\partial}{\partial t} (\rho u_i) + \frac{\partial}{\partial x_j} 
(\rho  u_i u_j) + \frac{k_B}{m} \frac{\partial}{\partial x_i} (\rho T_\text{tr})
\nonumber \\
& \qquad \qquad
- \epsilon \frac{\partial}{\partial x_j} 
\left[ \Lambda_\mu (\rho, T_\text{tr}, T_\text{int})
\left( \frac{\partial u_i}{\partial x_j} + \frac{\partial u_j}{\partial x_i}
- \frac{2}{3} \frac{\partial u_k}{\partial x_k} \delta_{ij} \right) \right]
 = 0,
\label{NS-b} \\
& \frac{\partial}{\partial t} \left[ \rho \left( \frac{3}{2} \frac{k_B}{m} T_\text{tr} + \frac{1}{2} |\bm{u}|^2 \right) \right]
+ \frac{\partial}{\partial x_j}  \left[ \rho u_j \left( \frac{5}{2} \frac{k_B}{m} T_\text{tr} + \frac{1}{2} |\bm{u}|^2 \right) \right]
\nonumber \\
& \qquad \qquad
- \epsilon \frac{\partial}{\partial x_j} \left[ u_i \Lambda_\mu (\rho, T_\text{tr}, T_\text{int})
\left( \frac{\partial u_i}{\partial x_j} + \frac{\partial u_j}{\partial x_i}
- \frac{2}{3} \frac{\partial u_k}{\partial x_k} \delta_{ij} \right) \right]
\nonumber \\
& \qquad \qquad
- \epsilon \frac{\partial}{\partial x_j} \left[
 \Lambda_\text{tr}^\text{tr} (\rho, T_\text{tr}, T_\text{int}) \frac{1}{T_\text{tr}} \frac{\partial T_\text{tr}}{\partial x_j}
+ \Lambda_\text{int}^\text{tr} (\rho, T_\text{tr}, T_\text{int}) \frac{1}{T_\text{int}} \frac{\partial T_\text{int}}{\partial x_j}
 \right]
\nonumber \\
& \quad = - \epsilon \kappa  \mathcal{F} (\rho, T_\text{tr}, T_\text{int}) (T_\text{tr} - T_\text{int}),
\label{NS-c} \\[2mm]
& \frac{\partial}{\partial t} \left( \frac{\delta}{2} \frac{k_B}{m} \rho T_\text{int} \right)
+ \frac{\partial}{\partial x_j} \left( \frac{\delta}{2} \frac{k_B}{m} \rho u_j T_\text{int} \right)
\nonumber \\
& \qquad \qquad
- \epsilon \frac{\partial}{\partial x_j} \left[
 \Lambda_\text{tr}^\text{int} (\rho, T_\text{tr}, T_\text{int}) \frac{1}{T_\text{tr}} \frac{\partial T_\text{tr}}{\partial x_j}
+ \Lambda_\text{int}^\text{int} (\rho, T_\text{tr}, T_\text{int}) \frac{1}{T_\text{int}} \frac{\partial T_\text{int}}{\partial x_j}
 \right]
 \nonumber \\
& \quad = \epsilon \kappa  \mathcal{F} (\rho, T_\text{tr}, T_\text{int}) (T_\text{tr} - T_\text{int}). 
\label{NS-d}
\end{align}
\end{subequations}
The system \eqref{NS} is the system of Navier--Stokes-type equations for two
temperatures and with relaxation terms.
Note that the viscous-stress terms, the heat-conduction terms, and
the relaxation terms are all of the order of $\epsilon$ for the scaling \eqref{scale-1},
unlike the system~\eqref{NS-2} that will appear for the scaling \eqref{scale-2} (Sec.~\ref{sec:two-temp-NS}).
One can readily show that the transport coefficients $\Lambda_\mu$
in equations~\eqref{NS-b} and \eqref{NS-c}, $\Lambda_\text{tr}^\text{tr}$ in~\eqref{NS-c}, and
$\Lambda_\text{int}^\text{int}$ in~\eqref{NS-d} are positive (see Appendix \ref{sec:positivity}).

Since the solutions $\mathcal{A} (|\bm{c}|, I)$, $\mathcal{B} (|\bm{c}|, I)$, and $\mathcal{C} (|\bm{c}|, I)$
to equations \eqref{EQ-decomp} are not obtained explicitly, system~\eqref{NS-2} is not completely explicit
in this sense. However, it is not difficult to obtain these solutions either numerically or approximately.
In addition, the coefficient $\mathcal{F}$ of the relaxation terms is explicit in terms of the
parameters included in the collision model \eqref{sigma-model}. Therefore,
we can claim that~\eqref{NS} is a system constructed explicitly.

\begin{remark}\label{rem:2}
Equations essentially similar to the system~\eqref{NS} were derived in a more
abstract form in \cite{W-CU-51} using
a different Boltzmann model with a single discrete energy variable under
the assumption that the difference $|T_\mathrm{tr} - T_\mathrm{int}|$ is small.
It should be emphasized that the assumption of smallness of $|T_\mathrm{tr} - T_\mathrm{int}|$
is not necessary to derive~\eqref{NS} here.
\end{remark}

\begin{remark}
Adding the factor $E^{\vartheta }$ with a constant $\vartheta$ to the scattering cross section
$\sigma_\mathrm{s}$~\eqref{sigma-model-a} [and correspondingly to $\sigma_\mathrm{r}$~\eqref{sigma-model-b}] makes the term $\left(
I,Q_\mathrm{s}(M_\mathrm{r}, M_\mathrm{r})\right) $ again of the form \eqref{Qs-Mr-Mr}.
However, $\mathcal{F}(\rho ,T_\mathrm{tr},T_\mathrm{int})$ is given only implicitly in
this case, as the integral corresponding to $\Omega$ in~\eqref{Omega} cannot be explicitly calculated. This is due to the fact that the mixed factor
$\left( \overline{q}T_\mathrm{tr}+\overline{v}T_\mathrm{int}\right) ^{\vartheta }$ appears
in the integral corresponding to the first line of~\eqref{Omega-q-v-bar}
and thus the integral with respect to $\overline{q}$ and that with respect to $\overline{v}$
are not separable. However, if $\vartheta $ is a nonnegative integer, then
one obtains a sum of such separable integrals, which can be calculated explicitly
as~\eqref{Omega-q-v-bar}.
In this case, $\mathcal{F}(\rho ,T_\mathrm{tr},T_\mathrm{int})$ is obtained explicitly.
\end{remark}

\begin{remark}
We started our analysis with equation~\eqref{BE-2},
which can be interpreted as both a dimensional and a dimensionless equation
(cf.~Sec.~\ref{sec:parameter}), Consequently, system~\eqref{NS} can also be
interpreted in both ways. To interpret equations~\eqref{NS} as dimensionless,
we need to set $m=k_\mathrm{B}=1$ and interpret all the independent and dependent variables, as well as the collision operators, as dimensionless, i.e., as the variables and collision operators with a hat $\;\widehat{}\;$
defined in Sec.~\ref{sec:dimless}.
In fact, the parameter setting \eqref{parameters} makes sense only for the dimensionless equations. On the other hand, to interpret~\eqref{NS}
as dimensional equations, we just need to let $\epsilon=1$ and $\kappa = \theta$. 
The same remark also applies to the equations in Sec.~\ref{sec:CE-2}.
\end{remark}

\subsubsection{Particular cases}

In the following, we will further investigate the transport coefficients in the system~\eqref{NS} using a collision kernel given by~\eqref{sigma-model-b} explicitly.
From identities~\eqref{Qr-BL}, \eqref{sigma-model-b}, and \eqref{Lr}, it follows that
\begin{align*}
\mathcal{L}_\text{r} h
= -C_\text{r} \int_{[0, 1] \times \mathbb{S}^2 \times \mathbb{R}^3 \times \mathbb{R}_+}
M_{\text{r}*} (h_*' + h' - h_* - h) \frac{ (I' I_*')^{\delta/2-1}}{(I + I_*)^{\delta-2-\alpha}}
|\bm{g}|^\beta \text{d}r \text{d}\bm{\sigma} \text{d}\bm{\xi}_* \text{d}I_*.
\end{align*}
Recalling that $\bm{c}=\bm{\xi}-\bm{u}$, let us put $\bm{c}_*=\bm{\xi}_*-\bm{u}$,
$\bm{c}'=\bm{\xi}'-\bm{u}$, and $\bm{c}_*'=\bm{\xi}_*'-\bm{u}$. Then, we have
\begin{align*}
\bm{g} = \bm{c} - \bm{c}_*, \qquad
\bm{c}' = \frac{\bm{c} + \bm{c}_*}{2} + \frac{|\bm{g}|}{2} \bm{\sigma}, \qquad
\bm{c}_*' = \frac{\bm{c} + \bm{c}_*}{2} - \frac{|\bm{g}|}{2} \bm{\sigma},
\end{align*}
and the above $\mathcal{L}_\text{r} h$ is transformed, using the
relation $I'I_*'=r(1-r) (I+I_*)^2$, into the following form:
\begin{align}\label{Lr-h}
\mathcal{L}_\text{r} h
= & -C_\text{r} \int_{[0, 1] \times \mathbb{S}^2 \times \mathbb{R}^3 \times \mathbb{R}_+}
M_{\text{r}*} (h_*' + h' - h_* - h) \, \nonumber \\ 
& \qquad \qquad \qquad \qquad \qquad \times [r (1-r)]^{\delta/2-1} (I + I_*)^\alpha
|\bm{g}|^\beta \text{d}r \text{d}\bm{\sigma} \text{d}\bm{c}_* \text{d}I_*.
\end{align}
Here, $M_\text{r}$ and $h$ are regarded as functions of $\bm{c}$
and $I$ rather than $\bm{\xi}$ and $I$ (the dependence on $t$ and $\bm{x}$,
if any, is omitted), and the conventional notation
$h_*=h(\bm{c}_*, I_*)$, $h' = h(\bm{c}', I')$, etc.~is used.
In the following,
the change of variables from $(\bm{\xi}, \bm{\xi}_*, \bm{\xi}', \bm{\xi}_*')$ to
$(\bm{c}, \bm{c}_*, \bm{c}', \bm{c}_*')$ is occasionally made, and the
corresponding notation, such as $h_*=h(\bm{c}_*, I_*)$, $h' = h(\bm{c}', I')$, is
used without any notice.

Now, we focus on the special case where $\alpha=0$. Then, equation~\eqref{Lr-h} reduces
to
\begin{align}\label{Lr-h-alpha=0}
\mathcal{L}_\text{r} h
= -C_\text{r} \int_{[0, 1] \times \mathbb{S}^2 \times \mathbb{R}^3 \times \mathbb{R}_+}
M_{\text{r}*} (h_*' + h' - h_* - h) \, [r (1-r)]^{\delta/2-1}
|\bm{g}|^\beta \text{d}r \text{d}\bm{\sigma} \text{d}\bm{c}_* \text{d}I_*.
\end{align}
If $h$ is a function of $\bm{c}$ only and does not depend on $I$, then equation~\eqref{Lr-h-alpha=0} is reduced to the following form (see Appendix \ref{sec:derivation}):
\begin{align}\label{Lr-h-alpha=0-2}
\mathcal{L}_\text{r} h
= - C_\text{r} \frac{\sqrt{m} \rho}{(2\pi k_\text{B} T_\text{tr})^{3/2}}
\frac{\Gamma^2 (\delta/2)}{\Gamma (\delta)}
\int_{\mathbb{R}^3 \times \mathbb{S}^2}
\! \! \exp \left( - \frac{m |\bm{c}_*|^2}{2k_\text{B} T_\text{tr}} \right)
|\bm{g}|^\beta (h_*' + h' - h_* - h) \text{d}\bm{c}_* \text{d}\bm{\sigma}.
\end{align}
That is, $\mathcal{L}_\text{r} h$ is also independent of $I$. Therefore,
noting that $A_{ij} (\bm{c})$
and $B_i (\bm{c})$ in~\eqref{EQ-decomp} are independent of $I$, we can consistently
assume that the functions
$\mathcal{A} (|\bm{c}|, I)$ and $\mathcal{B} (|\bm{c}|, I)$ are independent of $I$, namely,
\begin{align}\label{A0-B0}
\mathcal{A} (|\bm{c}|, I) = \mathcal{A}_0 (|\bm{c}|), \qquad
\mathcal{B} (|\bm{c}|, I) = \mathcal{B}_0 (|\bm{c}|).
\end{align}
On the other hand, if $h$ is of the form
$h = \left[ I/(k_\text{B} T_\text{int}) - \delta/2 \right]\widetilde{h}(\bm{c})$,
with $\widetilde{h}(\bm{c})$ being independent of $I$, then~\eqref{Lr-h-alpha=0}
is transformed into the following form (see Appendix \ref{sec:derivation}):
\begin{align}\label{Lr-h-alpha=0-3}
\mathcal{L}_\text{r} h
& =  - C_\text{r} \frac{\sqrt{m} \rho}{(2\pi k_\text{B} T_\text{tr})^{3/2}}
\frac{\Gamma^2 (\delta/2)}{\Gamma (\delta)}
\left( \frac{I}{k_\text{B} T_\text{int}} - \frac{\delta}{2} \right) 
\nonumber \\
& \qquad \qquad \qquad \qquad \times
\int_{\mathbb{R}^3 \times \mathbb{S}^2}
\exp \left( - \frac{m |\bm{c}_*|^2}{2k_\text{B} T_\text{tr}} \right)
|\bm{g}|^\beta (\widetilde{h}' - \widetilde{h}) \text{d}\bm{c}_* \text{d}\bm{\sigma},
\end{align}
which is also of the form
$\left[ I/(k_\text{B} T_\text{int}) - \delta/2 \right] \times (\text{function of } \bm{c}$).
Therefore, since $C_i (\bm{c}, I)$ in~\eqref{EQ-decomp} is of this form,
we can consistently assume that $\mathcal{C} (|\bm{c}|, I)$ is of the form
\begin{align}\label{C0}
\mathcal{C} (|\bm{c}|, I) = \left( \frac{I}{k_\text{B} T_\text{int}} - \frac{\delta}{2} \right)
\mathcal{C}_0 (|\bm{c}|).
\end{align}

Using expressions~\eqref{A0-B0} and \eqref{C0} in constraints~\eqref{sub-cond},
one finds that the second line of~\eqref{sub-cond} is automatically
satisfied and that the following condition needs to be imposed on
$\mathcal{B}_0 (|\bm{c}|)$:
\begin{align}\label{sub-cond-2}
\int_0^\infty c^4 \mathcal{B}_0 (c) \exp \left( -\frac{m c^2}{2k_\text{B} T_\text{tr}} \right)
\text{d}c= 0.
\end{align}
Let $\widetilde{\Lambda}_\mu$, $\widetilde{\Lambda}_\text{tr}^\text{tr}$, etc. denote
$\Lambda_\mu$, $\Lambda_\text{tr}^\text{tr}$, etc. in expressions~\eqref{Lambda-1}
and \eqref{Lambda-2} for the collision model \eqref{sigma-model-b} with $\alpha=0$,
i.e.,
\begin{align*}
\big( \widetilde{\Lambda}_\mu,\, \widetilde{\Lambda}_\text{tr}^\text{tr},\,
\widetilde{\Lambda}_\text{int}^\text{tr},\, \widetilde{\Lambda}_\text{tr}^\text{int},\,
\widetilde{\Lambda}_\text{int}^\text{int} \big)
= \big( \Lambda_\mu,\, \Lambda_\text{tr}^\text{tr},\,
\Lambda_\text{int}^\text{tr},\, \Lambda_\text{tr}^\text{int},\,
\Lambda_\text{int}^\text{int} \big) \;\,
\mbox{[for~\eqref{sigma-model-b} with $\alpha=0$]}.
\end{align*}
The substitution of the first identity of~\eqref{A0-B0} into~\eqref{Lambda-1}
gives
\begin{align}\label{Lambda-1-alpha=0}
\widetilde{\Lambda}_\mu (\rho, T_\text{tr}, T_\text{int}) = & \frac{8}{15\sqrt{\pi}} \rho
\left( \frac{m}{2k_\text{B} T_\text{tr}} \right)^{5/2} 
\int_0^\infty c^6 \mathcal{A}_0 (c) \exp\left( -\frac{mc^2}{2k_\text{B} T_\text{tr}} \right) \text{d}c,
\end{align}
and the substitution of the second identity of~\eqref{A0-B0} and identity~\eqref{C0}
into~\eqref{Lambda-2} gives
\begin{subequations}\label{Lambda-2-alpha=0}
\begin{align}
& \widetilde{\Lambda}_\text{tr}^\text{tr} (\rho, T_\text{tr}, T_\text{int})
= \frac{2}{3\sqrt{\pi}} \rho \left( \frac{m}{2 k_B T_\text{tr}} \right)^{3/2}
\int_0^\infty c^6 \mathcal{B}_0 (c) \exp\left( -\frac{mc^2}{2k_B T_\text{tr}} \right) \text{d}c,
\label{Lambda-2-alpha=0-a} \\
& \widetilde{\Lambda}_\text{int}^\text{int} (\rho, T_\text{tr}, T_\text{int})
= \frac{2}{3\sqrt{\pi}} \delta \frac{\rho}{m} \left( \frac{m}{2 k_\text{B} T_\text{tr}} \right)^{3/2}
k_\text{B} T_\text{int}
\int_0^\infty c^4 \mathcal{C}_0 (c) \exp\left( -\frac{mc^2}{2k_B T_\text{tr}} \right) \text{d}c,
\label{Lambda-2-alpha=0-b} \\
& \widetilde{\Lambda}_\text{int}^\text{tr} (\rho, T_\text{tr}, T_\text{int}) = 0,
\label{Lambda-2-alpha=0-c} \\
& \widetilde{\Lambda}_\text{tr}^\text{int} (\rho, T_\text{tr}, T_\text{int}) = 0.
\label{Lambda-2-alpha=0-d}
\end{align}
\end{subequations}
Identity \eqref{Lambda-2-alpha=0-d} is obvious from equality~\eqref{sub-cond-2}.

In summary, for collision models \eqref{sigma-model} with $\alpha=0$,
the two-temperature Navier--Stokes model \eqref{NS} reduces to the
following system:
\begin{subequations}\label{NS-alpha=0}
\begin{align}
& \frac{\partial \rho}{\partial t} + \frac{\partial}{\partial x_j} (\rho u_j) = 0,
\label{NS-alpha=0-a} \\
& \frac{\partial}{\partial t} (\rho u_i) + \frac{\partial}{\partial x_j} 
(\rho  u_i u_j) + \frac{k_B}{m} \frac{\partial}{\partial x_i} (\rho T_\text{tr})
\nonumber \\
& \qquad \qquad
- \epsilon \frac{\partial}{\partial x_j} 
\left[ \widetilde{\Lambda}_\mu (\rho, T_\text{tr}, T_\text{int})
\left( \frac{\partial u_i}{\partial x_j} + \frac{\partial u_j}{\partial x_i}
- \frac{2}{3} \frac{\partial u_k}{\partial x_k} \delta_{ij} \right) \right]
 = 0,
\label{NS-alpha=0-b} \\
& \frac{\partial}{\partial t} \left[ \rho \left( \frac{3}{2} \frac{k_B}{m} T_\text{tr} + \frac{1}{2} |\bm{u}|^2 \right) \right]
+ \frac{\partial}{\partial x_j}  \left[ \rho u_j \left( \frac{5}{2} \frac{k_B}{m} T_\text{tr} + \frac{1}{2} |\bm{u}|^2 \right) \right]
\nonumber \\
& \qquad \qquad
- \epsilon \frac{\partial}{\partial x_j} \left[ u_i \widetilde{\Lambda}_\mu (\rho, T_\text{tr}, T_\text{int})
\left( \frac{\partial u_i}{\partial x_j} + \frac{\partial u_j}{\partial x_i}
- \frac{2}{3} \frac{\partial u_k}{\partial x_k} \delta_{ij} \right) \right]
\nonumber \\
& \qquad \qquad
- \epsilon \frac{\partial}{\partial x_j} \left[
 \widetilde{\Lambda}_\text{tr}^\text{tr} (\rho, T_\text{tr}, T_\text{int}) \frac{1}{T_\text{tr}} \frac{\partial T_\text{tr}}{\partial x_j}
 \right]
\nonumber \\
& \quad = - \epsilon \kappa  \widetilde{\mathcal{F}} (\rho, T_\text{tr}) (T_\text{tr} - T_\text{int}),
\label{NS-alpha=0-c} \\[2mm]
& \frac{\partial}{\partial t} \left( \frac{\delta}{2} \frac{k_B}{m} \rho T_\text{int} \right)
+ \frac{\partial}{\partial x_j} \left( \frac{\delta}{2} \frac{k_B}{m} \rho u_j T_\text{int} \right)
- \epsilon \frac{\partial}{\partial x_j} \left[
\widetilde{\Lambda}_\text{int}^\text{int} (\rho, T_\text{tr}, T_\text{int}) \frac{1}{T_\text{int}} \frac{\partial T_\text{int}}{\partial x_j}
 \right]
 \nonumber \\
& \quad = \epsilon \kappa  \widetilde{\mathcal{F}} (\rho, T_\text{tr}) (T_\text{tr} - T_\text{int}),
\label{NS-alpha=0-d}
\end{align}
\end{subequations}
where $\widetilde{\mathcal{F}}$ indicates $\mathcal{F}$ for $\alpha=0$, that is,
\begin{align*}
\widetilde{\mathcal{F}} (\rho, T_\text{tr})
= C\, \frac{k_\text{B}^{1 +\beta/2}}{m^{2 + \beta/2}}\,
\rho^2 T_\text{tr}^{\beta/2},
\end{align*}
with
\begin{align*}
C = 2^{\beta +2} \sqrt{\pi }\,
\frac{\Gamma \left( \delta +1 \right)
\Gamma^{2}\left( \delta /2 \right) \Gamma \big( \left( \beta+5\right) /2\big)}
{\left[ \delta +\left( \beta +3 \right)/2 \right]
\Gamma^{2} \left( \delta \right) }\, C_\text{r}.
\end{align*}
It should be remarked that the so-called cross-diffusion terms in the
heat-flow vectors $q_{(\text{tr})i}$ and $q_{(\text{int})i}$ disappear
in this special case.

\begin{remark}
Let us consider the particular case $\alpha =\beta =0$ and note that
the following equalities hold:
\begin{align*}
& \int_{\mathbb{R}^{3}\times \mathbb{S}^{2}}
e^{-m\left\vert \bm{c}_{\ast}\right\vert ^{2}/\left( 2k_\mathrm{B}T_\mathrm{tr}\right) }
\left( \bm{c}-\bm{c}^{\prime } \right) \,\mathrm{d}\bm{c}_{\ast } \mathrm{d}\bm{\sigma }
\\
& \qquad = \int_{\mathbb{R}^{3}\times \mathbb{S}^{2}}
e^{-m\left\vert \bm{c}_{\ast } \right\vert^{2}/\left( 2k_\mathrm{B}T_\mathrm{tr}\right) }
\left( \bm{c}-\frac{\bm{c}^{\prime } + \bm{c}_{\ast }^{\prime }}{2} \right) \,
\mathrm{d}\bm{c}_{\ast } \mathrm{d}\bm{\sigma }
\\
& \qquad = \int_{\mathbb{R}^{3}\times \mathbb{S}^{2}}
e^{-m\left\vert \bm{c}_{\ast } \right\vert^{2}/\left( 2k_\mathrm{B}T_\mathrm{tr}\right) }
\left( \frac{\bm{c}-\bm{c}_{\ast }}{2} \right) \,\mathrm{d}\bm{c}_{\ast } \mathrm{d}\bm{\sigma }
\\
& \qquad = 2\pi \left( \frac{2\pi k_\mathrm{B} T_\mathrm{tr}}{m} \right)^{3/2}\bm{c}.
\end{align*}
Then, if $h$ is of the form $\left[ I/(k_\mathrm{B}T_\mathrm{int})-\delta /2\right] \bm{c}$,
$\mathcal{L}_\mathrm{r} h$ can be calculated as
\begin{equation*}
\mathcal{L}_\mathrm{r} h = 2\pi \frac{\rho }{m}
\dfrac{\Gamma^{2} \left( \delta/2\right) }{\Gamma \left( \delta \right) }
C_\mathrm{r} \left( \dfrac{I}{k_\mathrm{B}T_\mathrm{int}}-\dfrac{\delta }{2} \right) \bm{c}.
\end{equation*}
This means that $\mathcal{C}_{0} (|\bm{c}|)$ in equation~\eqref{C0} is constant when
$\alpha=\beta=0$ and is given by
\begin{equation*}
\mathcal{C}_{0} (|\bm{c}|) =2\pi \frac{\rho }{m}
\dfrac{\Gamma^{2} \left( \delta/2 \right) }{\Gamma \left( \delta \right) }
C_\mathrm{r}.
\end{equation*}
\end{remark}

\subsection{Case of $\theta = O(\epsilon)$}\label{sec:CE-2}

We next consider the case of $\theta = O(\epsilon)$ and let
\begin{align}\label{scale-2}
\theta = \overline{\kappa} \epsilon,
\end{align}
where $\overline{\kappa}$ is a positive constant \cite{BG-11}.
Then, equation~\eqref{BE-2} reads
\begin{align}\label{BE-4}
\frac{\partial f}{\partial t} + \xi_j \frac{\partial f}{\partial x_j}
= \frac{1-\overline{\kappa} \epsilon}{\epsilon}\, Q_\text{r} (f, f)
+ \overline{\kappa} Q_\text{s}(f, f).
\end{align}
%

\subsubsection{Chapman--Enskog expansion and zeroth-order solution}

Also here, we consider the Chapman--Enskog expansion~\eqref{CE-0} with~\eqref{CE-constraint}, and substitute it into equation~\eqref{BE-4}. Then,
the $O(1/\epsilon)$ term is the same as that in Sec.~\ref{sec:CE-1}, i.e.,
$Q_\text{r} (f^{(0)}, f^{(0)}) = 0$. Thus, $f^{(0)}$ is the same and is given by
$M_\text{r}$ [see~\eqref{f0}]. Therefore, $f=M_\text{r} + O(\epsilon)$,
and the stress tensor $p_{ij}$ and the heat-flow vector $q_i$ are
the same as in~\eqref{pq-0th}. On the other hand, by the use of the
expansion~\eqref{CE-0}, the term
$(\theta/\epsilon) \big( I, Q_\text{s} (f, f) \big)$ contained in equations~\eqref{transport-c}
and \eqref{transport-d} becomes
\begin{align}\label{Qs-f-f-2}
\frac{\theta}{\epsilon} \big( I, Q_\text{s} (f, f) \big)
& = \overline{\kappa}\, \big( I,\, Q_\text{s} (f^{(0)}, f^{(0)}) \big)
+ 2\overline{\kappa} \epsilon \big( I,\, Q_\text{s} (f^{(0)}, f^{(1)}) \big)
+ O(\epsilon^2)
\nonumber \\
& = \overline{\kappa}\, \big( I,\, Q_\text{s} (M_\text{r}, M_\text{r}) \big)
+ 2\overline{\kappa} \epsilon \big( I,\, Q_\text{s} (M_\text{r}, f^{(1)}) \big)
+ O(\epsilon^2).
\end{align}
The term $\big( I,\, Q_\text{s} (M_\text{r}, M_\text{r}) \big)$, which has already appeared
in Sec.~\ref{sec:CE-1} and was calculated in Appendix \ref{sec:source},
is given by the identity~\eqref{source-term-2} with \eqref{source-term-aux} for the
collision model \eqref{sigma-model}. Using
identities~\eqref{pq-0th}, equation~\eqref{Qs-f-f-2} in the form of the leading-order term plus
the error of $O(\epsilon)$, and identity~\eqref{source-term-2} with
\eqref{source-term-aux} in the transport equation~\eqref{transport} and
neglecting the $O(\epsilon)$ terms, we have the following equations:
\begin{subequations}\label{Euler-2}
\begin{align}
& \frac{\partial \rho}{\partial t} + \frac{\partial}{\partial x_j} (\rho u_j) = 0,
\label{Euler-2-a} \\
& \frac{\partial}{\partial t} (\rho u_i) + \frac{\partial}{\partial x_j} 
\left( \rho  u_i u_j + \frac{k_\text{B}}{m} \rho T_\text{tr} \delta_{ij} \right) = 0,
\label{Euler-2-b} \\
& \frac{\partial}{\partial t} \left[ \rho \left( \frac{3}{2} \frac{k_\text{B}}{m} T_\text{tr} + \frac{1}{2} |\bm{u}|^2 \right) \right]
+ \frac{\partial}{\partial x_j}  \left[ \rho u_j \left( \frac{5}{2} \frac{k_\text{B}}{m} T_\text{tr} + \frac{1}{2} |\bm{u}|^2 \right)
\right]
\nonumber \\
& \qquad \qquad 
= - \overline{\kappa}  \mathcal{F} (\rho, T_\text{tr}, T_\text{int}) (T_\text{tr} - T_\text{int}),
\label{Euler-2-c} \\
& \frac{\delta}{2} \frac{k_\text{B}}{m} \left[\frac{\partial}{\partial t} ( \rho T_\text{int} )
+ \frac{\partial}{\partial x_j} ( \rho u_j T_\text{int} ) \right]
= \overline{\kappa}  \mathcal{F} (\rho, T_\text{tr}, T_\text{int}) (T_\text{tr} - T_\text{int}). 
\label{Euler-2-d}
\end{align}
\end{subequations}
These are the Euler equations with relaxation terms proportional
to $T_\text{tr} - T_\text{int}$, which cause the interaction between the translational
and internal modes. The mathematical properties of systems of this type have been
studied in a more general framework \cite{Z-15}.
It should be noted that a system similar to system~\eqref{Euler-2} has been obtained
on the basis of extended thermodynamics \cite{ATRS-12}.

\subsubsection{First-order solution and constitutive laws at Navier--Stokes level}

From equation~\eqref{BE-4}, the equation containing the terms of $O(1)$ and higher is obtained. That is, by letting $f^{(1)}=M_\text{r}h$, we have
\begin{align}\label{EQ-h-3}
\mathcal{L}_\text{r} h = - \frac{1}{M_\text{r}} \left( \frac{\partial M_\text{r}}{\partial t}
+ \xi_j \frac{\partial M_\text{r}}{\partial x_j} \right)
+ \overline{\kappa} \frac{1}{M_\text{r}} Q_\text{s} (M_\text{r}, M_\text{r} ) + O(\epsilon).
\end{align}
Then, we take the same procedure as in Sec.~\ref{sec:CE-1} to calculate
the derivative terms on the right-hand side. To be more specific,
the time-derivative terms $\partial \rho/\partial t$, $\partial \bm{u}/\partial t$,
$\partial T_\text{tr}/\partial t$, and $\partial T_\text{int}/\partial t$ arising
from $\partial M_\text{r}/\partial t$ are replaced
with the space derivative terms and the relaxation term with the help of
equations~\eqref{Euler-2} [note that the system~\eqref{Euler-2} holds with the error of
$O(\epsilon)$]. As the result, neglecting the terms of $O(\epsilon)$, we obtain
from equation~\eqref{EQ-h-3} the integral equation for $h$ in the following form:
\begin{align*}
\mathcal{L}_\text{r} h = H_1 + H_2,
\end{align*}
where
\begin{subequations}\label{H1-H2}
\begin{align}
H_1 & = - \frac{m}{2k_\text{B} T_\text{tr}} 
\left( \frac{\partial u_i}{\partial x_j} + \frac{\partial u_j}{\partial x_i} \right) A_{ij} (\bm{c})
- \frac{1}{T_\text{tr}} \frac{\partial T_\text{tr}}{\partial x_j} B_j (\bm{c})
- \frac{1}{T_\text{int}} \frac{\partial T_\text{int}}{\partial x_j}
C_j (\bm{c}, I),
\label{H1-H2-a} \\
H_2 & = \frac{m}{k_\text{B} \rho} \left[  \frac{1}{T_\text{tr}} \left( \frac{m}{3k_\text{B} T_\text{tr}}
|\bm{c}|^2 - 1 \right) 
- \frac{1}{T_\text{int}}  \left( \frac{2I}{\delta k_\text{B} T_\text{int}} -1 \right) \right]
\overline{\kappa} \big (I,\, Q_\text{s} (M_\text{r}, M_\text{r}) \big)
\nonumber \\
& \quad
+ \overline{\kappa} \frac{1}{M_\text{r}} Q_\text{s} (M_\text{r}, M_\text{r} ).
\label{H1-H2-b}
\end{align}
\end{subequations}
Here, the relaxation term
$\overline{\kappa} \mathcal{F}(\rho, T_\text{tr}, T_\text{int}) (T_\text{tr}-T_\text{int})$
has been replaced with the original
$\overline{\kappa} \big (I,\, Q_\text{s} (M_\text{r}, M_\text{r}) \big)$
[cf.~\eqref{source-term-2}] for convenience.

Let us decompose the solution $h$ as
\begin{align}\label{decomp-h}
h = h_1 + h_2,\qquad \mathcal{L}_\text{r} h_1 = H_1, \qquad \mathcal{L}_\text{r} h_2 = H_2.
\end{align}
The equation for $h_1$ is the same as~\eqref{EQ-h-2} in Sec.~\ref{sec:CE-1}, so that
$h_1$ is given by the right-hand side of equality~\eqref{sol-h}. Therefore, we consider the
equation for $h_2$ below.

It can be easily seen that the right-hand side $H_2$ belongs to $( M_\text{r} \ker \mathcal{L}_\text{r})^\perp$.
Therefore, the solution $h_2$ is uniquely obtained in the same space
$( M_\text{r} \ker \mathcal{L}_\text{r})^\perp$ (cf.~Corollary \ref{Cor2} or Remark \ref{solvability}).

We now try to calculate the stress tensor $p_{ij}$ and the heat-flow vectors $q_{(\text{tr})i}$
and $q_{(\text{int})i}$ using $f = M_\text{r} (1 + \epsilon h) +O(\epsilon^2) =
M_\text{r} [1 + \epsilon (h_1 + h_2)]+O(\epsilon^2)$, i.e.,
\begin{align*}
& p_{ij} = m \Big( c_i c_j, \, M_\text{r} [1 + \epsilon (h_1 + h_2)] \Big) + O(\epsilon^2),
\\
& q_{(\text{tr})i} = \frac{m}{2} \Big( c_i |\bm{c}|^2,\,  M_\text{r} [1 + \epsilon (h_1 + h_2)] \Big)  + O(\epsilon^2),
\\
& q_{(\text{int})i} = \Big( c_i I,\,  M_\text{r} [1 + \epsilon (h_1 + h_2)] \Big)  + O(\epsilon^2).
\end{align*}
Actually, we need to consider only the contribution from $h_2$ because the other
contributions have already been obtained in Sec.~\ref{sec:CE-1}.
In other words, we just consider $(c_i c_j,\, M_\text{r} h_2)$, $(c_i |\bm{c}|^2,\, M_\text{r} h_2)$,
and $(c_i I,\, M_\text{r} h_2)$.

Using expressions~\eqref{def-ABC}, equations~\eqref{EQ-decomp},
Proposition \ref{prop:positivity-kernel}, and the decomposition~\eqref{decomp-h}, in addition to the
fact that $h_2 \in ( M_\text{r} \ker \mathcal{L}_\text{r})^\perp$, we obtain
the following equalities:
\begin{subequations}\label{moment-h2}
\begin{align}
(c_i c_j,\, M_\text{r} h_2) & = \left( \Big(c_i c_j - \frac{1}{3} |\bm{c}|^2 \delta_{ij} \Big),\, M_\text{r} h_2 \right)
= \big( A_{ij} (\bm{c}), \,  M_\text{r} h_2 \big)
\nonumber \\
&
=  \big( \mathcal{L}_\text{r} (A_{ij} (\bm{c}) \mathcal{A} (|\bm{c}|, I) ),\, M_\text{r} h_2 \big)
=  \big( M_\text{r} A_{ij} (\bm{c}) \mathcal{A} (|\bm{c}|, I),\, \mathcal{L}_\text{r} h_2 \big)
\nonumber \\
& =  \big( M_\text{r} A_{ij} (\bm{c}) \mathcal{A} (|\bm{c}|, I),\, H_2 \big),
\label{moment-h2-a} \\
\frac{m}{2k_\text{B} T_\text{tr}} (c_i |\bm{c}|^2,\, M_\text{r} h_2)
& = \left( c_i \Big( \frac{m|\bm{c}|^2}{2k_\text{B} T_\text{tr}} - \frac{5}{2} \Big),\, M_\text{r} h_2 \right)
= \big( B_i (\bm{c}), \,  M_\text{r} h_2 \big)
\nonumber \\
&
=  \big( \mathcal{L}_\text{r} (c_i \mathcal{B} (|\bm{c}|, I) ),\, M_\text{r} h_2 \big)
=  \big( M_\text{r} c_i \mathcal{B} (|\bm{c}|, I),\, \mathcal{L}_\text{r} h_2 \big)
\nonumber \\
& =  \big( M_\text{r} c_i \mathcal{B} (|\bm{c}|, I),\, H_2 \big),
\label{moment-h2-b} \\
\frac{1}{k_\text{B} T_\text{int}} (c_i I,\, M_\text{r} h_2)
& = \left( c_i \Big( \frac{I}{k_\text{B} T_\text{int}} - \frac{\delta}{2} \Big),\, M_\text{r} h_2 \right)
= \big( C_i (\bm{c}), \,  M_\text{r} h_2 \big)
\nonumber \\
&
=  \big( \mathcal{L}_\text{r} (c_i \mathcal{C} (|\bm{c}|, I) ),\, M_\text{r} h_2 \big)
=  \big( M_\text{r} c_i \mathcal{C} (|\bm{c}|, I),\, \mathcal{L}_\text{r} h_2 \big)
\nonumber \\
& =  \big( M_\text{r} c_i \mathcal{C} (|\bm{c}|, I),\, H_2 \big).
\label{moment-h2-c} 
\end{align}
\end{subequations}
It should be noted here that $Q_\text{r} (M_\text{r}, M_\text{r})$ is a function
of $|\bm{c}|$ and $I$, as shown in Appendix \ref{sec:symmetry}, and
thus, $H_2$ is also a function of $|\bm{c}|$ and $I$. On the other hand,
$\int_{\mathbb{R}^3} A_{ij} (\bm{c}) F(|\bm{c}|)\, \text{d}\bm{c} = 0$ and
$\int_{\mathbb{R}^3} c_i F(|\bm{c}|)\, \text{d}\bm{c} = 0$ hold for an
arbitrary function $F(|\bm{c}|)$ of $|\bm{c}|$ for which the integrals make sense.
Therefore, the last line of~\eqref{moment-h2-a}, that of~\eqref{moment-h2-b},
and that of~\eqref{moment-h2-c} are all zero, so that $(c_i c_j,\, M_\text{r} h_2)$,
$(c_i |\bm{c}|^2,\, M_\text{r} h_2)$, and $(c_i I,\, M_\text{r} h_2)$ all vanish.
This means that the contributions of $h_2$ to $p_{ij}$, $q_{(\text{tr}) i}$, and
$q_{(\text{int}) i}$ are zero.

In summary, $p_{ij}$, $q_{(\text{tr}) i}$, and $q_{(\text{int}) i}$ are given by the expressions~\eqref{pij}, \eqref{qi-a}, and \eqref{qi-b}, respectively. When $\alpha=0$
[cf.~\eqref{sigma-model}], they are given by the same expressions~\eqref{pij},
\eqref{qi-a}, and \eqref{qi-b} with $\Lambda_\mu = \widetilde{\Lambda}_\mu$,
$\Lambda_\text{tr}^\text{tr} = \widetilde{\Lambda}_\text{tr}^\text{tr}$,
$\Lambda_\text{int}^\text{int} = \widetilde{\Lambda}_\text{int}^\text{int}$,
$\Lambda_\text{int}^\text{tr} = \Lambda_\text{tr}^\text{int}=0$
[cf.~identities~\eqref{Lambda-1-alpha=0} and \eqref{Lambda-2-alpha=0}].

\smallskip
\subsubsection{Source term}

The remaining task is to investigate the $O(\epsilon)$-term in the source term
in the system~\eqref{Qs-f-f-2}. Since $f^{(1)} = M_\text{r} h = M_\text{r} (h_1 + h_2)$,
$\big( I, Q_\text{s} (M_\text{r}, f^{(1)}) \big)$ is written as
\begin{align}\label{decomp-Qs}
\big( I, Q_\text{s} (M_\text{r}, f^{(1)}) \big)
= \big( I, Q_\text{s} (M_\text{r}, M_\text{r} h_1) \big)
+ \big( I, Q_\text{s} (M_\text{r}, M_\text{r} h_2) \big).
\end{align}

As shown in Appendix \ref{sec:derivation-2}, we have the following expression for the first term on the
right-hand side:
\begin{align}\label{moment-h_1}
& \big( I, Q_\text{s} (M_\text{r}, M_\text{r} h_1) \big)
\nonumber \\
& \;\;\;
= \frac{m\rho ^{2}C_\text{s}}{4 \pi^2 \left( k_\text{B}T_\text{tr}\right)^{3}
\left( k_\text{B}T_\text{int}\right)^{\delta } \Gamma \left( \delta \right) }
\nonumber \\
& \qquad \quad \times
\int_{\lbrack 0,1]\times \left( \mathbb{R}^{3}\times \mathbb{R}_{+} \right) ^{2}}
h_{1} (\bm{c}, I) \left\vert \mathbf{c-c}_{\ast }\right\vert^{\beta}
e^{-m \left( \left\vert \mathbf{c}\right\vert^{2}+\left\vert \mathbf{c}_{\ast }\right\vert^{2}\right)
/\left( 2k_\text{B}T_\text{tr}\right) }
\nonumber \\
& \qquad \qquad \qquad \qquad \times
\left[ \frac{m}{4} (1-R)\left\vert \mathbf{c-c}_{\ast } \right\vert^{2}-R\left( I+I_{\ast}\right) \right]
\left( I+I_{\ast }\right) ^{\alpha } \left( II_{\ast }\right)^{\delta /2-1}
\nonumber \\
& \qquad \qquad \qquad \qquad \times
e^{-\left( I+I_{\ast }\right) /\left( k_\text{B}T_\text{int}\right)}
R^{\left( \beta +1\right) /2} \left( 1-R\right) ^{\alpha +\delta -1}
\text{d}R \text{d}\bm{c}_{\ast } \text{d} I_{\ast } \text{d}\bm{c} \text{d}I,
\end{align}
where the arguments $t$ and $\bm{x}$ are omitted in $h_1$. If we consider the
integral
\begin{align*}
\mathcal{I}_s = \int_{\mathbb{R}^3} |\bm{c}-\bm{c}_*|^s e^{-m(|\bm{c}|^2 + |\bm{c}_*|^2)/(2k_\text{B} T_\text{tr})}
\text{d} \bm{c}_*,
\end{align*}
with a positive constant $s$, then it is seen that $\mathcal{I}_s$ is spherically symmetric in $\bm{c}$,
that is, a function of $|\bm{c}|$, for the same reason as Appendix \ref{sec:symmetry}.
Now, let us consider the integral
\begin{align*}
\int_{\mathbb{R}^3} h_1 (\bm{c}, I) \mathcal{I}_s (|\bm{c}|) \text{d} \bm{c},
\end{align*}
and recall that $h_1 (\bm{c}, I)$ is given by the right-hand side of equality~\eqref{sol-h}.
Then, it is expressed as
\begin{align*}
\int_{\mathbb{R}^3} h_1 (\bm{c}, I) \mathcal{I}_{s} (|\bm{c}|) \text{d} \bm{c}
=  &- \frac{m}{2k_\text{B} T_\text{tr}} 
\left( \frac{\partial u_i}{\partial x_j} + \frac{\partial u_j}{\partial x_i} \right)
\int_{\mathbb{R}^3} A_{ij} (\bm{c}) \mathcal{A} (|\bm{c}|, I) \mathcal{I}_s (|\bm{c}|) \text{d}\bm{c}
\\
& - \frac{1}{T_\text{tr}} \frac{\partial T_\text{tr}}{\partial x_j}
\int_{\mathbb{R}^3} c_j \mathcal{B} (|\bm{c}|, I) \mathcal{I}_s (|\bm{c}|) \text{d}\bm{c}
\\
& - \frac{1}{T_\text{int}} \frac{\partial T_\text{int}}{\partial x_j}
\int_{\mathbb{R}^3} c_j \mathcal{C} (|\bm{c}|, I) \mathcal{I}_s (|\bm{c}|) \text{d}\bm{c}.
\end{align*}
However, for the same reason as for equalities~\eqref{moment-h2}, all three integrals
on the right-hand side of the above equation vanish.
From this fact and the equality~\eqref{moment-h_1}, it follows that
\begin{align}\label{I-Qs-Mrh1}
\big( I, Q_\text{s} (M_\text{r}, M_\text{r} h_1) \big) = 0.
\end{align}

Next, we consider the second term on the right-hand side of the decomposition~\eqref{decomp-Qs}.
Let us write $H_2$ [see~\eqref{H1-H2-b}] in a slightly different way, that is,
\begin{align*}
H_2 = \overline{\kappa} \big (I,\, Q_\text{s} (M_\text{r}, M_\text{r}) \big) D,
\end{align*}
where
\begin{align}\label{D}
D & = \frac{m}{k_\text{B} \rho} \left[  \frac{1}{T_\text{tr}} \left( \frac{m}{3k_\text{B} T_\text{tr}}
|\bm{c}|^2 - 1 \right) 
- \frac{1}{T_\text{int}}  \left( \frac{2I}{\delta k_\text{B} T_\text{int}} -1 \right) \right]
+ \frac{ Q_\text{s} (M_\text{r}, M_\text{r} ) }
{ M_\text{r} \big (I,\, Q_\text{s} (M_\text{r}, M_\text{r}) \big) }.
\end{align}
Then it is easily seen that
$D \in ( M_\text{r} \ker \mathcal{L}_\text{r})^\perp$.
The last term is seemingly divergent as $T_\text{tr} \to T_\text{int}$ because
$\big(I, Q_\text{s} (M_\text{r}, M_\text{r}) \big)$ is proportional to $T_\text{tr}-T_\text{int}$.
However, we will see that it is not the case below.
Therefore, the integral equation
\begin{align}\label{eq-D}
\mathcal{L}_\text{r} \widetilde{D} = D,
\end{align}
has a unique solution $\widetilde{D}$ such that
$\widetilde{D} \in ( M_\text{r} \ker \mathcal{L}_\text{r})^\perp$
because $\mathcal{L}_\text{r}$ is a Fredholm operator in $L^2 (M_\text{r} \text{d}\bm{\xi} \text{d}I)$
(Corollary \ref{Cor2} or Remark \ref{solvability} for $\theta=0$). Thus, $h_2$ can be expressed as
\begin{align*}
h_2 = \overline{\kappa} \big (I,\, Q_\text{s} (M_\text{r}, M_\text{r}) \big) \widetilde{D}.
\end{align*}
so that it follows that
\begin{equation}\label{I-Qs-Mrh2}
\big( I, Q_\text{s} \left( M_\text{r}, M_\text{r} h_{2} \right) \big)
= \overline{\kappa} \big( I, Q_\text{s} \left( M_\text{r}, M_\text{r} \right) \big)
\left( I, Q_\text{s} \big( M_\text{r}, M_\text{r} \widetilde{D} \big) \right).
\end{equation}
If $\big( I, Q_\text{s} \big( M_\text{r}, M_\text{r} \widetilde{D} \big) \big)$ is bounded,
then we can conclude that $\big( I, Q_\text{s} \left( M_\text{r}, M_\text{r} h_{2} \right) \big)$
is proportional to $ \big( I, Q_\text{s} \left( M_\text{r}, M_\text{r} \right) \big)$, or
equivalently proportional to $T_\text{tr} - T_\text{int}$ [cf.~\eqref{Qs-Mr-Mr}].
Therefore, we proceed by proving the boundedness of
$\big( I, Q_\text{s} \big( M_\text{r}, M_\text{r} \widetilde{D} \big) \big)$ in the following.

The Fredholmness of $\mathcal{L}_{r}$ ($\mathcal{L}_{r}$ being a closed
linear operator with a closed range) indicates that, for any function
$g(\bm{\xi }, I) \in \left( M_\text{r} \ker \mathcal{L}_{r} \right)^{\perp }
\cap \mathcal{D}(\mathcal{L}_\text{r})$,
there exists a constant $\mu >0$ such that 
\begin{equation*}
\left( \mathcal{L}_{r}g, M_{r} \mathcal{L}_{r}g \right) \geq \mu \left(g, M_{r}g\right) 
\end{equation*}
holds (cf.~\cite{Go-66}; Chap.~IV, Sec.~5.1 in \cite{K-80}). Thus, we have 
\begin{equation}\label{DMD}
\left( \widetilde{D},M_{r}\widetilde{D}\right) \leq \frac{1}{\mu }\left( 
\mathcal{L}_{r}\widetilde{D},M_{r}\mathcal{L}_{r}\widetilde{D}\right) =\frac{%
1}{\mu }\left( D,M_{r}D\right) \text{.}
\end{equation}
With the help of this inequality, one can show the following inequality (see Appendix \ref{sec:proof-ineq}):
\begin{align}\label{inequality}
\left( I, Q_\text{s} \big( M_\text{r}, M_\text{r} \widetilde{D} \big) \right)^2
\le  C_\text{g}   (D, M_\text{r} D).
\end{align}
Here and in what follows, $C_\text{g}$ indicates a generic positive constant depending on
the macroscopic quantities $\rho$, $T_\text{tr}$, and $T_\text{int}$. Therefore,
we have to prove that $(D, M_\text{r} D)$ is bounded.

For this purpose, we consider $M_\text{r}^{-1} Q_\text{s} (M_\text{r}, M_\text{r})$,
which occurs in $D$ [see~\eqref{D}],  using
the first line of~\eqref{Qs-Mr-Mr-3}. Let us first estimate the factor
$e^{-\eta (1-R) E/\zeta} - e^{-\eta(I+I_*)/\zeta}$ in~\eqref{Qs-Mr-Mr-3},
noting that [cf.~identities~\eqref{zeta-eta}]
\begin{align*}
& e^{-\eta (1-R) E/\zeta} - e^{-\eta(I+I_*)/\zeta}
\\
& \qquad = \left\{
\begin{array}{l}
e^{-(1-R)E/\zeta} - e^{-(I+I_*)/\zeta},  \qquad \qquad \qquad
(\text{for}\; \eta=1,\; \text{i.e.,}\; T_\text{tr} > T_\text{int}), \\
e^{E/\zeta} \left( e^{-RE/\zeta} - e^{-m|\bm{c}-\bm{c}_*|^2/(4\zeta)} \right), \qquad
(\text{for}\; \eta=-1,\; \text{i.e.,}\; T_\text{tr} < T_\text{int}).
\end{array}
\right.
\end{align*}
Here, note that if $0 \le s_1 \le s_2$, then
\begin{align*}
0 < e^{-s_1} - e^{-s_2} \le e^{-s_1} (s_2-s_1) \le s_1 + s_2.
\end{align*}
Thus, for any nonnegative $s_1$ and $s_2$, it holds that
\begin{align*}
|e^{-s_1} - e^{-s_2}| \le s_1 + s_2.
\end{align*}
Using this relation and recalling that
\[
E = \frac{m}{4} |\bm{c}-\bm{c}_*|^2 + I + I_*,
\]
one obtains
\begin{align}\label{estimate:e-e}
\begin{aligned}
\left|  e^{-(1-R)E/\zeta }-e^{-\left( I+I_{\ast }\right) /\zeta}\right|
& \leq \zeta^{-1} \left[ (1-R)E+I+I_{\ast } \right] \leq 2\zeta^{-1} E,
\\
\left| e^{-RE/\zeta } -e^{-m\left( \left| \bm{c}-\bm{c}_* \right|^{2}\right) /\left( 4\zeta \right)} \right|
& \leq \zeta^{-1} \left( RE+\frac{m}{4} |\bm{c}-\bm{c}_*|^{2} \right)
\leq 2\zeta^{-1} E.
\end{aligned}
\end{align}
Incidentally, it is noted that $E$ is estimated as
\begin{align}\label{estimate:E}
E & = \frac{m}{2} \left\vert \bm{c} \right\vert^{2}
+ \frac{m}{2} \left\vert \bm{c}_* \right\vert^{2}
- \frac{m}{4} \left\vert \bm{c} + \bm{c}_* \right\vert^{2} + I + I_{\ast }
\nonumber \\
& \leq \frac{m}{2}\left\vert \bm{c} \right\vert^{2} + \frac{m}{2} \left\vert \bm{c}_{\ast } \right\vert^{2} + I + I_{\ast}.
\end{align}

With these results, $M_\text{r}^{-1} Q_\text{s} (M_\text{r}, M_\text{r})$ can easily be estimated as follows
(see Appendix \ref{sec:proof-estimate}):
\begin{align}\label{estimate-Qs}
& M_\text{r}^{-1} |Q_\text{s} (M_\text{r}, M_\text{r})|
\nonumber \\
& \qquad \le \left\{
\begin{array}{l}
C_\text{g} |T_\text{tr} - T_\text{int}| (1 + |\bm{c}|^2 + |\bm{c}|^{\beta+2})
(1 + I + I^{\alpha+1}) e^{I/\zeta}, \qquad (T_\text{tr} > T_\text{int}),
\\[2mm]
C_\text{g} |T_\text{tr} - T_\text{int}| (1 + |\bm{c}|^2 + |\bm{c}|^{\beta+2})
(1 + I + I^{\alpha+1}) e^{m|\bm{c}|^2/(2\zeta)},
\qquad (T_\text{tr} < T_\text{int}),
\end{array}
\right.
\end{align}
where, as mentioned above, $C_\text{g}$ indicates a generic positive constant depending
on the macroscopic quantities. This estimate shows that the last term
on the right-hand side of equation~\eqref{D} is bounded as $T_\text{tr} \to T_\text{int}$.

Now, we try to estimate $(D, M_\text{r} D)$. It has implicitly been assumed that
$T_\text{tr}$ and $T_\text{int}$ are strictly positive and bounded. Here, we write
it explicitly as
$0 < C_l \le T_\text{tr} \le C_u < \infty$ and $0 < C_l \le T_\text{int} \le C_u < \infty$;
then, we additionally assume that $|T_\text{tr} - T_\text{int}| \le C_l/3$. 
Thus, we have the following inequalities
\begin{align*}
& 2\frac{T_\text{tr}-T_\text{int}}{T_\text{tr}T_\text{int}}-\frac{1}{T_\text{int}}
\leq \frac{2}{3} \frac{C_l}{T_\text{tr} T_\text{int}} - \frac{1}{T_\text{int}}
\leq \frac{2}{3T_\text{int}}-\frac{1}{T_\text{int}}
= -\frac{1}{3T_\text{int}},
\qquad (T_\text{tr} > T_\text{int}),
\\
& \frac{T_\text{int}-T_\text{tr}}{T_\text{tr}T_\text{int}}-\frac{1}{2T_\text{tr}}
\leq \frac{1}{3} \frac{C_l}{T_\text{tr} T_\text{int}} - \frac{1}{2T_\text{tr}}
\leq \frac{1}{3T_\text{tr}}-\frac{1}{2T_\text{tr}}
= -\frac{1}{6T_\text{tr}},
\qquad \;\;\, (T_\text{tr} < T_\text{int}),
\end{align*}
which, respectively, indicate that
\begin{align}\label{bound}
\begin{aligned}
& e^{2I/\zeta } M_\text{r}
\leq C_\text{g} I^{\delta /2-1} e^{-m\left\vert \bm{c} \right\vert^{2}/\left( 2k_\text{B}T_\text{tr}\right) }
e^{-I/\left(3k_\text{B} T_\text{int} \right)}, \qquad \;\;\; (T_\text{tr}>T_\text{int}),
\\
& e^{m\left\vert \bm{c} \right\vert^{2}/\zeta }M_\text{r}
\leq C_\text{g} I^{\delta/2-1} e^{-m\left\vert \bm{c} \right\vert^{2}/\left( 6k_\text{B}T_\text{tr} \right) }
e^{-I/\left( k_\text{B}T_\text{int} \right) }, \qquad (T_\text{tr} <T_\text{int}).
\end{aligned}
\end{align}
Let us decompose $D$ as
\begin{align*}
& D = D_1 + D_2,
\\
& D_1 = \frac{m}{k_\text{B} \rho} \left[  \frac{1}{T_\text{tr}} \left( \frac{m}{3k_\text{B} T_\text{tr}}
|\bm{c}|^2 - 1 \right) 
- \frac{1}{T_\text{int}}  \left( \frac{2I}{\delta k_\text{B} T_\text{int}} -1 \right) \right],
\\
& D_2 = \frac{ Q_\text{s} (M_\text{r}, M_\text{r} ) }
{ M_\text{r} \big (I,\, Q_\text{s} (M_\text{r}, M_\text{r}) \big) }
=  \frac{1}{\mathcal{F} M_\text{r}} \frac{ Q_\text{s} (M_\text{r}, M_\text{r}) }{T_\text{tr}-T_\text{int}},
\end{align*}
where identity~\eqref{Qs-Mr-Mr} has been used in the last equality for $D_2$, and, for convenience, let 
\begin{align*}
\mathcal{I} ( |\bm{c}|, I ) = (1 + |\bm{c}|^2 + |\bm{c}|^{\beta+2}) (1 + I + I^{\alpha+1}), 
\end{align*}
so that it holds that
\begin{align*}
|D_1| \le C_\text{g} \mathcal{I} (|\bm{c}|, I),
\end{align*}
with a generic positive constant $C_\text{g}$. Then, we have
\begin{align*}
&\left| (D, M_\text{r} D) \right|  \le \left| (D_1, M_\text{r} D_1) \right| + 2 \left| (D_1, M_\text{r} D_2) \right|
+ \left| (D_2, M_\text{r} D_2) \right|
\\
=& \left| (D_1, M_\text{r} D_1) \right|
+ \frac{2}{\left| \mathcal{F} \right|} \frac{ \left| \big( D_1, Q_\text{s} (M_\text{r}, M_\text{r}) \big) \right|}{ \left| T_\text{tr}-T_\text{int} \right| }
+ \frac{1}{\mathcal{F}^2} \frac{ \left| \big( M_\text{r}^{-1} Q_\text{s} (M_\text{r}, M_\text{r}),\, Q_\text{s} (M_\text{r}, M_\text{r}) \big) \right| }
{ \left| T_\text{tr}-T_\text{int} \right|^2 }.
\end{align*}
It is obvious that $|(D_1, M_\text{r} D_1)|$ is bounded. In addition, with the help of estimates~\eqref{estimate-Qs}, the following inequalities follow:
\begin{align*}
& \frac{ \left| \big( D_1,\, Q_\text{s} (M_\text{r}, M_\text{r}) \big) \right|}{ \left| T_\text{tr}-T_\text{int} \right| }
\\
& \qquad \le
\left\{
\begin{array}{l}
C_\text{g}\, \big( \mathcal{I},\, e^{I/\zeta} M_\text{r} \mathcal{I} \big)\,
\le \, C_\text{g} \big( \mathcal{I},\, e^{2I/\zeta} M_\text{r} \mathcal{I} \big),
\qquad \qquad \qquad (T_\text{tr} > T_\text{int}),
\\[2mm]
C_\text{g}\, \big( \mathcal{I},\, e^{m|\bm{c}|^2/(2\zeta)} M_\text{r} \mathcal{I} \big)\,
\le \, C_\text{g} \big( \mathcal{I},\, e^{m|\bm{c}|^2/\zeta} M_\text{r} \mathcal{I} \big),
\qquad \, (T_\text{tr} <\, T_\text{int}),
\end{array}
\right.
\\
& \frac{ \left| \big( M_\text{r}^{-1} Q_\text{s} (M_\text{r}, M_\text{r}),\, Q_\text{s} (M_\text{r}, M_\text{r}) \big) \right| }
{ \left| T_\text{tr}-T_\text{int} \right|^2 }
\\
& \qquad \le
\left\{
\begin{array}{l}
C_\text{g}\, \big( e^{I/\zeta} \mathcal{I},\, e^{I/\zeta} M_\text{r} \mathcal{I} \big)\,
\le \, C_\text{g} \big( \mathcal{I},\, e^{2I/\zeta} M_\text{r} \mathcal{I} \big),
\qquad \qquad \qquad \qquad \;\;\, (T_\text{tr} > T_\text{int}),
\\[2mm]
C_\text{g}\, \big( e^{m|\bm{c}|^2/(2\zeta)} \mathcal{I},\, e^{m|\bm{c}|^2/(2\zeta)} M_\text{r} \mathcal{I} \big)\,
\le \, C_\text{g} \big( \mathcal{I},\, e^{m|\bm{c}|^2/\zeta} M_\text{r} \mathcal{I} \big),
\qquad \, (T_\text{tr} <\, T_\text{int}).
\end{array}
\right.
\end{align*}
In view of estimate \eqref{bound}, both 
$\left| \left( D_1,\, Q_\text{s} (M_\text{r}, M_\text{r}) \right) \right| / \left| T_\text{tr}-T_\text{int} \right|$
and
$\left| \big( M_\text{r}^{-1} Q_\text{s} (M_\text{r}, M_\text{r}),\, Q_\text{s} (M_\text{r}, M_\text{r}) \big) \right|$
$/$$|T_\text{tr} - T_\text{int}|^2$
are seen to be bounded. In consequence, $(D, M_\text{r} D)$ is bounded. 

From the estimate~\eqref{inequality}, it is concluded that
$\big( I, Q_\text{s} ( M_\text{r}, M_\text{r} \widetilde{D} ) \big)$ is bounded.
Letting
\begin{align}\label{K}
\mathcal{K} (\rho, T_\text{tr}, T_\text{int})
= 2 \big( I, Q_\text{s} ( M_\text{r}, M_\text{r} \widetilde{D} ) \big),
\end{align}
and taking account of identities~\eqref{I-Qs-Mrh1}, \eqref{I-Qs-Mrh2}, and \eqref{Qs-Mr-Mr}
in decomposition~\eqref{decomp-Qs}, one obtains
\begin{equation*}
\big( I, Q_\text{s} ( M_\text{r}, f^{(1)} ) \big)
= \frac{1}{2} \overline{\kappa} \mathcal{F} \big( \rho, T_\text{tr}, T_\text{int} \big)\,
\mathcal{K} \big( \rho, T_\text{tr}, T_\text{int} \big)
(T_\text{tr} - T_\text{int}).
\end{equation*}
Therefore, expression~\eqref{Qs-f-f-2}, i.e., the source term
$(\theta/\epsilon) \big(I, Q_\text{s} (f, f) \big)$ included in equations~\eqref{transport-c} and
\eqref{transport-d}, is recast as
\begin{align}\label{source-upto-eps}
\frac{\theta}{\epsilon} \big( I, Q_\text{s} (f, f) \big)
= \overline{\kappa} \mathcal{F} (\rho, T_\text{tr}, T_\text{int})
[1 + \epsilon \mathcal{K} (\rho, T_\text{tr}, T_\text{int}) ] (T_\text{tr} - T_\text{int})
+ O(\epsilon^2),
\end{align}
where $\mathcal{F}$ and $\mathcal{K}$ are, respectively, given by~\eqref{source-term-aux}
and \eqref{K}.

\subsubsection{Two-temperature Navier--Stokes equations}\label{sec:two-temp-NS}

Recall that the stress tensor $p_{ij}$ and heat-flow vectors $q_{\text{(tr)}i}$ and
$q_{\text{(int)}i}$ are the same as those for $\theta= O(\epsilon^2)$ and are given
by the expressions~\eqref{pij}, \eqref{qi-a}, and \eqref{qi-b}, respectively. If we use these
results as well as identity~\eqref{source-upto-eps} in the transport equations \eqref{transport}
and neglect the terms of $O(\epsilon^2)$, we obtain the following equations:
\begin{subequations}\label{NS-2}
\begin{align}
& \frac{\partial \rho}{\partial t} + \frac{\partial}{\partial x_j} (\rho u_j) = 0,
\label{NS-2-a} \\
& \frac{\partial}{\partial t} (\rho u_i) + \frac{\partial}{\partial x_j} 
(\rho  u_i u_j) + \frac{k_B}{m} \frac{\partial}{\partial x_i} (\rho T_\text{tr})
\nonumber \\
& \qquad \qquad
- \epsilon \frac{\partial}{\partial x_j} 
\left[ \Lambda_\mu (\rho, T_\text{tr}, T_\text{int})
\left( \frac{\partial u_i}{\partial x_j} + \frac{\partial u_j}{\partial x_i}
- \frac{2}{3} \frac{\partial u_k}{\partial x_k} \delta_{ij} \right) \right]
 = 0,
\label{NS-2-b} \\
& \frac{\partial}{\partial t} \left[ \rho \left( \frac{3}{2} \frac{k_B}{m} T_\text{tr} + \frac{1}{2} |\bm{u}|^2 \right) \right]
+ \frac{\partial}{\partial x_j}  \left[ \rho u_j \left( \frac{5}{2} \frac{k_B}{m} T_\text{tr} + \frac{1}{2} |\bm{u}|^2 \right) \right]
\nonumber \\
& \qquad \qquad
- \epsilon \frac{\partial}{\partial x_j} \left[ u_i \Lambda_\mu (\rho, T_\text{tr}, T_\text{int})
\left( \frac{\partial u_i}{\partial x_j} + \frac{\partial u_j}{\partial x_i}
- \frac{2}{3} \frac{\partial u_k}{\partial x_k} \delta_{ij} \right) \right]
\nonumber \\
& \qquad \qquad
- \epsilon \frac{\partial}{\partial x_j} \left[
 \Lambda_\text{tr}^\text{tr} (\rho, T_\text{tr}, T_\text{int}) \frac{1}{T_\text{tr}} \frac{\partial T_\text{tr}}{\partial x_j}
+ \Lambda_\text{int}^\text{tr} (\rho, T_\text{tr}, T_\text{int}) \frac{1}{T_\text{int}} \frac{\partial T_\text{int}}{\partial x_j}
 \right]
\nonumber \\
& \quad = - \overline{\kappa}  \mathcal{F} (\rho, T_\text{tr}, T_\text{int})
[1 + \epsilon \mathcal{K} (\rho, T_\text{tr}, T_\text{int}) ]
(T_\text{tr} - T_\text{int}),
\label{NS-2-c} \\[2mm]
& \frac{\partial}{\partial t} \left( \frac{\delta}{2} \frac{k_B}{m} \rho T_\text{int} \right)
+ \frac{\partial}{\partial x_j} \left( \frac{\delta}{2} \frac{k_B}{m} \rho u_j T_\text{int} \right)
\nonumber \\
& \qquad \qquad
- \epsilon \frac{\partial}{\partial x_j} \left[
 \Lambda_\text{tr}^\text{int} (\rho, T_\text{tr}, T_\text{int}) \frac{1}{T_\text{tr}} \frac{\partial T_\text{tr}}{\partial x_j}
+ \Lambda_\text{int}^\text{int} (\rho, T_\text{tr}, T_\text{int}) \frac{1}{T_\text{int}} \frac{\partial T_\text{int}}{\partial x_j}
 \right]
 \nonumber \\
& \quad =  \overline{\kappa}  \mathcal{F} (\rho, T_\text{tr}, T_\text{int})
[1 + \epsilon \mathcal{K} (\rho, T_\text{tr}, T_\text{int}) ] (T_\text{tr} - T_\text{int}),
\label{NS-2-d}
\end{align}
\end{subequations}
where $\Lambda_\mu$, $\Lambda_\text{tr}^\text{tr}$, $\Lambda_\text{int}^\text{tr}$,
$\Lambda_\text{tr}^\text{int}$, and $\Lambda_\text{int}^\text{int}$ are given
by~\eqref{Lambda-1} and \eqref{Lambda-2}, and $\mathcal{F}$ and $\mathcal{K}$
are, respectively, given by~\eqref{source-term-aux} and \eqref{K}, as
mentioned above. These equations are basically of the same form as  equations~\eqref{NS}
when $\theta = O(\epsilon^2)$. The only difference appears in the relaxation terms.
To be more specific, the right-hand sides of equations~\eqref{NS-c} and \eqref{NS-d}
when $\theta=O(\epsilon^2)$ are of $O(\epsilon)$, whereas those of
equations~\eqref{NS-2-c} and \eqref{NS-2-d} contain terms of $O(1)$ and
$O(\epsilon)$. Although the boundedness of $|T_\text{tr} - T_\text{int}|$ is
assumed for the estimate \eqref{bound}, it should be emphasized that
its smallness is not required to derive the system~\eqref{NS-2}.

It should be mentioned that the two-temperature Navier--Stokes equations
of the form of~\eqref{NS-2} [i.e., with the relaxation terms of $O(1)$, not of $O(\epsilon)$]
have been derived from the ES model for a polyatomic gas \cite{ALPP-00} by an appropriate
parameter setting \cite{ABGK-20}. In the present study, it is shown that
the two-temperature Navier--Stokes system with relaxation terms of $O(1)$
can also be derived from the Boltzmann equation \eqref{BE-1} with~\eqref{col-op}
for a particular collision kernel~\eqref{Wtheta-2}, \eqref{Ws-Wr-2}, and \eqref{sigma-model}.

\begin{remark}
In order to calculate the first-order coefficient $\mathcal{K}$ in the relaxation
terms of the equations~\eqref{NS-2-c} and \eqref{NS-2-d}, one has to obtain the
solution $\widetilde{D}$ to the integral equation \eqref{eq-D}. This may be harder
than obtaining $\mathcal{A}$, $\mathcal{B}$, and $\mathcal{C}$ [see the paragraph
before Remark \ref{rem:2}] due to the complexity of the right-hand side
of equation \eqref{eq-D}. Nevertheless, it should be possible, in principle, to obtain
$\widetilde{D}$ numerically or approximately.

\end{remark}

\begin{remark}
If the two-temperature Navier--Stokes system derived in \cite{ABGK-20} is
compared with the system~\eqref{NS-2}, the difference is as follows.
In the former system,  
$\Lambda_\mathrm{int}^\mathrm{tr} (\rho, T_\mathrm{tr}, T_\mathrm{int})
= \Lambda_\mathrm{tr}^\mathrm{int} (\rho, T_\mathrm{tr}, T_\mathrm{int}) = 0$,
that is, the cross-diffusion terms disappear as in the system~\eqref{NS-alpha=0}.
In addition,  the $O (\epsilon)$ term $\mathcal{K} (\rho, T_\mathrm{tr}, T_\mathrm{int})$
is identically zero. Furthermore, $\Lambda_\mu (\rho, T_\mathrm{tr}, T_\mathrm{int})$,
$\Lambda_\mathrm{tr}^\mathrm{tr} (\rho, T_\mathrm{tr}, T_\mathrm{int})$, and
$\Lambda_\mathrm{int}^\mathrm{int} (\rho, T_\mathrm{tr}, T_\mathrm{int})$ are
simple and explicit functions of $T_\mathrm{tr}$ and $T_\mathrm{int}$.
Thanks to the simplicity, the system derived in \cite{ABGK-20} has been
successfully applied to the problem of shock-wave structure \cite{ABGK-20}, and
its boundary conditions have been derived \cite{KABGM-21}. 
\end{remark}

\begin{remark}
A two-temperature fluid model at the Navier--Stokes level, corresponding to the system~\eqref{NS-2}, is also discussed in
 \cite{BG-11,BG-22}, with a scaling corresponding to \eqref{scale-2},
on the basis of the Boltzmann equation with
discrete energy variables for the internal modes. However, the Boltzmann equation is presented
in a more general and abstract form, and the forms of the transport coefficients
corresponding to $\Lambda_\mu$, $\Lambda_\mathrm{tr}^\mathrm{tr}$, etc.
are not explicitly shown. Furthermore, it is not clear if the first-order source term
corresponding to $\big( I, Q_\mathrm{s} (M_\mathrm{r}, f^{(1)}) \big)$ is proportional
to $T_\mathrm{tr} - T_\mathrm{int}$ or vanishes.

\end{remark}

\section{Concluding remarks}\label{sec:concluding}

In the present paper, we focus our attention on the systematic derivation of 
fluid-dynamic equations with two temperatures, i.e., translational temperature $T_\text{tr}$
and internal one $T_\text{int}$, and with relaxation terms, from the Boltzmann equation for
a polyatomic gas. It was a common understanding that such fluid equations hold
when the interaction
between the translational and internal modes is weak, that is, when resonant (or elastic)
collisions occur much more frequently than standard (or inelastic) collisions. 
In order to describe this situation,
we proposed a Boltzmann-type model in which the collision kernel is a linear combination
of a resonant collision kernel with coefficient $1 - \theta$ and a standard
collision kernel with coefficient $\theta$, where $\theta$ is a parameter ($0 \le \theta \le 1$). 
Furthermore, we adopted specific forms of collision kernels for both resonant and standard
collisions. These collision kernels were chosen mainly for mathematical convenience rather
than physical realism. Then, using the Chapman--Enskog expansion, we performed
a systematic analysis for small $\theta$ and for small Knudsen numbers Kn.

First, we consider the case when $\theta$ is of the order of $\text{Kn}^2$, that is,
the interaction between the translational and internal modes is very weak.
In this case, an Euler system without interaction between the translational and
internal modes is obtained at the leading order, and a two-temperature
Navier--Stokes system with relaxation terms proportional to $T_\text{tr} – T_\text{int}$
is obtained at the first order in Kn. In this system, the relaxation terms, 
viscosity terms, and heat-conduction terms are all of the order of Kn.
Moreover, the coefficients of
the relaxation terms and the transport coefficients are expressed in terms of
the parameters included in the assumed collision kernels.

The case we consider next is when $\theta$ is of the order of $\text{Kn}$,
that is, the interaction between the translational and internal modes is still weak,
but not extremely weak. In this case, at the leading order, one obtains an Euler system
with relaxation terms proportional to $T_\text{tr} – T_\text{int}$, through
which the internal modes interact with the translational mode. At the order of Kn,
a two-temperature Navier—Stokes system, similar to that derived for
$\theta = O(\text{Kn}^2)$, is obtained. The difference is that
the relaxation terms in this case include $O(1)$ terms as well as $O(\text{Kn})$
terms, both being proportional to $T_\text{tr} - T_\text{int}$.
It had been known that this type of Navier—Stokes equations
[with $O(1)$ relaxation terms]
could be derived from model kinetic equations such as the ES model by an
appropriate parameter setting \cite{ABGK-20}. However, it was far from obvious whether a similar system of equations could be derived explicitly from the Boltzmann equation.
The present study provides a positive answer to this question, even though the used collision operator is a particular model.

It would be worthwhile to apply the current two-temperature Navier—Stokes system, in both cases of $\theta=O(\text{Kn}^2)$ and $\theta = O(\text{Kn})$, to some
fundamental problems, such as the problem of shock-wave structure \cite{ABGK-20}.
It would also be interesting to consider different types of collision operators
and to see if the same type of two-temperature Navier—Stokes system can
be derived from them. These will be topics of future research.

\section*{Acknowledgment}

This work was initiated while K.A.~was visiting Karlstad University. The authors are
grateful for the travel support from the Royal Swedish Academy of Sciences, grant MA2020-0007.

\appendix

\section{Proof of the bounds \eqref{cocfa}}\label{sec:nu-bounds}

Let us first note that
\begin{equation*}
\nu _{\theta }=4\pi C_{\text{r}}\frac{\Gamma ^{2}\left( \delta /2\right) }
{\Gamma \left( \delta \right) }\int_{\mathbb{R}^{3}\times \mathbb{R}_{+}}\!
\left( I+I_{\ast }\right) ^{\alpha }\left\vert \bm{\xi }-\bm{\xi }_{\ast } \right\vert ^{\beta }
M_{\ast }\,\text{d}\bm{\xi }_{\ast }\,\text{d}I_{\ast },
\end{equation*}
since
\begin{eqnarray*}
\nu _{\text{s}} &=& \int_{[0,1]^{2}\times \mathbb{S}^{2}\times \mathbb{R}^{3}\times \mathbb{R}_{+}}
\!C_{\text{s}}\left( I+I_{\ast }\right) ^{\alpha}\left\vert \bm{\xi }-\bm{\xi }_{\ast }\right\vert ^{\beta}
M_{\ast }
\\
&& \times R^{\left( \beta +1\right) /2}\left( 1-R\right) ^{\delta +\alpha-1}
[r(1-r)]^{\delta /2-1} \text{d}R\,\text{d}r\,\text{d}\bm{\sigma }\,
\text{d}\bm{\xi }_{\ast }\,\text{d}I_{\ast }
\\
&=&4\pi C_{\text{r}}\frac{\Gamma ^{2}\left( \delta /2\right) }{\Gamma \left(\delta \right) }
\int_{\mathbb{R}^{3}\times \mathbb{R}_{+}}\!\left( I+I_{\ast}\right) ^{\alpha }
\left\vert \bm{\xi }-\bm{\xi }_{\ast}\right\vert^{\beta }
M_{\ast }\,\text{d}\bm{\xi }_{\ast }\,\text{d} I_{\ast },
\end{eqnarray*}
and
\begin{eqnarray*}
\nu _{\text{r}} &=&\int_{[0,1]\times \mathbb{S}^{2}\times \mathbb{R}^{3}\times \mathbb{R}_{+}}
\!C_{\text{r}}\left( I+I_{\ast }\right) ^{\alpha}
\left\vert \bm{\xi }-\bm{\xi }_{\ast }\right\vert ^{\beta}
M_{\ast }[r(1-r)]^{\delta /2-1}\text{d}r\,\text{d}\bm{\sigma }\,
\text{d}\bm{\xi }_{\ast }\,\text{d}I_{\ast }
\\
&=&4\pi C_{\text{r}}\frac{\Gamma ^{2}\left( \delta /2\right) }{\Gamma \left(\delta \right) }
\int_{\mathbb{R}^{3}\times \mathbb{R}_{+}}\!\left( I+I_{\ast}\right) ^{\alpha }
\left\vert \bm{\xi }-\bm{\xi }_{\ast}\right\vert ^{\beta }
M_{\ast }\,\text{d}\bm{\xi }_{\ast }\,\text{d}I_{\ast }.
\end{eqnarray*}
Moreover, it is clear that
\begin{align*}
\int_{\substack{ \mathbb{R}^{3}\times \mathbb{R}_{+} \\ \left\vert 
\bm{\xi }_{\ast }\right\vert \leq 1/2}}\!I_{\ast }^{s_{1}}\left\vert 
\bm{\xi }_{\ast }\right\vert ^{s_{2}}M_{\ast }\,\text{d}\bm{\xi }_{\ast }
\,\text{d}I_{\ast } &= C_\text{g}>0,
\\
\int_{\substack{ \mathbb{R}^{3}\times \mathbb{R}_{+} \\ \left\vert 
\bm{\xi }_{\ast }\right\vert \geq 2}}\!I_{\ast }^{s_{1}}\left\vert 
\bm{\xi }_{\ast }\right\vert ^{s_{2}}M_{\ast }\,\text{d}\bm{\xi }_{\ast }
\,\text{d}I_{\ast } &= C_\text{g}>0,
\\
\int_{\mathbb{R}^{3}\times \mathbb{R}_{+}}\!\left( 1+I_{\ast }\right)^{s_{1}}
\left( 1+\left\vert \bm{\xi }_{\ast }\right\vert \right)^{s_{2}}
M_{\ast }\,\text{d}\bm{\xi }_{\ast }\,\text{d}I_{\ast }
&= C_\text{g}<\infty \;\;\;\; (\text{for any }s_{1},s_{2}\geq 0),
\end{align*}
where $C_\text{g}$ denotes a generic constant.

The bounds for the collision frequency now follow by the following estimates
\begin{align*}
\!\left( I+I_{\ast }\right) ^{\alpha } &\leq \left( 1+I\right) ^{\alpha
}\left( 1+I_{\ast }\right) ^{\alpha },
\\
\left\vert \bm{\xi }-\bm{\xi }_{\ast }\right\vert ^{\beta }
&\leq \left( \left\vert \bm{\xi }\right\vert +\left\vert 
\bm{\xi }_{\ast }\right\vert \right) ^{\beta }\leq \left(
1+\left\vert \bm{\xi }\right\vert \right) ^{\beta }\left(
1+\left\vert \bm{\xi }_{\ast }\right\vert \right) ^{\beta },
\end{align*}
for the upper bound, and 
\begin{align*}
\left( I+I_{\ast }\right) ^{\alpha } &\geq \left\{ 
\begin{array}{l}
I^{\alpha }\geq \left( 1/2\right) ^{\alpha }\left( 1+I\right) ^{\alpha }%
\;\;\;\; (\text{if}\;\, I\geq 1),
\\[1mm]
I_{\ast }^{\alpha }\geq \left( I_{\ast }/2\right) ^{\alpha }\left(
1+I\right) ^{\alpha }\;\;\;\; (\text{if}\;\, I\leq 1),
\end{array}%
\right. 
\\[2mm]
\left\vert \bm{\xi }-\bm{\xi }_{\ast }\right\vert ^{\beta }
& \geq \left\vert \left\vert \bm{\xi }\right\vert -\left\vert 
\bm{\xi }_{\ast }\right\vert \right\vert ^{\beta }
\\
& \geq \left\{ 
\begin{array}{l}
\left( 1/2\right) ^{\beta }\left\vert \bm{\xi }\right\vert ^{\beta
}\geq \left( 1/4\right)^{\beta } \left( 1+\left\vert \bm{\xi }\right\vert
\right) ^{\beta } \;\;\; (\text{for}\;\, \left\vert \bm{\xi }_{\ast }\right\vert
\leq 1/2\;\; \text{if}\; \left\vert \bm{\xi }\right\vert \geq 1),
\\ [1mm]
\left( \left\vert \bm{\xi }_{\ast }\right\vert/2\right) ^{\beta }\geq \left(
\left\vert \bm{\xi }_{\ast }\right\vert /4\right) ^{\beta }\left(
1+\left\vert \bm{\xi }\right\vert \right) ^{\beta }\;\;\;
(\text{for}\;\, \left\vert \bm{\xi }_{\ast }\right\vert \geq 2\;\; \text{if}\; \left\vert 
\bm{\xi }\right\vert \leq 1),
\end{array}%
\right.
\end{align*}
for the lower bound.

\section{Calculation of $Q_\text{s} (M_\text{r}, M_\text{r})$ and
$( I, Q_\text{s} (M_\text{r}, M_\text{r}) )$ }\label{sec:source}

Using identities~\eqref{col-op}, \eqref{Ws-Wr-2-a}, and \eqref{Qs-Qr},
we have the following expression of $Q_\text{s} (M_\text{r}, M_\text{r})$:
\begin{align}\label{Qs-Mr-Mr-1}
Q_\text{s} (M_\text{r}, M_\text{r}) &= \int_{\left( \mathbb{R}^{3}\times \mathbb{R}_{+} \right)^{3}}
\sigma_\text{s} \frac{\left\vert \bm{g}\right\vert }{\left\vert \bm{g}^{\prime } \right\vert }
\frac{m}{2} \left[ M_{\text{r}}' M_{\text{r}*}'
\left( \frac{II_{\ast }}{I^{\prime }I_{\ast }^{\prime }}%
\right) ^{\delta /2-1}- M_\text{r} M_{\text{r}*} \right]
\nonumber \\
& \qquad \qquad
\times \bm{\delta} _{3}\left( \mathbf{G}-\mathbf{G}^{\prime }\right) \bm{\delta}
_{1}\left( E-E^{\prime }\right) \,\text{d}\boldsymbol{\xi }_{\ast }\text{d}\boldsymbol{\xi 
}^{\prime }\text{d}\boldsymbol{\xi }_{\ast }^{\prime }\text{d}I_{\ast }\text{d}I^{\prime
}\text{d}I_{\ast }^{\prime }.
\end{align}
Let $M_\text{s}$ [see~\eqref{Ms}] with
$T = T_\text{tr}$ be denoted by $M_\text{s}^\text{tr}$, i.e.,
\begin{equation*}
M_\mathrm{s}^\text{tr} = \frac{n I^{\delta/2-1}}{(2\pi k_\mathrm{B} T_\text{tr}/m)^{3/2}
(k_\mathrm{B} T_\text{tr})^{\delta/2}
\Gamma (\delta/2)}
\exp \left( - \frac{m |\bm{\xi} - \bm{u}|^2 + 2I}{2k_\mathrm{B} T_\text{tr}}
\right).
\end{equation*}
Then, we have the relation
\begin{equation*}
\frac{M_\text{r}}{M_\text{s}^\text{tr}}=\dfrac{T_\text{tr}^{\delta /2}}{T_\text{int}^{\delta /2}}
\exp \left( - \frac{I}{ k_\text{B}T_\text{int} } + \frac{I}{ k_\text{B}T_\text{tr} } \right)
=\dfrac{T_\text{tr}^{\delta /2}}{T_\text{int}^{\delta /2}}e^{-\eta I/\zeta },
\end{equation*}
with
\begin{equation}\label{zeta-eta}
\zeta =k_\text{B}\frac{T_\text{tr}T_\text{int}}{\left\vert T_\text{tr}-T_\text{int}\right\vert }, \qquad
\eta =\frac{T_\text{tr}-T_\text{int}}{\left\vert T_\text{tr}-T_\text{int}\right\vert }
=\left\{ 
\begin{array}{c}
~\text{~}1 \;\;\; \text{ if }\; T_\text{tr}>T_\text{int}, \\ 
-1 \;\;\; \text{ if }\;T_\text{tr}<T_\text{int}.
\end{array}
\right.
\end{equation}
Moreover, since
\begin{equation*}
\frac{(M_\text{s}^\text{tr})^\prime (M_\text{s}^\text{tr})_*^\prime}
{\left( I^{\prime }I_{\ast}^{\prime }\right)^{\delta /2-1}}
= \frac{M_\text{s}^\text{tr} M_{s\ast }^\text{tr}}{\left( II_{\ast}\right)^{\delta /2-1}},
\end{equation*}
holds, it follows that
\begin{align}\label{relation}
M_\text{r}' M_{\text{r}*}'
\left( \frac{II_{\ast }}{I^{\prime}
I_{\ast }^{\prime }}\right) ^{\delta /2-1}- M_\text{r} M_{\text{r}*}
& =
M_\text{s}^\text{tr} M_{\text{s}\ast }^\text{tr} \left( \frac{M_\text{r}'}{(M_\text{s}^\text{tr})^{\prime }}
\frac{M_{\text{r}*}'}{ (M_\text{s}^\text{tr} )_\ast^{\prime }}-\frac{M_\text{r}}{M_\text{s}^\text{tr}}
\frac{M_{\text{r}*}}{M_{\text{s}\ast }^\text{tr}}\right)
\nonumber \\
& = M_\text{s}^\text{tr} M_{\text{s}\ast }^\text{tr} T_\text{tr}^{\delta }T_\text{int}^{-\delta }\left( e^{-\eta \left(
I^{\prime }+I_{\ast }^{\prime }\right) /\zeta }-e^{-\eta \left( I+I_{\ast
}\right) /\zeta }\right) \text{.}
\end{align}

Now, a series of changes of integration variables is performed. More specifically,

\begin{itemize}

\item
$(\bm{\xi}_*,\, \bm{\xi}',\, \bm{\xi}_*',\, I_*,\, I',\, I_*') \to  
(\bm{\xi}_*,\, \bm{g}',\, \bm{G}',\, I_*,\, r,\, s)$ with the help of
$\bm{g}' = \bm{\xi}' - \bm{\xi}_*'$, $\bm{G}' =(\bm{\xi}' + \bm{\xi}_*')/2$,
$r=I'/(I' + I_*')$ [cf.~relations~\eqref{BL-r}], and $s=(I' + I_*')/\zeta$; 

\item
$(\bm{\xi}_*,\, \bm{g}',\, \bm{G}',\, I_*,\, r,\, s) \to  
(\bm{\xi}_*,\, |\bm{g}'|,\, \bm{\sigma},\, \bm{G}',\, I_*,\, r,\, s)$
with the help of $\bm{\sigma} = \bm{g}'/|\bm{g}'|$ (spherical coordinates
for $\bm{g}'$);

\item
$(\bm{\xi}_*,\, |\bm{g}'|,\, \bm{\sigma},\, \bm{G}',\, I_*,\, r,\, s) \to
(\bm{\xi}_*,\, w,\, \bm{\sigma},\, \bm{G}',\, I_*,\, r,\, s)$ with
the help of $w=m|\bm{g}'|^2/4\zeta$.

\end{itemize}

\noindent
The calculation of the Jacobian at each step leads to
\begin{align*}
\text{d}\bm{\xi}_* \text{d}\bm{\xi}' \text{d}\bm{\xi}_*'
\text{d} I_* \text{d} I' \text{d} I_*'
& = |\bm{g}'|^2 \zeta^2 s\, \text{d}\bm{\xi}_* \text{d} I_* \text{d} |\bm{g}'|
\text{d}\bm{\sigma} \text{d}\bm{G}'\text{d}r \text{d}s
\\
& = \frac{4}{m^{3/2}} \zeta^{7/2} s \sqrt{w}\, \text{d}\bm{\xi}_* \text{d} I_*
\text{d}w\, \text{d}\bm{\sigma} \text{d}\bm{G}'\text{d}r \text{d}s,
\end{align*}
and the domain of integration in the variables
$(\bm{\xi}_*,\, w,\, \bm{\sigma},\, \bm{G}',\, I_*,\, r,\, s)$
becomes $\bm{\xi}_* \in \mathbb{R}^3$,\, $w \in \mathbb{R}_+$, $\bm{\sigma} \in \mathbb{S}^2$,
$\bm{G}' \in \mathbb{R}^3$, $I_* \in \mathbb{R}_+$,\, $r \in [0, 1]$, and $s \in \mathbb{R}_+$.
Here, we introduce some additional variables for later convenience:
\begin{align}\label{new-var}
\begin{aligned}
& r' = \frac{I}{I + I_*}, \quad v = \frac{I + I_*}{\zeta}, \quad u=\frac{m|\bm{g}|^2}{4\zeta},
\quad \vartheta_\text{int} = \frac{k_\text{B} T_\text{tr}}{\zeta}
= \frac{|T_\text{tr} - T_\text{int}|}{T_\text{int}},
\\
& \widetilde{\bm{\xi}} = \bm{\xi}-\bm{u}, \quad \widetilde{\bm{\xi}}_* = \bm{\xi}_*-\bm{u}, \quad
\bm{g} = \widetilde{\bm{\xi}}-\widetilde{\bm{\xi}}_* =\bm{\xi}-\bm{\xi}_*, \quad
\bm{\sigma}' = \frac{\bm{g}}{|\bm{g}|},
\\
& \widetilde{\bm{G}} = \frac{\widetilde{\bm{\xi}}+\widetilde{\bm{\xi}_*}}{2}, \quad
\widehat{\bm{G}} = \frac{\sqrt{m}\, \widetilde{\bm{G}}}{\sqrt{k_\text{B} T_\text{tr}}}.
\end{aligned}
\end{align}
Then, we have 
\begin{align*}
\frac{|\bm{g}|}{|\bm{g}'|} = \frac{\sqrt{u}}{\sqrt{w}}, \qquad E = \zeta(u + v), \qquad
E' = \zeta(w + s).
\end{align*}
In consequence,~\eqref{Qs-Mr-Mr-1} is transformed as
\begin{align}\label{Qs-Mr-Mr-2}
Q_{s} (M_\text{r}, M_\text{r}) &= \frac{2}{\sqrt{m}} \zeta^{7/2} \frac{T_\text{tr}^\delta}{T_\text{int}^\delta}
\int_{\mathbb{S}^2} \text{d} \bm{\sigma} \cdot
\int_{\mathbb{R}^3} \bm{\delta}_3 (\bm{G} - \bm{G}') \text{d} \bm{G}'
\nonumber \\
& \quad \;\; \times
\int_{\mathbb{R}^3 \times \mathbb{R}_+^3 \times [0, 1]} \sigma_\text{s} s \sqrt{u}\,
M_\text{s}^\text{tr} M_{\text{s} *}^\text{tr} (e^{-\eta s} - e^{-\eta v})
\nonumber \\
& \qquad \qquad \qquad \times
\bm{\delta}_1 \big( \zeta (u + v) - \zeta (w + s) \big)
\text{d}\bm{\xi}_* \, \text{d} I_*\, \text{d}s\, \text{d}w\, \text{d} r.
\nonumber \\
& = \frac{ m^{5/2} \zeta^{\delta + 1/2} n^2}
{\pi^2 k_\text{B}^{\delta+3} T_\text{tr}^3 T_\text{int}^\delta \Gamma^2 (\delta/2)}
\nonumber \\
& \qquad \times
\int_{\mathbb{R}^3 \times \mathbb{R}_+^3 \times [0, 1]}
\sigma_\text{s} s\sqrt{u}\, [r' (1-r')]^{\delta/2 -1} v^{\delta-2}
\bm{\delta}_1 (u + v - w -s)
\nonumber \\
& \qquad \qquad \qquad \qquad \times
e^{- |\widehat{\bm{G}}|^2} e^{-(u+v)/\vartheta_\text{int}} (e^{-\eta s} - e^{-\eta v})\,
\text{d} \bm{\xi}_* \text{d}I_* \text{d}s\, \text{d}w\, \text{d}r.
\end{align}

Now, we consider the model \eqref{sigma-model-a} for $\sigma_\text{s}$. With some
rearrangement and then with some new variables, it can be rewritten as follows:
\begin{align}\label{scs}
\sigma _\text{s} &=C_\text{s}
\frac{\left\vert \bm{g}\right\vert ^{\beta -1}}{E^{2}}
\left( I+I_{\ast }\right)^{\alpha }\left( \frac{m\left\vert \bm{g}^{\prime }\right\vert ^{2}}
{4E}\right) ^{\left( \beta +1\right) /2}
\left( \frac{I^{\prime }+I_{\ast }^{\prime }}{E}\right) ^{\alpha }
\left( \frac{I^{\prime }}{E}\frac{I_{\ast }^{\prime }}{E}\right) ^{\delta /2-1}
\notag \\
&=C_\text{s}\frac{\left( 4\zeta u/m\right) ^{\left( \beta -1\right) /2}}
{\left[ \zeta \left( u+v\right) \right] ^{2}}\left( \zeta v\right) ^{\alpha }
\left( \frac{w}{u+v}\right) ^{\left( \beta +1\right) /2}\!\!
\left( \frac{s}{u+v}\right)^{\delta +\alpha -2}\!\!
\left( \frac{I^{\prime }}{I^{\prime }+I_{\ast }^{\prime }}\cdot
\frac{I_{\ast }^{\prime }}{I^{\prime }+I_{\ast }^{\prime }}\right) ^{\delta/2-1}
\notag \\
&=\frac{C_\text{s}}{\zeta ^{(5-\beta )/2-\alpha }}
\left( \frac{4}{m}\right)^{\left( \beta -1\right) /2}
\dfrac{u^{\left( \beta -1\right) /2}v^{\alpha}s^{\delta +\alpha -2} w^{\left( \beta +1\right) /2}}
{\left( u+v\right)^{\delta +\alpha +\left( \beta +1\right) /2}}
\left[ r\left( 1-r\right) \right] ^{\delta /2-1},
\end{align}
for $\left( \alpha ,\beta \right) \in \left[ 0,\delta /2\right) \times \left[ 0,1\right]$,
where
\begin{align*}
C_\text{s}=\frac{\Gamma \left( \delta +\alpha +\left( \beta +3\right) /2\right) }{%
\Gamma \left( \left( \beta +3\right) /2\right) \Gamma \left( \delta +\alpha
\right) }C_\text{r}.
\end{align*}

Let us consider
\begin{align*}
\big( I,\, Q_\text{s} (M_\text{r}, M_\text{r}) \big) = \int_{\mathbb{R}^3\times\mathbb{R}_+} 
Q_\text{s} (M_\text{r}, M_\text{r}) I \text{d}\bm{\xi} \text{d}I.
\end{align*}
We substitute expressions~\eqref{Qs-Mr-Mr-2} and \eqref{scs} into the above
equation and carry out a series of changes of integration variables
using some new variables defined by relations~\eqref{new-var}. That is,

\begin{itemize}

\item
$(\bm{\xi},\, \bm{\xi}_*,\, I,\, I_*,\, s,\, w,\, r) \to
(\widetilde{\bm{\xi}},\, \widetilde{\bm{\xi}}_*,\, r',\, v,\, s,\, w,\, r)$;

\item
$(\widetilde{\bm{\xi}},\, \widetilde{\bm{\xi}}_*,\, r',\, v,\, s,\, w,\, r) \to
(\bm{g},\, \widetilde{\bm{G}},\,  r',\, v,\, s,\, w,\, r)$;

\item
$(\bm{g},\, \widetilde{\bm{G}},\,  r',\, v,\, s,\, w,\, r) \to
(|\bm{g}|,\, \bm{\sigma}',\, \widetilde{\bm{G}},\,  r',\, v,\, s,\, w,\, r)$;

\item
$(|\bm{g}|,\, \bm{\sigma}',\, \widetilde{\bm{G}},\,  r',\, v,\, s,\, w,\, r) \to
(u,\, \bm{\sigma}',\, \widehat{\bm{G}},\, r',\, v,\, s,\, w,\, r)$.

\end{itemize}

\noindent
As the result, we obtain
\begin{align*}
\text{d}\bm{\xi}\, \text{d}\bm{\xi}_*\, \text{d}I\, \text{d}I_*\, \text{d}s\, \text{d}w\, \text{d}r
& = |\bm{g}|^{2} \text{d} |\bm{g}|\, \text{d} \bm{\sigma}' \text{d}\widetilde{\bm{G}} \cdot
\zeta^2 v\,\text{d}r' \text{d}v \cdot \text{d}s\, \text{d}w\, \text{d}r
\\
& = \frac{4}{m^{3}} \zeta^{7/2} \left( k_\text{B} T_\text{tr} \right)^{3/2} \sqrt{u} v\,
\text{d}r' \text{d}\widehat{\bm{G}}\,\text{d}\bm{\sigma}' \text{d}u\, \text{d}v\,
\text{d}s\, \text{d}w\, \text{d}r, 
\end{align*}
and the domain of integration in the variables $(r',\, \widehat{\bm{G}},\, \bm{\sigma}',\, u,\, v,\, s,\, w,\, r)$
is as follows: $r' \in [0, 1]$, $\widehat{\bm{G}} \in \mathbb{R}^3$, $\bm{\sigma}' \in \mathbb{S}^2$,
$u \in \mathbb{R}_+$, $v \in \mathbb{R}_+$, $s \in \mathbb{R}_+$, $w \in \mathbb{R}_+$, and
$r \in [0, 1]$.

With these changes of variables, the following expression is obtained:
\begin{align}\label{source-term}
& \big( I, Q_\text{s} (M_\text{r}, M_\text{r}) \big)
\nonumber \\
& \;\; = \dfrac{4^{\left( \beta +1\right) /2} n^{2} \zeta^{\delta +\alpha +(\beta +5)/2} C_\text{s}}
{m^{\beta /2}\pi^{2}k_\text{B}^{\delta+3/2} T_\text{tr}^{3/2} T_\text{int}^{\delta }}
\int_{\mathbb{S}^{2}} \text{d} \bm{\sigma}' \int_{\mathbb{R}^{3}}e^{-| \widehat{\bm{G}} |^{2}}
\text{d}\widehat{\bm{G}}
\nonumber \\
& \;\; \quad \times \dfrac{1}{\Gamma^{2}\left( \delta /2\right) }
\int_{0}^{1} r^{\delta/2-1} \left( 1-r \right)^{\delta /2-1} \text{d}r
\int_{0}^{1} \left( r' \right)^{\delta /2} \left( 1-r' \right)^{\delta /2-1} \text{d}r'
\nonumber \\
& \;\; \quad \times \int_{\mathbb{R}_{+}^{4}} \frac{\left( uw \right)^{\left( \beta+1\right) /2}
s^{\delta +\alpha -1} v^{\delta +\alpha }}
{\left( u+v \right)^{\delta +\alpha +\left( \beta +1\right) /2}}
e^{-\left( u+v \right) /\vartheta _\text{int}} \left( e^{-\eta s}-e^{-\eta v} \right)
\nonumber \\
& \qquad \qquad \qquad \qquad \qquad \qquad \qquad \qquad \times
\bm{\delta} _{1} \left( u+v-w-s \right) \text{d}s\, \text{d}w\, \text{d}u\, \text{d}v
\nonumber \\
& \;\; =\frac{4^{\left( \beta +3 \right) /2}}{m^{\beta /2}} C_\text{s} n^{2}\sqrt{\pi }
\frac{\zeta^{\delta +\alpha +(\beta +5)/2}}{k_\text{B}^{\delta+3/2}
T_\text{tr}^{3/2} T_\text{int}^{\delta}}
\dfrac{\Gamma \left( \delta /2 \right) \Gamma \left( \delta /2+1 \right) }
{\Gamma \left( \delta \right) \Gamma \left( \delta +1 \right) }
\nonumber \\
& \times \int_{0}^{\infty} \!\! \int_{0}^{\infty} \!\! \int_{0}^{u+v}
\frac{\left[ u \left( u+v-s \right) \right]^{\left( \beta +1 \right) /2}}
{\left( u+v \right)^{\delta +\alpha +\left( \beta +1\right) /2}}\,
s^{\delta +\alpha -1} v^{\delta+\alpha }
e^{-\left( u+v \right) /\vartheta _\text{int}} 
\nonumber \\
& \qquad \qquad \qquad \qquad \qquad \qquad \qquad \qquad \times\left( e^{-\eta s}-e^{-\eta v} \right)
\text{d}s\, \text{d}u\, \text{d}v
\nonumber \\
& \;\; =\frac{2^{\beta +2}}{m^{\beta /2}} C_\text{s} n^{2} \sqrt{\pi }
\frac{\zeta^{\delta +\alpha +(\beta +5)/2}}{k_\text{B}^{\delta+3/2} T_\text{tr}^{3/2} T_\text{int}^{\delta}}
\dfrac{\Gamma^{2} \left( \delta /2 \right)}{\Gamma^{2} \left( \delta \right) }\, \Omega,
\end{align}
where $\Omega$ is expressed as
\begin{align}\label{Omega}
\Omega = \int_{0}^{\infty} \!\! \int_{v}^{\infty} \!\! \int_{0}^{q}
\frac{\left( q-s \right)^{\left( \beta +1 \right) /2} \left( q-v \right)^{\left( \beta+1 \right)/2}}
{q^{\delta +\alpha +\left( \beta +1 \right) /2}}\, s^{\delta+\alpha -1} v^{\delta +\alpha }
e^{-q/\vartheta _\text{int}} \nonumber \\
\times \left( e^{-\eta s}-e^{-\eta v} \right)
\text{d}s\, \text{d}q\, \text{d}v,
\end{align}
after changing the integration variables from $(s, u, v)$ to $(s, q, v)$ with $q=u+v$.

We carry out further transformation of $\Omega$. By changing the order of integrations
with respect to $q$ and $v$, it can be expressed in the following form:
\begin{align*}
\Omega = \int_0^\infty \!\! \int_0^q \!\! \int_0^q
F (s, v, q) \left( \frac{e^{-\eta s}}{s} - \frac{e^{-\eta v}}{s} \right)\, \text{d}s\, \text{d}v\, \text{d}q,
\end{align*}
where
\begin{align*}
F (s, v, q) 
= \frac{\left( q-s \right)^{\left( \beta +1 \right) /2} \left( q-v \right)^{\left( \beta+1 \right)/2}}
{q^{\delta +\alpha +\left( \beta +1 \right) /2}}\, s^{\delta+\alpha} v^{\delta +\alpha }
e^{-q/\vartheta _\text{int}}. 
\end{align*}
By changing the labels of the integration variables $(s, v)$ to $(v, s)$ and changing
the order of integrations with respect to $v$ and $s$, we have
\begin{align*}
\int_0^\infty \!\! \int_0^q \!\! \int_0^q F (s, v, q) \frac{e^{-\eta s}}{s} \, \text{d}s\, \text{d}v\, \text{d}q
= \int_0^\infty \!\! \int_0^q \!\! \int_0^q F (v, s, q) \frac{e^{-\eta v}}{v} \, \text{d}s\, \text{d}v\, \text{d}q.
\end{align*}
Since $F (v, s, q)  = F(s, v, q)$, it follows that
\begin{align*}
\Omega & = \int_0^\infty \!\! \int_0^q \!\! \int_0^q
F (s, v, q) \left( \frac{e^{-\eta v}}{v} - \frac{e^{-\eta v}}{s} \right)\, \text{d}s\, \text{d}v\, \text{d}q
\\
& = \int_0^\infty \!\! \int_0^q \left[ \int_0^q
 \left( q-s \right)^{\left( \beta +1 \right) /2} 
\,(s - v) s^{\delta+\alpha-1} 
\text{d}s\, \right] 
\frac{v^{\delta +\alpha-1 } \left( q-v \right)^{\left( \beta+1 \right)/2}}{q^{\delta +\alpha +\left( \beta +1 \right) /2}}
\nonumber \\
& \qquad \qquad \qquad \qquad \qquad \qquad \qquad \qquad \qquad \qquad \qquad \times
e^{-q/\vartheta _\text{int}} e^{-\eta v}\,  \text{d}v\, \text{d}q.
\end{align*}
Let us consider the following integral, which is part of the above expression of $\Omega$:
\begin{align*}
\mathcal{J} = \frac{1}{q^{\delta + \alpha + (\beta+1)/2}} \int_0^q
(q-s)^{(\beta+1)/2} (s-v)\, s^{\delta+\alpha-1} \text{d}s.
\end{align*}
Letting $\widetilde{s}=s/q$ and using the definition of the beta function
$B(x, y) = \int_0^1 t^{x-1} (1-t)^{y-1} \text{d}t$ and its relation to the
gamma function $B(x, y) = \Gamma(x) \Gamma(y)/\Gamma(x+y)$, we have
\begin{align*}
\mathcal{J} & =\int_{0}^{1}\left( 1-\widetilde{s}\,\right)^{\left( \beta +1\right)/2}
\left( q\, \widetilde{s}\,^{\delta +\alpha}-v\, \widetilde{s}\,^{\delta +\alpha-1}\right)
\text{d}\,\widetilde{s}
\\
& = \widetilde{C} \left[ q\dfrac{\Gamma \left( \delta +\alpha +1\right) }
{\Gamma \left( \delta +\alpha \right) }
-v\dfrac{\Gamma \big( \delta +\alpha +\left( \beta +5\right) /2 \big) }
{\Gamma \big( \delta +\alpha +\left( \beta +3\right) /2 \big) }\right]
\\
& = \widetilde{C} \left[ \left( \delta +\alpha \right) \left( q-v\right) -%
\frac{\beta +3}{2}v\right],
\end{align*}
with
\begin{align*}
\widetilde{C} =\frac{\Gamma ( \delta +\alpha )
\Gamma \big( \left( \beta +3\right) /2 \big) }
{\Gamma \big( \delta+\alpha +\left( \beta +5\right) /2 \big) }.
\end{align*}
With this expression of $\mathcal{J}$, further transformation of $\Omega$
can be made as follows:
\begin{align*}
\Omega
& = \widetilde{C} \int_{0}^{\infty} \!\! \int_{0}^{q}
\left[ \left( \delta +\alpha \right) (q-v)^{\left( \beta +3\right)/2}
v^{\delta +\alpha -1}
-\frac{\beta +3}{2} (q-v)^{\left( \beta +1 \right)/2} v^{\delta +\alpha} \right]
\\
& \qquad \qquad \qquad \qquad \qquad \qquad \qquad \qquad 
\qquad \qquad \qquad \times
e^{-q/\vartheta _\text{int}} e^{-\eta v} \text{d}v\, \text{d}q
\\
& = \widetilde{C} \int_{0}^{\infty} \!\! \int_{0}^{q}
\frac{\partial}{\partial v}
\left[ (q-v)^{\left( \beta +3\right) /2}v^{\delta +\alpha} \right]
e^{-q/\vartheta_\text{int}} e^{-\eta v} \text{d}v\, \text{d}q
\\
& = \eta\, \widetilde{C} \int_{0}^{\infty} \!\! \int_0^q
(q-v)^{\left( \beta+3 \right)/2} v^{\delta +\alpha} e^{-q/\vartheta_\text{int}}
e^{-\eta v} \text{d}v\, \text{d}q
\\
& = \eta\, \widetilde{C} \int_{0}^{\infty} \!\! \int_{v}^{\infty}
(q-v)^{\left( \beta+3 \right)/2} v^{\delta +\alpha} e^{-q/\vartheta_\text{int}}
e^{-\eta v} \text{d}q\, \text{d}v.
\end{align*}
Changing the integration variables from $(q, v)$ to $(\overline{q}, \overline{v})$,
where
\begin{align*}
\overline{q} = \frac{q - v}{\vartheta_\text{int}}, \qquad
\overline{v} = \frac{v}{\vartheta_\text{tr}}, \qquad
\vartheta_\text{tr} = \frac{k_\text{B} T_\text{int}}{\zeta}
= \frac{|T_\text{tr} - T_\text{int}|}{T_\text{tr}},
\end{align*}
with $\vartheta_\text{int}$ being defined in relations~\eqref{new-var}, and
noting that $\vartheta_\text{tr}/\vartheta_\text{int}=T_\text{int}/T_\text{tr}$
and $\eta \vartheta_\text{tr} = 1 - T_\text{int}/T_\text{tr}$, we
have
\begin{align}\label{Omega-q-v-bar}
\Omega&  = \eta\, \widetilde{C} \vartheta_\text{tr}^{\delta+\alpha+1}
\vartheta_\text{int}^{(\beta+5)/2}
\int_0^\infty \overline{q}^{(\beta+3)/2} e^{-\overline{q}} \text{d}\overline{q} \cdot
\int_0^\infty \overline{v}^{\delta+\alpha} e^{-\overline{v}} \text{d}\overline{v}
\nonumber \\
& = \eta\, \widetilde{C} \vartheta_\text{tr}^{\delta+\alpha+1}
\vartheta_\text{int}^{(\beta+5)/2}
\Gamma \big( (\beta+5)/2 \big) \Gamma (\delta+\alpha+1).
\end{align}

By substituting this $\Omega$ into expression~\eqref{source-term} and using
the explicit forms of $\vartheta_\text{tr}$ and $\vartheta_\text{int}$,
the following expression of $\big(I,\, Q_\text{s} (M_\text{r}, M_\text{r}) \big)$
is obtained:
\begin{align}\label{source-term-2}
\big( I, Q_\text{s} (M_\text{r}, M_\text{r}) \big) = \mathcal{F} (\rho, T_\text{tr}, T_\text{int})\,
(T_\text{tr} - T_\text{int}),
\end{align}
where
\begin{subequations}\label{source-term-aux}
\begin{align}
& \mathcal{F} (\rho, T_\text{tr}, T_\text{int})
= C\, \frac{k_\text{B}^{\alpha + 1 +\beta/2}}{m^{2 + \beta/2}}\,
\rho^2 T_\text{tr}^{\beta/2} T_\text{int}^\alpha,
\\
& C = 2^{\beta +2} \sqrt{\pi }\,
\frac{\Gamma \left( \delta +\alpha +1 \right)
\Gamma^{2}\left( \delta /2 \right) \Gamma \big( \left( \beta+5\right) /2\big)}
{\left[ \delta +\alpha +\left( \beta +3 \right)/2 \right]
\Gamma^{2} \left( \delta \right) }\, C_\text{r}.
\end{align}
\end{subequations}
%

\section{Positivity of $\Lambda_\mu$, $\Lambda_\text{tr}^\text{tr}$,
and $\Lambda_\text{int}^\text{int}$}\label{sec:positivity}

Let us first recall that (Proposition \ref{prop:positivity-kernel})
\begin{align}\label{L-positivity}
(\mathcal{L}_\text{r} h,\, M_\text{r} h) > 0,
\end{align}
for $h$ in $( M_\text{r} \text{ker} \mathcal{L}_\text{r})^\perp$.

Next, let us put $h = A_{ij} (\bm{c}) \mathcal{A} (|\bm{c}|, I)$,
which is in $( M_\text{r} \text{ker} \mathcal{L}_\text{r})^\perp$. Then, using
the first of equations~\eqref{EQ-decomp}, we have
\begin{align*}
0 & < \Big( \mathcal{L}_\text{r} \big( A_{ij} (\bm{c}) \mathcal{A}(|\bm{c}|, I) \big),
\, M_\text{r} A_{ij} (\bm{c}) \mathcal{A}(|\bm{c}|, I) \Big)
\\
& = \big( A_{ij} (\bm{c}), \, M_\text{r} A_{ij} (\bm{c}) \mathcal{A} (|\bm{c}|, I) \big)
\\
& = \frac{2}{3} \int_{\mathbb{R}^3 \times \mathbb{R}_+} |\bm{c}|^4 M_\text{r} \mathcal{A} (|\bm{c}|, I) \, \text{d}\bm{\xi} \text{d}I,
\end{align*}
so that
\begin{align*}
0 & < \int_0^\infty
\left[ \int_{\mathbb{R}^3} |\bm{c}|^4 \mathcal{A} (|\bm{c}|, I) \exp \left( - \frac{m |\bm{c}|^2}{2k_\mathrm{B} T_\mathrm{tr}} \right)
\text{d} \bm{c} \right]
I^{\delta/2-1} \exp \left( - \frac{I}{k_\mathrm{B} T_\mathrm{int}} \right)
\text{d} I
\\
& = 4\pi \int_0^\infty
\left[ \int_0^\infty c^6  \mathcal{A} (c, I) \exp \left( - \frac{m c^2}{2k_\mathrm{B} T_\mathrm{tr}} \right)
\text{d} c \right]
I^{\delta/2-1} \exp \left( - \frac{I}{k_\mathrm{B} T_\mathrm{int}} \right)
\text{d} I.
\end{align*}
This means by equalities~\eqref{Lambda-1} that $\Lambda_\mu (\rho, T_\text{tr}, T_\text{int}) > 0$.

Next, we let $h = c_i \mathcal{B} (|\bm{c}|, I)$. It belongs to $( M_\text{r} \text{ker} \mathcal{L}_\text{r})^\perp$
due to identities~\eqref{sub-cond}. Then, by the use of the second equation of~\eqref{EQ-decomp}, it follows from
the bound~\eqref{L-positivity} that
\begin{align*}
0 & <  \Big( \mathcal{L}_\text{r} \big( c_i \mathcal{B} (|\bm{c}|, I) \big), \, M_\text{r} c_i \mathcal{B} (|\bm{c}|, I) \Big)
\\
& =  \big( B_i (\bm{c}), \, M_\text{r} c_i \mathcal{B} (|\bm{c}|, I) \big)
\\
& = \int_{\mathbb{R}^3 \times \mathbb{R}_+} |\bm{c}|^2
\Big( \frac{m |\bm{c}|^2}{2k_\text{B} T_\text{tr}} - \frac{5}{2} \Big) M_\text{r} \mathcal{B} (|\bm{c}|, I)
\text{d}\bm{\xi} \text{d}I.
\end{align*}
Thus, taking account of identities~\eqref{sub-cond}, we obtain
\begin{align*}
0 & < \int_0^\infty
\left[ \int_{\mathbb{R}^3} |\bm{c}|^2 \Big( \frac{m |\bm{c}|^2}{2k_\text{B} T_\text{tr}} - \frac{5}{2} \Big)
\mathcal{B} (|\bm{c}|, I) \exp \!  \left( - \frac{m |\bm{c}|^2}{2k_\mathrm{B} T_\mathrm{tr}} \right)
\text{d} \bm{c} \right]\!
I^{\delta/2-1} \exp \! \left( - \frac{I}{k_\mathrm{B} T_\mathrm{int}} \right)
\text{d} I
\\
& = \frac{4\pi m}{2k_\mathrm{B} T_\mathrm{tr}} \int_0^\infty
\left[ \int_0^\infty c^6  \mathcal{B} (c, I) \exp \left( - \frac{m c^2}{2k_\mathrm{B} T_\mathrm{tr}} \right)
\text{d} c \right] 
I^{\delta/2-1} \exp \left( - \frac{I}{k_\mathrm{B} T_\mathrm{int}} \right)
\text{d} I.
\end{align*}
This shows from equalities~\eqref{Lambda-2-a} that
$\Lambda_\text{tr}^\text{tr} (\rho, T_\text{tr}, T_\text{int}) > 0$.
 Letting $h = c_i\, \mathcal{C} (|\bm{c}|, I)$ and making a similar argument,
 one can readily show that $\Lambda_\text{int}^\text{int} (\rho, T_\text{tr}, T_\text{int}) > 0$.

\section{Derivation of expressions~\eqref{Lr-h-alpha=0-2} and \eqref{Lr-h-alpha=0-3}}\label{sec:derivation}

If $h$ is a function of $\bm{c}$ only and does not depend on $I$, then equation~\eqref{Lr-h-alpha=0}
can be transformed as
\begin{align*}
\mathcal{L}_\text{r} h & = - C_\text{r} \frac{\rho/m}{(2\pi k_\text{B} T_\text{tr}/m)^{3/2}}
\int_{\mathbb{R}^3 \times \mathbb{S}^2} \exp \left( - \frac{m|\bm{c}_*|^2}{2k_\text{B} T_\text{tr}} \right)
( h_*' + h' - h_* - h ) |\bm{g}|^\beta \text{d} \bm{c}_* \text{d} \bm{\sigma}
\\
& \qquad \times \int_0^1 [r(1-r)]^{\delta/2-1} \text{d}r
\cdot \frac{1}{(k_\text{B} T_\text{int})^{\delta/2} \Gamma (\delta/2)}
\int_0^\infty I_*^{\delta/2-1} \exp \left( - \frac{I_*}{k_\text{B} T_\text{int}} \right) \text{d} I_*.
\end{align*}
Since $\int_0^1 [r(1-r)]^{\delta/2-1} \text{d}r = B(\delta/2, \delta/2) = \Gamma^2(\delta/2)/\Gamma (\delta)$,
where $B(x, y)$ is the beta function, and
$\int_0^\infty I_*^{\delta/2-1} \exp \left( - I_*/(k_\text{B} T_\text{int} ) \right) \text{d} I_*
= (k_\text{B} T_\text{int})^{\delta/2} \Gamma (\delta/2)$,~\eqref{Lr-h-alpha=0-2} follows.

Next, let us assume that $h$ is of the form $h = \left[ I/(k_\text{B} T_\text{int}) - \delta/2 \right]\widetilde{h}(\bm{c})$,
with $\widetilde{h}(\bm{c})$ being independent of $I$. Then, equation~\eqref{Lr-h-alpha=0}
can be written as
\begin{align*}
\mathcal{L}_\text{r} h & = - C_\text{r}
\frac{\sqrt{m} \rho}{(2\pi k_\text{B} T_\text{tr})^{3/2} (k_\text{B} T_\text{int})^{\delta/2} \Gamma (\delta/2)}
\\
& \quad \;\; \times \int_{[0, 1] \times \mathbb{S}^2 \times \mathbb{R}^3 \times \mathbb{R}_+}
\exp \left( - \frac{m|\bm{c}_*|^2}{2k_\text{B} T_\text{tr}} \right) \exp \left( - \frac{I_*}{k_\text{B} T_\text{int}} \right)
[r(1-r)]^{\delta/2-1} I_*^{\delta/2-1} |\bm{g}|^\beta
\\
& \qquad  \quad \times
\left[ \frac{1}{k_\text{B} T_\text{int}} ( I_*' \widetilde{h}_*' + I' \widetilde{h}' - I_*\widetilde{h}_* -I\widetilde{h} )
- \frac{\delta}{2} ( \widetilde{h}_*' + \widetilde{h}' - \widetilde{h}_* - \widetilde{h} ) \right]
\text{d}r \text{d}\bm{\sigma} \text{d}\bm{c}_* \text{d}I_*.
\end{align*}
One can replace $\widetilde{h}_*'$ with $\widetilde{h}'$ in the above equation because
$\bm{c}_*'$ becomes $\bm{c}'$ by $\bm{\sigma} \to -\bm{\sigma}$ [cf.~relations~\eqref{new-var}].
In addition, the integral
$\int_0^\infty  I_*^{\delta/2-1} \left[ I_*/(k_\text{B} T_\text{int}) - \delta/2 \right]
\exp \left( - I_*/(k_\text{B} T_\text{int}) \right) \text{d}I_*$ vanishes. Therefore, using
the relation $I'+I_*'=I+I_*$, we have
\begin{align*}
\mathcal{L}_\text{r} h & = - C_\text{r}
\frac{\sqrt{m} \rho}{(2\pi k_\text{B} T_\text{tr})^{3/2} (k_\text{B} T_\text{int})^{\delta/2} \Gamma (\delta/2)}
\\
& \quad \;\; \times \left\{
\int_{\mathbb{R}^3 \times \mathbb{S}^2} \exp \left( - \frac{m|\bm{c}_*|^2}{2k_\text{B} T_\text{tr}} \right)
|\bm{g}|^\beta  \left[ \left( \frac{I}{k_\text{B} T_\text{int}} - \frac{\delta}{2} \right) \widetilde{h}
+ \delta \widetilde{h}' \right] \text{d}\bm{c}_* \text{d}\bm{\sigma}
\right.
\\
& \qquad \qquad \qquad \times
\int_0^1 [r(1-r)]^{\delta/2-1} \text{d}r \cdot \int_0^\infty I_*^{\delta/2-1}
\exp \left( - \frac{I_*}{k_\text{B} T_\text{int}} \right) \text{d}I_*
\\
& \qquad \quad
- \int_{\mathbb{R}^3 \times \mathbb{S}^2}
\exp \left( - \frac{m|\bm{c}_*|^2}{2k_\text{B} T_\text{tr}} \right) |\bm{g}|^\beta
\, \widetilde{h}' 
\text{d} \bm{c}_* \text{d} \bm{\sigma} 
\\
& \qquad \qquad \qquad \left. \times 
\int_0^1 [r(1-r)]^{\delta/2-1} \text{d}r \cdot
\int_0^\infty \frac{I+I_*}{k_\text{B} T_\text{int}}  I_*^{\delta/2-1}
\exp \left( - \frac{I_*}{k_\text{B} T_\text{int}} \right) \text{d}I_*
\right\}.
\end{align*}
Expressing the integral with respect to $r$ and that with respect to $I_*$ in terms
of the gamma functions as was done above and using basic properties
of the gamma function, one obtains~\eqref{Lr-h-alpha=0-3}.

\section{Spherical symmetry of $Q_\text{s} (M_\text{r}, M_\text{r})$}\label{sec:symmetry}

Substituting equality~\eqref{sigma-model-a} into expression~\eqref{Qs-BL}, using
identity~\eqref{relation}, and taking account of relations~\eqref{variables} in $\bm{c}$ variables, i.e.,
\begin{align}\label{variables-c}
\begin{aligned}
& \bm{g} = \bm{c} - \bm{c}_*, \qquad \bm{g}' = \bm{c}' - \bm{c}_*'
\\
& \bm{c}' = \frac{\bm{c} + \bm{c}_*}{2} + \sqrt{\frac{RE}{m}}\, \bm{\sigma}, \qquad
\bm{c}_*' = \frac{\bm{c} + \bm{c}_*}{2} - \sqrt{\frac{RE}{m}}\, \bm{\sigma},
\\
& I' = r (1 - R) E, \qquad I_*' = (1 - r)(1 - R)E,
\\
& E = \frac{m}{4} |\bm{c} - \bm{c}_*|^2 + I + I_*,
\end{aligned}
\end{align}

\noindent
we obtain the following expression of $Q_\text{s} (M_\text{r}, M_\text{r})$:
\begin{align}\label{Qs-Mr-Mr-3}
& Q_\text{s}(M_\text{r}, M_\text{r}) (\bm{c}, I)
\nonumber \\
& \;\;\;
=C_\text{s} \int_{[0,1]^{2} \times \mathbb{S}^{2} \times \mathbb{R}^{3}\times \mathbb{R}_{+}}
M_\text{s}^\text{tr} M_{\text{s} \ast}^\text{tr} \frac{T_\text{tr}^{\delta}}{T_\text{int}^{\delta}}
\left\vert \bm{c} - \bm{c}_{\ast}\right\vert^{\beta}
\left( I+I_{\ast}\right)^{\alpha} \left( e^{-\eta \left( 1-R \right) E/\zeta}
-e^{-\eta \left( I+I_{\ast} \right) /\zeta} \right) 
\nonumber \\
& \qquad \qquad \qquad \qquad \qquad
\times \left[ r\left( 1-r\right) \right] ^{\delta /2-1}R^{\left( \beta+1\right) /2}
\left( 1-R\right) ^{\alpha +\delta -1}
\text{d}R\, \text{d}r\, \text{d} \bm{\sigma}\, \text{d}\bm{\xi }_{\ast}\, \text{d}I_{\ast}
\nonumber \\
& \;\;\;
= \frac{m\rho^{2}C_\text{s}}{2\pi^2 \left( k_\text{B}T_\text{tr}\right)^{3}
\left( k_\text{B}T_\text{int}\right)^{\delta }\Gamma \left( \delta \right) }
\nonumber \\
& \quad \times
\int_{[0,1]\times \mathbb{R}^{3}\times \mathbb{R}_{+}}
\left\vert \bm{c}-\bm{c}_{\ast } \right\vert ^{\beta }
e^{-m\left( \left\vert \bm{c} \right\vert^{2}
+\left\vert \bm{c}_{\ast } \right\vert^{2} \right) /
\left( 2k_{B}T_\text{tr} \right) }
\left( I+I_{\ast} \right)^{\alpha} \left( II_{\ast} \right)^{\delta /2-1}e^{-\left( I+I_{\ast }  \right) /
\left( k_{B} T_\text{tr} \right) }
\nonumber \\
& \qquad \qquad \qquad \qquad \times 
\left( e^{-\eta \left( 1-R\right) E/\zeta }-e^{-\eta \left(I+I_{\ast }\right) /\zeta }\right)
R^{\left( \beta +1\right) /2}\left(1-R\right) ^{\alpha +\delta -1}
\text{d}R\, \text{d}\bm{c}_{\ast}\, \text{d}I_{\ast},
\end{align}
where the relation $\int_0^1 [r(1-r)]^{\delta/2-1} \text{d}r = \Gamma^2 (\delta/2)/\Gamma (\delta)$
has been used (cf.~Appendix \ref{sec:source}).
Thus, $Q_\text{s} (M_\text{r},M_\text{r}) (S\bm{c},I)=Q_\text{s} (M_\text{r},M_\text{r}) (\bm{c}, I)$
for any isometry $S\in O(3)$, and, hence, $Q_\text{s} (M_\text{r},M_\text{r})$ is a function
of $\left\vert \bm{c} \right\vert $ and $I$.

\section{Derivation of expression~\eqref{moment-h_1}}\label{sec:derivation-2}

Let $\text{d}A_\text{s}$ denote the measure $\text{d}A_\theta$ with $\theta=1$ [see~\eqref{dAtheta}]. From Proposition~\ref{P1} and Lemma~\ref{L0},
we have
\begin{align*}
\big( I,\, Q_\text{s} (M_\text{r}, M_\text{r} h_1) \big)
& = \frac{1}{8} \int_{ (\mathbb{R}^3 \times \mathbb{R}_+)^4 }
\left[ \frac{M_\text{r}' M_{\text{r}*}' (h_1' + h_{1*}')}{(I' I_*')^{\delta/2-1}}
- \frac{M_\text{r} M_{\text{r}*} (h_1 + h_{1*})}{(I I_*)^{\delta/2-1}} \right] 
\nonumber \\
& \qquad \qquad \qquad \qquad \qquad \qquad \times
(I + I_* - I' - I_*') \text{d} A_\text{s}
\nonumber \\
& = \frac{1}{4} \int_{ (\mathbb{R}^3 \times \mathbb{R}_+)^4 }
\frac{M_\text{r} M_{\text{r}*} (h_1 + h_{1*})}{(I I_*)^{\delta/2-1}}
(I' + I_*' - I - I_*) \text{d} A_\text{s}
\nonumber \\
& = \frac{1}{2} \int_{ (\mathbb{R}^3 \times \mathbb{R}_+)^4 }
\frac{M_\text{r} M_{\text{r}*} }{(I I_*)^{\delta/2-1}} h_1
(I' + I_*' - I - I_*) \text{d} A_\text{s}.
\end{align*}
If we change the integration variables as in~\eqref{Qs-BL}, noting that there are additional integrations with respect to $\bm{\xi}$ and $I$ here,
and change the integration variables $(\bm{\xi}, \bm{\xi_*})$ to
$(\bm{c}, \bm{c}_*)$, where $\bm{c}=\bm{\xi}-\bm{u}$ and $\bm{c}_*=\bm{\xi}_*-\bm{u}$,
then we obtain
\begin{align}\label{I-Qs-Mr-Mrh1}
\big( I,\, Q_\text{s} (M_\text{r}, M_\text{r} h_1) \big)
& = \frac{1}{2} \int_{ [0, 1]^2 \times \mathbb{S}^2 \times (\mathbb{R}^3 \times \mathbb{R}_+)^2 }
\frac{M_\text{r} M_{\text{r}*} }{(I I_*)^{\delta/2-1}} h_1 (\bm{c}, I)
(I' + I_*' - I - I_*)
\nonumber \\
& \qquad \qquad \qquad \times
\left( I I_{\ast } \right)^{\delta /2-1}
|\bm{g}| \sigma_\text{s} (1-R) E^2
\text{d}R\, \text{d}r\, \text{d}\bm{\sigma}\, \text{d}\bm{c}_* \, \text{d}I_*
\text{d}\bm{c}\, \text{d}I.
\end{align}
where the arguments $t$ and $\bm{x}$ in $h_1$ are omitted. On the other hand,
by the use of the relations \eqref{variables-c}, the following expressions of
$\sigma_\text{s}$ [see~\eqref{sigma-model-a}] and $I' + I_*' - I - I_*$ are
obtained: 
\begin{align*}
& \sigma_\text{s} = C_\text{s} (I + I_*)^\alpha [r(1 - r)]^{\delta/2-1} (1 - R)^{\alpha + \delta -2}
R^{(\beta+1)/2} E^{-2} |\bm{g}|^{\beta-1},
\\
& I' + I_*' - I - I_* = \frac{m}{4} (1 - R) |\bm{c} - \bm{c}_*|^2 - R(I + I_*).
\end{align*}
If we substitute these results, as well as the explicit forms of $M_\text{r}$ and $M_{\text{r}*}$,
into equality~\eqref{I-Qs-Mr-Mrh1} and carry out the integrations
with respect to $r$ and $\bm{\sigma}$, we obtain~\eqref{moment-h_1}.

\section{Proof of inequality \eqref{inequality}}\label{sec:proof-ineq}

By replacing $h_1$ with $\widetilde{D}$ in~\eqref{moment-h_1} and expressing it
using $M_\text{r}$ and $M_{\text{r}*}$ instead of
$e^{-(|\bm{c}|^2 + |\bm{c}_*|^2)/(2k_\text{B} T_\text{tr})}$$\times$
$e^{-(I+I_*)/(k_\text{B}T_\text{int})}$$\times$$(II_*)^{\delta/2-1}$,
one obtains
\begin{align*}
& \left( I, Q_\text{s} \big( M_\text{r}, M_\text{r} \widetilde{D} \big) \right)^{2}
\\
& \quad =\frac{4\pi^{2} C_\text{s}^{2} \Gamma^{4}\left( \delta /2 \right) }
{\Gamma^{2}\left( \delta \right)}
\bigg\{ \int_{[0,1] \times \left( \mathbb{R}^{3} \times \mathbb{R}_{+} \right)^{2}}
\widetilde{D} |\bm{c} - \bm{c}_*|^{\beta } M_\text{r} M_{\text{r}*} \left( I+I_* \right)^{\alpha}
\\
& \qquad \;\; \times \left[ \frac{m}{4} (1-R) |\bm{c} - \bm{c}_*|^2 - R ( I+I_{\ast } ) \right]
R^{\left( \beta+1\right)/2} \left( 1-R \right)^{\alpha +\delta -1}
\text{d}R \text{d} \bm{c}_* \text{d} I_{\ast} \text{d} \bm{c}\, \text{d} I \bigg\}^{2}.
\end{align*}
Then, with the help of the Cauchy-Schwarz inequality, the following inequality
is obtained:
\begin{align*}
& \left( I, Q_\text{s} \big( M_\text{r}, M_\text{r} \widetilde{D} \big) \right)^{2}
\leq \frac{4\pi^{2} C_\text{s}^{2} \Gamma^{4} \left( \delta /2 \right) }
{\Gamma^{2} \left( \delta \right) }\,
S_1 \times S_2,
\end{align*}
where
\begin{align*}
& S_1 = \int_{[0,1] \times \left( \mathbb{R}^{3} \times \mathbb{R}_{+} \right) ^{2}}
\widetilde{D}^{2} M_\text{r} M_{\text{r}\ast}
\text{d}R \text{d} \bm{c}_* \text{d} I_{\ast} \text{d} \bm{c}\, \text{d} I
\\
& S_2 = \int_{[0,1] \times \left( \mathbb{R}^{3} \times \mathbb{R}_{+} \right)^{2}}
|\bm{c} - \bm{c}_*|^{2\beta} M_\text{r }M_{\text{r}\ast} \left( I+I_{\ast }\right)^{2\alpha}
\left[ \frac{m}{4} (1-R) |\bm{c} - \bm{c}_*|^{2} - R ( I+I_{\ast} ) \right]^{2}
\\
& \qquad \qquad \qquad \qquad \qquad \qquad \times
R^{\beta +1} \left( 1-R \right)^{2\left( \alpha +\delta -1 \right) }
\text{d}R \text{d} \bm{c}_* \text{d} I_{\ast} \text{d} \bm{c}\, \text{d} I.
\end{align*}
Using estimate~\eqref{DMD} and the H\"older inequality, the factors
$S_1$ and $S_2$ are estimated as follows:
\begin{align*}
& S_1 = \int_0^1 \text{d}R \cdot
\int_{ \mathbb{R}^{3}\times \mathbb{R}_{+}} M_{\text{r}\ast } \text{d}\bm{c}_* \text{d}I_* \cdot
\int_{ \mathbb{R}^{3}\times \mathbb{R}_{+}}
\widetilde{D}^{2} M_\text{r} \text{d}\bm{c}\, \text{d}I
= \frac{\rho}{m} \left( \widetilde{D}, M_\text{r} \widetilde{D} \right)
\\
& \quad \leq \frac{\rho}{m \mu} \left( D, M_\text{r} D \right),
\\
& S_2 \leq \int_{0}^{1} \text{d} R \cdot
\int_{\left( \mathbb{R}^{3}\times \mathbb{R}_{+}\right)^{2}}
M_\text{r} M_{\text{r}\ast }\, 4^{\alpha +\beta }
\left( |\bm{c}|^{2\beta} + |\bm{c}_*|^{2\beta } \right)
\left( I^{2\alpha } + I_{\ast }^{2\alpha } \right)
\\
& \qquad \qquad \qquad \qquad \qquad \qquad \times
\left[ m \left( |\bm{c}|^{2} + |\bm{c}_*|^{2} \right)
+\left( I+I_{\ast } \right) \right]^{2}
\text{d}\bm{c}_{\ast } \text{d}I_{\ast } \text{d}\bm{c}\, \text{d}I < +\infty
\end{align*}
Thus, inequality~\eqref{inequality} follows.

\section{Proof of inequalities \eqref{estimate-Qs}}\label{sec:proof-estimate}

Let us first recall that $1/\zeta=(1/k_\text{B})(1/T_\text{int} - 1/T_\text{tr})$
if $T_\text{tr} > T_\text{int}$ (i.e., $\eta=1$) and that
$1/\zeta=(1/k_\text{B})(1/T_\text{tr} - 1/T_\text{int})$ if $T_\text{tr} < T_\text{int}$
(i.e., $\eta=-1$) [cf.~\eqref{zeta-eta}]. Then expression~\eqref{Qs-Mr-Mr-3} and estimates \eqref{estimate:e-e}
and \eqref{estimate:E} lead to the following inequalities for
$M_\text{r}^{-1}$$Q_\text{s} (M_\text{r}, M_\text{r})$
(note that $C_\text{g}$ is a generic positive constant depending on the macroscopic quantities):
\begin{itemize}

\item
For $T_\text{tr} > T_\text{int}$:
\begin{align*}
&M_\text{r}^{-1} \left\vert Q_\text{s} (M_\text{r}, M_\text{r}) \right\vert
\notag \\
& \quad \leq C_\text{g} \frac{M_\text{s}^\text{tr}}{M_\text{r}}
\int_{[0,1] \times \mathbb{R}^{3} \times \mathbb{R}_{+}}
\left( \left\vert \bm{c} \right\vert^{\beta} + \left\vert \bm{c}_{\ast } \right\vert^{\beta} \right)
e^{-m \left\vert \bm{c}_{\ast} \right\vert^{2}/\left( 2k_\text{B}T_\text{tr} \right) }
\notag \\
& \qquad \qquad \qquad \;\; \times
\left( I^{\alpha } + I_{\ast}^{\alpha} \right) I_{\ast }^{\delta /2-1}
e^{-I_{\ast }/\left(k_\text{B}T_\text{tr}\right) }
\left\vert e^{-(1-R)E/\zeta } - e^{-\left( I+I_{\ast } \right) /\zeta} \right\vert
\text{d}R \text{d}\bm{c}_{\ast} \text{d}I_{\ast}
\notag \\
& \quad \leq C_\text{g} \zeta ^{-1} e^{I/\zeta } \int_{0}^{1} \text{d}R \cdot
\int_{\mathbb{R}^{3} \times \mathbb{R}_{+}}
\left( \left\vert \bm{c} \right\vert^{\beta } + \left\vert \bm{c}_{\ast } \right\vert^{\beta } \right)
e^{-m \left\vert \bm{c}_{\ast} \right\vert^{2}/\left( 2k_\text{B} T_\text{tr} \right) }
\notag \\
& \qquad \qquad \qquad \qquad\;\;\;\; \times
\left( I^{\alpha }+I_{\ast}^{\alpha } \right) I_{\ast }^{\delta /2-1}
e^{-I_{\ast }/\left(k_\text{B} T_\text{tr} \right) }
\left( \left\vert \bm{c} \right\vert^{2} + \left\vert \bm{c}_{\ast } \right\vert^{2} + I + I_{\ast } \right)
\text{d}\bm{c}_{\ast } \text{d}I_{\ast}
\notag \\
& \quad \leq C_\text{g} \left\vert T_\text{tr} - T_\text{int} \right\vert
\left( 1 + \left\vert \bm{c} \right\vert^{2} + \left\vert \bm{c} \right\vert^{\beta +2} \right)
\left( 1 + I + I^{\alpha +1} \right) e^{I/\zeta }.
\end{align*}

\item
For $T_\text{tr} < T_\text{int}$:
\begin{align*}
&M_\text{r}^{-1} \left\vert Q_\text{s} (M_\text{r}, M_\text{r}) \right\vert
\notag \\
& \quad \leq C_\text{g} \frac{M_\text{s}^\text{tr}}{M_\text{r}}
\int_{[0,1] \times \mathbb{R}^{3} \times \mathbb{R}_{+}}
\left( \left\vert \bm{c} \right\vert^{\beta} + \left\vert \bm{c}_{\ast } \right\vert^{\beta} \right)
e^{-m \left\vert \bm{c}_{\ast} \right\vert^{2}/\left( 2k_\text{B}T_\text{tr} \right) }
\notag \\
& \qquad \qquad \times
\left( I^{\alpha } + I_{\ast}^{\alpha} \right) I_{\ast }^{\delta /2-1}
e^{-I_{\ast }/\left(k_\text{B}T_\text{tr}\right) }
e^{E/\zeta}  \left\vert e^{-RE/\zeta} - e^{-m|\bm{c}-\bm{c}_*|^2/(4\zeta)} \right\vert
\text{d}R \text{d}\bm{c}_{\ast} \text{d}I_{\ast}
\notag \\
& \quad \leq C_\text{g} \zeta ^{-1} e^{m|\bm{c}|^2/(2\zeta) } \int_{0}^{1} \text{d}R \cdot
\int_{\mathbb{R}^{3} \times \mathbb{R}_{+}}
\left( \left\vert \bm{c} \right\vert^{\beta } + \left\vert \bm{c}_{\ast } \right\vert^{\beta } \right)
e^{-m \left\vert \bm{c}_{\ast} \right\vert^{2}/\left( 2k_\text{B} T_\text{int} \right) }
\notag \\
& \qquad \qquad \qquad \qquad \;\;\; \times
\left( I^{\alpha }+I_{\ast}^{\alpha } \right) I_{\ast }^{\delta /2-1}
e^{-I_{\ast }/\left(k_\text{B} T_\text{int} \right) }
\left( \left\vert \bm{c} \right\vert^{2} + \left\vert \bm{c}_{\ast } \right\vert^{2} + I + I_{\ast } \right)
\text{d}\bm{c}_{\ast } \text{d}I_{\ast}
\notag \\
& \quad \leq C_\text{g} \left\vert T_\text{tr} - T_\text{int} \right\vert
\left( 1 + \left\vert \bm{c} \right\vert^{2} + \left\vert \bm{c} \right\vert^{\beta +2} \right)
\left( 1 + I + I^{\alpha +1} \right) e^{m|\bm{c}|^2/(2\zeta)}.
\end{align*}

\end{itemize}

\end{document}